\documentstyle[amssymb,12pt]{amsart}
\newtheorem{prop}{Proposition}[section]
\newtheorem{prop:def}{Proposition-Definition}[section]

\newtheorem{lemma}{Lemma}[section]
\newtheorem{thm}{Theorem}[section]
\newtheorem{cor}{Corollary}[section]

\theoremstyle{remark}
\newtheorem{remark}{Remark}
\newtheorem{notation}{Notation}

\begin{document}
\newcommand{\nc}{\newcommand} \nc{\on}{\operatorname}
\nc{\pa}{\partial} \nc{\cA}{{\cal A}}\nc{\cV}{{\cal V}} \nc{\cB}{{\cal
    B}}\nc{\cC}{{\cal C}} \nc{\cE}{{\cal E}}\nc{\cG}{{\cal
    G}}\nc{\cH}{{\cal H}} \nc{\cX}{{\cal X}}\nc{\cR}{{\cal
    R}}\nc{\cL}{{\cal L}} \nc{\cK}{{\cal K}}
\nc{\sh}{\on{sh}}\nc{\Id}{\on{Id}}\nc{\Diff}{\on{Diff}}
\nc{\ad}{\on{ad}}\nc{\Der}{\on{Der}}\nc{\End}{\on{End}}
\nc{\res}{\on{res}}\nc{\ddiv}{\on{div}}
\nc{\card}{\on{card}}\nc{\dimm}{\on{dim}}
\nc{\Jac}{\on{Jac}}\nc{\Ker}{\on{Ker}}\nc{\Vect}{\on{Vect}}
\nc{\Hom}{\on{Hom}}
\nc{\Spec}{\on{Spec}}\nc{\Cl}{\on{Cl}}
\nc{\Imm}{\on{Im}}\nc{\limm}{\on{lim}}\nc{\Ad}{\on{Ad}}
\nc{\ev}{\on{ev}} \nc{\corr}{\on{corr}} 
\nc{\Hol}{\on{Hol}}\nc{\Det}{\on{Det}}
\nc{\Bun}{\on{Bun}}\nc{\diag}{\on{diag}}
\nc{\de}{\delta}
\nc{\si}{\sigma}\nc{\ve}{\varepsilon}\nc{\z}{\zeta}
\nc{\al}{\alpha}\nc{\vp}{\varphi} \nc{\vpi}{\varpi} 
\nc{\CC}{{\mathbb
    C}}\nc{\ZZ}{{\mathbb Z}} \nc{\NN}{{\mathbb N}}\nc{\VV}{{\mathbb
    V}} \nc{\zz}{{\mathbf z}} \nc{\WW}{{\mathbb W}} 
\nc{\mmm}{{\mathbf m}} 
\nc{\ddelta}{{\underline{\delta}}} 
\nc{\AAA}{{\mathbb A}}\nc{\cO}{{\cal O}}
\nc{\cF}{{\cal F}}\nc{\cM}{{\cal M}}\nc{\cT}{{\cal T}}\nc{\cW}{{\cal W}}
\nc{\la}{{\lambda}}\nc{\G}{{\mathfrak g}}\nc{\mm}{{\mathfrak m}}
\nc{\A}{{\mathfrak a}} \nc{\HH}{{\mathfrak h}} \nc{\N}{{\mathfrak
    n}}\nc{\B}{{\mathfrak b}} \nc{\La}{\Lambda}
\nc{\g}{\gamma}\nc{\eps}{\epsilon}\nc{\wt}{\widetilde}
\nc{\wh}{\widehat} \nc{\bn}{\begin{equation}}\nc{\en}{\end{equation}}
\nc{\SL}{{\mathfrak{sl}}}\nc{\ttt}{{\mathfrak{t}}} \nc{\td}{\tilde}

%
%
%

\newcommand{\ldar}[1]{\begin{picture}(10,50)(-5,-25)
\put(0,25){\vector(0,-1){50}}
\put(5,0){\mbox{$#1$}} 
\end{picture}}

\newcommand{\lrar}[1]{\begin{picture}(50,10)(-25,-5)
\put(-25,0){\vector(1,0){50}}
\put(0,5){\makebox(0,0)[b]{\mbox{$#1$}}}
\end{picture}}

\newcommand{\luar}[1]{\begin{picture}(10,50)(-5,-25)
\put(0,-25){\vector(0,1){50}}
\put(5,0){\mbox{$#1$}}
\end{picture}}

\title[Solutions of the KZB equations in genus greater than one] {Solutions
  of the KZB equations in genus greater than one}

\author{B. Enriquez}

\address{B.E.: DMA, CNRS, ENS, 45 rue d'Ulm, 75005 Paris, 
France}

\author{G. Felder}

\address{G.F.: D-Math, ETH-Zentrum, HG G44, CH-8092 Zurich,
 Switzerland}

\date{January 2000}

\maketitle

\subsection*{Introduction}

The Knizhnik-Zamolodchikov-Bernard connection (\cite{KZ,Bernard,TUY})
can be viewed  as a flat connection  over the complement
$\det_{g,\nu}^*$ of the zero-section in the total space of the
determinant line bundle over the moduli space $\cM_{g,\nu\cdot 1^2}$ of
genus $g$ curves with $\nu$ marked points and tangent vectors
(\cite{BFM}). Projectivization of this connection  is the pull-back of a
flat projective connection on $\cM_{g,\nu\cdot 1^2}$. 

The aim of this paper is two-fold. We first construct a flat version of
the KZB connection  over the moduli space $\cM_{g,\nu\cdot 1^2}$, when
$\nu\geq 1$ (we will set $\nu = n +1$). We then give integral  formulas
for flat sections of this connection, using the functional
parametrization of  conformal blocks and the KZB connection introduced
in \cite{EF}.

\subsubsection*{Construction of a flat connection} More precisely, we
consider the moduli space $\cM_{g,1^{\infty},n\cdot 1^2}$ of  systems $m
= (X,P_0,t,P_i,v_i)$ consisting of a curve $X$ of genus $g$, a point
$P_0$, a jet of coordinate $t$ at $P_0$, and $n$ points $P_i$ with 
tangent vectors $v_i$ ($i = 1, \ldots,n$). We associate to the  data of
a semisimple Lie algebra $\bar\G$, a representation  $\VV$ of the
Kac-Moody loop algebra $\G$ of $\bar\G$ and irreducible  representations
$V_i$ of $\bar\G$, the sheaf $\cB_{\VV,(V_i)}$ of coinvariants over
$\cM_{g,1^{\infty},n\cdot 1^2}$.  We denote by $\cW^\infty$ the subspace
of $\cM_{g,1^{\infty},n\cdot 1^2}$  formed of the systems $m$ such that
$P_0$ is not a Weierstrass point of $X$.  The main result of section
\ref{sect:1}  is the construction of a flat connection
$\nabla^{\VV,(V_i)}$ on this sheaf, which is smooth over the complement
of $\cW^{\infty}$   (Thm.\ \ref{thm:flatness}). $\nabla^{\VV,(V_i)}$ has
 a logarithmic singularity around $\cW^\infty$, which we compute in
sect.\ \ref{comp:log}. Moreover,  the projectivization of
$\nabla^{\VV,(V_i)}$ is smooth and  is isomorphic to the connection of
\cite{TUY}. 

The proof of these facts is contained in sect.\ \ref{sect:1}.  Our main
tool is a normalization of the Sugawara tensor (formula (\ref{def:TR})),
 which makes use of direct sum decompositions of spaces of the formal functions 
and formal one-forms at a non-Weierstrass point of $X$ (Lemma
\ref{lemma:decompositions}). 

We can therefore think that the curvature of $\nabla^{\VV,(V_i)}$ is a 
scalar delta-function concentrated on the Weierstrass locus
$\cW^\infty$.

\subsubsection*{Functional parametrization of $\nabla^{\VV,(V_i)}$}

Our next aim is to give explicit expressions for the flat sections of 
$\nabla^{\VV,(V_i)}$, in the special case where $\bar\G = \SL_2$.  For
this, we use the functional parametrization of this connection that we
introduced in our work \cite{EF}. In sect. \ref{sect:3}, we recall the
main results of \cite{EF} for the  $\SL_2$ case. 

Our work \cite{EF} was inspired by \cite{FS}, where conformal blocks are
parametrized by  correlation functions of generating currents of the Lie
subalgebra $\bar \N_+\otimes \CC((t))$ of $\G$, where $\bar \N_+$ is the
positive  nilpotent subalgebra of $\bar\G$.  In \cite{EF}, we  consider
twisted analogues of these correlation functions; this means that the
nilpotent currents are multiplied by exponentials of  elements of
$\bar\HH\otimes\CC((t))$, where $\bar\HH$ is the Cartan subalgebra of
 $\bar\G$. This enabled us in \cite{EF} to construct a ``functional'' sheaf 
$\cF_{\ell,(\La_i)}$ with flat connection $\nabla^{\cF_{\ell,(\La_i)}}$ 
over a covering 
$\cM^{(a)}_{g,1^{\infty},n\cdot 1^2}$, and a morphism of sheaves
with connections
$$
\corr : (\cB_{\VV,(V_i)},\nabla^{\VV,(V_i)}) \to (\cF_{\ell,(\La_i)}, 
\nabla^{\cF_{\ell,(\La_i)}}) . 
$$ 
We specialize $\VV$ to an induced representation 
$\VV_{-\ell k}$ and the $V_i$ to a lowest weight Verma module $V_{-\La_i}$
(see \ref{sect:modules}). 

The fiber of $\cF_{\ell,(\La_i)}$ over a point of the moduli space is a
space of functions $f(\la_1,\ldots,\la_g|z_1,\ldots,z_N)$,
defined on a subset of 
$\CC^g \times (\wt X - \sigma^{-1}(P_0))^N$, satisfying regularity, symmetry and 
transformation properties (see sect.\ \ref{corr}). Here $\sigma:
\wt X \to X$ is the universal cover map of $X$, and $N$ is given by 
$$
N = {1\over 2} [k\ell + \sum_{i}\La_i] .
$$ 
The connection $\nabla^{\cF_{\ell,(\La_i)}}$ is expressed by
(\ref{operator:KZ:moduli}), (\ref{KZB:moduli})  and
(\ref{expr:conn:cF}).   In the case $N=0$, the formula for
$\nabla^{\cF_{\ell,(\La_i)}}$ is  an analogue of the Lam\'e equation.

\subsubsection*{Connections on $\Omega_X$}

The expression of $\nabla^{\cF_{\ell,(\La_i)}}$ 
involves a connection $D^{(\la)}$ on the canonical bundle
$\Omega_X$ on $X$, depending on $\la \in\CC^g$, which is defined
by the formula
$$
(D^{(\la)}\omega)(z) = - \limm_{z'\to z}\left( G_\la(z,z')\omega(z') 
+ G_\la(z',z)\omega(z) \right) ,  
$$
where $G_\la(z,z')$ is the twisted Green function on $X$ 
corresponding to $\la$ (see (\ref{twisted:Green})). 

By a connection $\nabla_\Omega$ on $\Omega_X$, we understand the data
of a divisor $D$ on $X$ and a collection of compatible maps
$\Gamma(U,\Omega_X) \to \Gamma(U,\Omega_X^{\otimes 2}(D))$ for  each
open subset of $X$, satisfying the Leibniz rule. 
 
Any connection $\nabla_\Omega$ on $\Omega_X$ can be expressed in the form
$$
\nabla_\Omega\omega = \al_\nabla d(\omega / \al_\nabla),  
$$
where $\al_\nabla$ is some form on the universal cover of $X - D$. 

The aim of sect.\ \ref{sect:2} is to compute the form
$\al_\nabla$ associated with the connection $D^{(\la)}$. This is 
done in Prop.\ \ref{pestalozzi}.

\subsubsection*{Flat sections of  $\cF_{\ell,(\La_i)}$}

As we have seen, to describe flat connections of $\cB_{\VV,(V_i)}$ is 
the same as 

(i) describing  the flat sections of $\cF_{\ell,(\La_i)}$, and

(ii) characterizing the image of $\cB_{\VV,(V_i)}$ in $\cF_{\ell,(\La_i)}$. 

In this paper, we only say some words about problem (ii) in the 
case of integrable modules (see Remark \ref{rem:integrable}).  

We solve (i) completely in the case $N = 0$ (sections 
\ref{sect:N=0:n=0} and \ref{sect:N=0:n:arb}).  Out main results are 
Thm.\ \ref{thm:no:point} and Thm.\ \ref{thm:pts}, where we show that
flat sections correspond  bijectively to solutions of a heat equation.
This result is analogous to  that of \cite{FW} in the elliptic case
$g=1$; it is also close to the formulas of Verlinde and Verlinde 
(\cite{Verlinde}).  Our formulas also involve functions $\al(m,P_0,v_0)$ and
$\beta(m,P_0,v_0)$ on the moduli space $\cM_{g,1^2}$, which are
determined by differential  equations (Thm.\ \ref{thm:no:point}, Prop.\
\ref{koche}). They are homogeneous of fixed degree in the variable
$v_0$. Their exact meaning is unclear to us, but we propose conjectures
about them in Remark \ref{conjectures}.

In the case $N>0$, we construct flat sections of $\cF_{\ell,(\La_i)}$ as
integrals over twisted cycles over powers of $X$ of the functions
obtained in  Thm.\ \ref{thm:pts} (Thm.\ \ref{jahrzeit}). Our main tool
is the result that some flat sections  of $\cB_{\VV,(V_i)}$ can be
constructed as integrals of flat sections of $\cB_{\VV,(V_i),V_{-2},
\ldots,V_{-2}}$ (Prop.\ \ref{freud}).  This result generalizes those of
Cherednik (\cite{Ch}) in the case of genus $0$ and Varchenko and the 
second author (\cite{FeVa}) in the case $g = 1$. (An abstract version of this 
result can also be found in the work \cite{BFS}). We then obtain our
expression of flat sections  of $\cF_{\ell,(\La_i)}$ in Thm.\
\ref{jahrzeit}. This expression involves differential-evaluation
operators $\Phi_{(\La_i)}(w_1,\ldots,w_n|z_1,\ldots,z_n)$,  which
correspond to the insertion of a product $\prod_{j=1}^p f(w_j)
\prod_{j=1}^p e(z_j)$ in correlation functions (Lemma \ref{gardel}). 
It seems that comparing our formulas with explicit formulas in the 
case $k=1$ would yield integral identities, as it was done in 
\cite{FeVa}. 

In section \ref{sect:example}, we write an example of our formulas 
for a simple case ($g = 2$).

\subsubsection*{}
Let us now say some words on the topological aspects of our work. 

{\it Representations of the mapping class group.}
The following three types of representations of the mapping class group (MCG)
are usually thought to be closely related:  
 
(a) in \cite{Kohno}, Kohno constructed projective representations of the
MCG using braiding and fusing matrices of the genus zero  KZ equations,
and the method of Moore and Seiberg \cite{MS} for constructing such
representations; 

(b) in \cite{Lyub}, Lyubashenko, relying on work of Reshetikhin-Turaev
(\cite{Resh-Tur}), constructed projective representations of the MCG 
using the representation theory of quantum groups at roots of unity; 

(c) the monodromy of the KZB connection also provides projective
representations of the  MCG. 

The identification between representations of the MCG of (a) and (b)
should follow from the identification of the fusing and braiding
matrices, which can be checked in genera $0$ and $1$. 

On the other hand, our integral formulas might serve to understand the
connection (b)-(c).  This could be done using the geometric
interpretation of quantum groups of \cite{CFW}. The identification
(b)-(c) would be an analogue of the Kazhdan-Lusztig equivalence in genus
$\geq 1$ (\cite{KL}).

{\it Extensions of the mapping class group.} As we have seen, the KZB
connection constructed in \cite{BFM} gives representations of
$\pi_1(\det_{g,\nu}^*)$, which  is an  extension of the MCG. In
\cite{Harer}, Harer showed that the MCG is perfect and computed its
universal central extension (the central extensions of the MCG arising
in  natural projective representations were  later  computed by Masbaum
and Roberts \cite{Masbaum}). It is natural to believe that
$\pi_1(\det_{g,\nu}^*)$ identifies with this universal central
extension. 

On the other hand, if we denote by $\cW^{(2)}$ the subspace of
$\cM_{g,\nu\cdot 1^2}$ formed of the moduli such that the first point 
is not Weierstrass, Thm.\ \ref{thm:flatness} yields representations of
$\pi_1(\cM_{g,\nu\cdot 1^2} - \cW^{(2)})$. $\pi_1(\cM_{g,\nu\cdot 1^2} -
\cW^{(2)} )$ is also an extension of MCG.  

Projectivizations of both representations yield coinciding 
representations of the MCG. It is therefore natural to think that 

(a) the fundamental groups $\pi_1(\det_{g,\nu}^*)$ and 
$\pi_1(\cM_{g,\nu\cdot 1^2} - \cW^{(2)})$ 
coincide

and

(b) the representations of these groups provided by \cite{BFM}
and Thm.\ \ref{thm:flatness} coincide.

\subsubsection*{} We would like to thank B. Feigin for discussions about
this paper.  The first author would also like to thank A.-S.\ Sznitman
for  invitations to the FIM (ETHZ) in 1998 and June 1999, as well as
A.\ Alekseev for an invitation to ESI (Vienna) in August 1999, 
during which a part of this work was done.

\section{Flat version of the KZB connection and Weierstrass points}
\label{sect:1}

\begin{notation} For $P_0$ a point of a curve $X$, and $\omega$ a formal 
differential at $P_0$, we denote by $\langle \omega\rangle$ its
residue $\res_{P_0}(\omega)$ at the point $P_0$. For $\xi$ a vector and
$\omega^{(2)}$ a quadratic differential at the neighborhood of $P_0$, we
set $\langle \xi, \omega^{(2)} \rangle = \langle
\xi\omega^{(2)}\rangle$.  In the same way, if $f$ is a function at the
neighborhood of $P_0$, we set $\langle f,\omega \rangle = \langle
f\omega \rangle$. If $f$ and $\omega$ depend on variables  $z$ and $t$,
we will write $\langle f, \omega \rangle_z$ for $\res_{z = P_0}[(f\omega)(z,t)]$.  

If $\al$ is a differential of order $k$ at $P_0$ and $v_0$ is a tangent
vector at $P_0$, we denote by $(\al,v_0)$ the evaluation of $v_0$ on
$\al$. This expression is homogeneous of degree $k$ in $v_0$. 

Finally, if $\al$ is a formal differential of order $k$ at $P_0$, we define
its valuation as the smallest integer $k$ such that $t^k\al$ is regular at 
$P_0$, where $t$ is a local coordinate at $P_0$, and we denote it by 
val$(\al)$. 

\end{notation}

\subsection{Weierstrass points}

Let $X$ be a compact complex curve; let us denote its genus by $g$.  A
point $P_0$ of $X$ is not a Weierstrass point if there are holomorphic
differentials on $X$ with valuations $0,1,\cdots,g-1$ at $P_0$; in other
words, the Wronskian of holomorphic differentials does not vanish at $P_0$. 

Let $P_0$ be an arbitrary point of $X$, and let $t$ be a formal
coordinate on $X$ at $P_0$. Let $R$ be the subring of $\CC((t))$, formed
of the Laurent expansions at $P_0$ of the functions on $X$ regular
outside $P_0$. Let $\Omega^{out}$ be the subspace of $\CC((t))dt$ formed
by the forms on $X$, regular outside $P_0$. 

\begin{lemma} \label{lemma:decompositions}
The following statements are equivalent: 

1) $P_0$ is not a Weierstrass point of $X$. 

2) We have direct sum decompositions
\begin{equation} \label{decomp:forms}
 \CC((t))dt = \Omega^{out} \oplus \left( \CC t^{-1}dt \oplus
  t^{g}\CC[[t]] dt\right) 
\end{equation}
and 
\begin{equation} \label{decomp:funs}
  \CC((t)) = R \oplus \left( \oplus_{i = -g}^{-1} \CC t^{i} \oplus
    t\CC[[t]]\right).
\end{equation}
\end{lemma}

{\em Proof.} Let us show the equivalence of 1) and decomposition 
(\ref{decomp:forms}). 
Let us study $\Omega^{out} \cap (\CC t^{-1}dt + 
t^g \CC[[t]]dt)$. Let $\omega$ belong to this intersection. 
Then $\omega$ has the expansion $\omega_{-1} t^{-1} dt + \sum_{i\geq g}
\omega_i t^i dt$. Since $\res_{P_0}(\omega) = 0$, we have $\omega_{-1} = 0$. 
Therefore the valuation at $P_0$  of $\omega$ is $\geq g$. 
It follows that $\Omega^{out} \cap (\CC t^{-1}dt + 
t^g \CC[[t]]dt)$ consists of all regular forms on $X$ whose valuation at $P_0$
is $\geq 0$. This is zero iff $P_0$ is not Weierstrass. 

Let us now show that (\ref{decomp:forms}) implies (\ref{decomp:funs}).  
We will first show that 
$$
R \cap \left( \oplus_{i = -g}^{-1} \CC t^{i}
  \oplus t\CC[[t]]\right) = 0.
$$ 
Assume that this intersection is
nonzero, and let $\phi$ be a nonzero function in this intersection.
Let $j$ be its valuation at $P_0$. $j$ belongs to $\{-1, \ldots,
-g\}$. Since $P_0$ is Weierstrass, we have a holomorphic differential $\omega$ on $X$ with
valuation at $P_0$ equal to $-j-1$. Then $\phi\omega$ is a holomorphic
differential on $X$ with nonzero residue at $P_0$, a contradiction. It
follows that $R \cap \left( \oplus_{i = -g}^{-1} \CC t^{i}
  \oplus t\CC[[t]]\right) = 0$. 

(\ref{decomp:funs}) now follows from the facts that
$\Omega^{out}$ and $R$ on one side, and $ \left( \CC t^{-1}dt \oplus
  t^{g}\CC[[t]] dt\right)$ and $\left( \oplus_{i = -g}^{-1} \CC t^{i}
  \oplus t\CC[[t]]\right)$ on the other side, are each other's
annihilators for the residue pairing.  \hfill \qed \medskip

\subsection{``Weierstrass'' normalization of the Sugawara tensor}

The residue defines a nondegenerate pairing between $\CC((t))$ and
$\CC((t))dt$. Assume that $P_0$ is not a Weierstrass point of $X$. As we
have seen, the spaces $R$ and $\Omega^{out}$, and $$ \left( \oplus_{i =
-g}^{-1} \CC t^{i} \oplus t\CC[[t]]\right)\  \on{and} \ \left( \CC
t^{-1}dt \oplus t^{g}\CC[[t]] dt\right) $$ are annihilators of each
other for this pairing.

Let $(\omega_i)$ and $(r_i)$ be bases of 
$$ \left( \oplus_{i = -g}^{-1} \CC t^{i} \oplus t\CC[[t]]\right)\ 
\on{and} \ \left( \CC t^{-1}dt \oplus t^{g}\CC[[t]] dt\right) , 
$$ 
such that the sequences (val$(r_i))_i$ and (val$(\omega_i))_i$ tend to 
infinity, and let $(\omega^i)$ and $(r^i)$ be the dual bases of
$\Omega^{out}$ and $R$.  We have the equality (in $\CC[[z^{\pm1},w^{\pm1}]]dw$)
$$ \delta(z,w) dz = \sum_i \omega^i(z) r_i(w) + \sum_i \omega_i(z)
r^i(w) , 
$$
where $\delta(z,w)dw = \sum_{i\in\ZZ} z^i w^{-i-1} dw$.

Let $\bar\G$ be a semisimple Lie algebra with nondegenerate invariant bilinear 
form $\langle , \rangle_{\bar\G}$, and define $\G$ as the canonical 
central extension of $\bar\G \otimes
\CC((t))$; it is the direct sum $\bar \G \otimes \CC((t)) \oplus \CC
K$, endowed with the Lie bracket such that $K$ is central and
$$[x\otimes f, x'\otimes f'] = [x,x']\otimes ff' + \langle
x,x'\rangle_{\bar\G} \langle df f' \rangle K.$$ Let $\VV$ be any
$\G$-module on which the Lie subalgebra $\bar\G \otimes \CC[[t]]$ acts
locally finitely. We will denote by $x[f]$ the element $(x\otimes
f,0)$ of $\G$ for $x$ in $\bar\G$, $f$ in $\CC((t))$.

Define, for $\xi$ in $\CC((t)){\pa\over{\pa t}}$, 
\begin{equation} \label{def:TR}
T_{R}[\xi] = {{1}\over{2\kappa}} \sum_\al \sum_i \left( x^\al
  [\xi\omega^i] x^\al[r_i] + x^\al[r^i]x^\al[\xi\omega_i]\right) ,
\end{equation} 
with $(x^\al)$ an orthonormal basis of $\bar\G$ and $\kappa = k +
h^\vee$, $h^\vee$ the dual Coxeter number of $\bar\G$. We call $T_R[\xi]$
the normalized Sugawara tensor. 

Then we have, for any $f$ in $\CC((t))$,
\begin{equation} \label{vf}
 [T_R[\xi], x[f]] = - x[\xi(f)],
\end{equation}
 where $\xi(f)$ denotes the action of the vector field $\xi$ on the
function $f$ (the product of a form $\omega$ and a vector $\xi$, which
is a function, is simply denoted $\xi\omega$; we have therefore
$\xi(f) = \xi df$).

Define $\bar\G(R)$ as image of $\bar\G \otimes R$ in $\G$ by the map
$x\otimes f \mapsto x[f]$; this is a Lie subalgebra of $\G$.

One of the useful features of $T_R[\xi]$ is: 

\begin{lemma}
If $\xi$ is a vector field on $X$, regular outside $P_0$, 
$T_R[\xi]$ belongs to $\bar\G(R)U\G$. 
\end{lemma}

{\em Proof.} The follows from the fact that $\xi\omega^i$
belongs to $R$. \hfill \qed\medskip

\subsubsection{Another expression of the normalized Sugawara tensor}

Here and below, $d_{z'}$ denotes the partial exterior differentiation
with respect to the variable $z'$.

\begin{lemma}
  Set
\begin{equation} \label{bar:omega}
 \bar\omega(z,z') = \sum_i \omega^i(z) dr_i(z') .  
 \end{equation}  
Define
  $T_{\bar\omega}(z)$ as 
  $$ {1 \over{2\kappa}} \lim_{z'\to z} \left[ \sum_\al x^\al(z)
  x^\al(z') - k (\dimm\bar\G) \bar\omega(z,z') \right] .$$ Then
  the matrix elements of $T_R[\xi]$ and $\langle \xi , 
  T_{\bar\omega}(z) \rangle$ in any highest weight representation 
  of $\G$ are well-defined and coincide. 
\end{lemma}

{\em Proof.} 
Set
$$ x(z) = \sum_i x[\al^i]\eps_i(z), \quad x(z)_{z\to \Omega^{in}} =
\sum_i x[r^i]\omega_i(z), \quad x(z)_{z\to \Omega^{out}} = \sum_i
x[r_i]\omega^i(z),
$$ for $(\al^i),(\eps_i)$ dual bases of $\CC((t))dt$ and $\CC((t))$.
We have 
\begin{align*}
  2 \kappa T_{\bar\omega}(z) = & \sum_\al x^\al(z)_{z\to
    \Omega^{out}} x^\al(z) + x^\al(z) x^\al(z)_{z\to \Omega^{in}} \\ &
  + \limm_{z'\to z}\sum_\al [x^\al(z)_{z\to \Omega^{out}}, x^\al(z')]
  - k (\dimm \bar\G)\bar\omega(z,z').
\end{align*}
As the limit vanishes, we have 
$$ T_{\bar\omega}(z) = {1\over{2\kappa}} \sum_\al x^\al(z)_{z\to
  \Omega^{out}} x^\al(z) + x^\al(z) x^\al(z)_{z\to \Omega^{in}}.
$$
We have then 
\begin{align*}
  & \langle T_{\bar\omega}(z) , \xi \rangle = {1\over{2\kappa}}
  \sum_\al x^\al[r^i] x^\al[\eps^j] \langle \omega_i \al_j, \xi
  \rangle + x^\al[\eps^j] x^\al[r_i] \langle \omega^i \al_j, \xi
  \rangle \\ & = {1\over{2\kappa}} \sum_i x^\al[r^i] x^\al[\xi
  \omega_i] + x^\al[\xi\omega^i] x^\al[r_i] = T_R[\xi].
\end{align*}
\hfill \qed \medskip

\subsection{Vector fields on moduli spaces}

Let $g,n$ be integers $\geq 0$. Let $\cM_{g,1^\infty,n\cdot 1^2}$ be
the space of isomorphism classes of systems 
$$
m = (X,t,P_1,v_1,\ldots,P_n,v_n)
$$ 
of a genus $g$ curve $X$, a marked point $P_0$ with a coordinate $t$ at
$P_0$, and other marked points $P_1,\ldots,P_n$ distinct from $P_0$
together with tangent vectors $v_1,\ldots,v_n$ at these points.
We will denote $\cM_{g,1^\infty,0\cdot 1^2}$ simply as $\cM_{g,1^\infty}$; 
it is the space of isomorphism classes of pairs 
$$
m_\infty = (X,t) ;  
$$
$\cM_{g,1^\infty}$ can also be viewed as the space of subrings $R$ 
of the Laurent series field $\CC((t))$, such that $R$ is isomorphic to a 
coordinate ring $H^0(X - \{P_0\},\cO_X)$, with $X$ a curve of genus $g$, 
and the inclusion $R\subset \CC((t))$ is the inclusion of $R$
in the local field of $X$ at $P_0$. 

In the same way, $\cM_{g,1^\infty,n \cdot 1^2}$ can be identified with the space of
systems $m = [R, \bar\chi_1,\ldots,\bar\chi_n]$, where $R\subset
\CC((t))$  is the same subring of $\CC((t))$ as above, and  
$\bar\chi_i$ are morphisms from $R$ to $\CC[\eta] /
(\eta^2)$. (Indeed, a morphism $\bar\chi$ from $R$ to $\CC[\eta]  / (\eta^2)$
is the same as the data of a pair $(\chi,\pa_{\bar\chi})$ of a morphism
$\chi$ from $R$ to $\CC$ and a map $\pa_{\bar\chi}$ from $R$ to $\CC$, 
such that $\pa_{\bar\chi}(fg) = f\pa_{\bar\chi}(g) + \pa_{\bar\chi}(f)g$; $\chi$
corresponds to a point of Spec$(R)$ and $\pa_{\bar\chi}$ to a tangent 
vector at this point.)

Define $\xi\mapsto [\xi]$ as the map from $\CC((t)) {{\pa}\over{\pa t}}$ to 
$\Vect(\cM_{g,1^\infty, n\cdot 1^2} )$, such that for $\xi$ in 
$\CC((t)) {{\pa}\over{\pa t}}$, the infinitesimal translate of 
$[R, \bar\chi_1,\ldots,\bar\chi_n]$ by $[\xi]$ is 
$[(1+\eps\xi) (R)\subset \CC((t)), \bar\chi_1\circ (1-\eps\xi),\ldots,
\bar\chi_n\circ (1-\eps\xi)]$.

\begin{prop} (\cite{Concini,BS,KNMY,Konts})
  $\xi\mapsto [\xi]$ defines a Lie algebra morphism from $\CC((t)){\pa \over{\pa
  t}}$ to $\Vect(\cM_{g,1^{\infty},n\cdot 1^2})$.   Let $\ev_m$ be the
  evaluation of a vector field over   $\cM_{g,1^\infty,n\cdot 1^2}$ at the
  point  $m\in \cM_{g,1^\infty,n\cdot 1^2}$.   The kernel of $\xi\mapsto \ev_m([\xi])$ 
  is the space of vector fields   $\xi$, which admit a regular
  continuation to $X - \{P_0\}$ and vanish  to second order at the $P_s$. 
\end{prop}

\subsection{The flat connection $(\cB_{\VV,(V_s)},\nabla^{\VV,(V_s)})$} 
\label{sect:before}

\subsubsection{The sheaf $\cB_{\VV,(V_s)}$} \label{class:C}

Let $\cM_{g,1^\infty,n}$ be the space of isomorphism classes of systems
$m = (X,t,P_s)$ of a genus $g$ curve $X$, a point $P_0$ on it, with a
local coordinate $t$, and of $n$ distinct points $P_s$, distinct from
$P_0$. $\cM_{g,1^\infty,n}$ can be viewed as the space of systems $m = [R,
\chi_s]$ of a subring $R$ of $\CC((t))$ as above,  and of a collection
of $n$ morphisms $(\chi_s)_{s = 1,\ldots,n}$  from $R$ to $\CC$. Let $p$
be the morphism from $\CC[\eta]/(\eta^2)$ to $\CC$, such that $p(\eta) =
0$. We have a natural projection $\pi_{2\to 1}$ from
$\cM_{g,1^\infty,n\cdot 1^2}$ to $\cM_{g,1^\infty,n}$, such that
$\pi_{2\to 1}([R, \bar\chi_s]) = [R,p\circ\bar\chi_s]$.

Say that a representation $(\rho_\VV,\VV)$ of $\G$ is of the class (C) if the following
conditions are satisfied: 

\medskip 

1) for any $v\in \VV$, there exists an integer $n(v)$ such that 
if $n\geq n(v)$, $\rho_\VV(\bar\G\otimes t^n)(v) = 0$; (this implies that the
components of the Sugawara tensor act as well-defined endomorphisms of $\VV$); 

2) let $T_0$ be the zero-mode of the Sugawara tensor. For $n$ and $d$ 
integer numbers, let us set $\VV[n,d] = \Ker(\rho_\VV(T_0) - d) \cap 
\Ker(\rho_\VV(h[1]) - n)$; then $\VV = \oplus_{n,d\in\ZZ} \VV[n,d]$, and 
for any integer $n$, there is $d(n)$ such that for $d\geq d(n)$, $\VV[n,d] = 0$. 

\medskip 

Modules of the class (C) contain both twisted Weyl modules over $\G$ and
integrable representation of $\G$. Moreover, they are such that the positive
Fourier modes of the Sugawara tensor act  locally nilpotently.

Fix a representation $\VV$ of $\G$ of the class (C) and  representations
$(\pi_s,V_s)$ of $\bar\G$. We will assume that the Casimir element
$\sum_\al (x_{\al})^2$ acts on $V_s$ as a scalar $\Delta_s id_{V_s}$
(this is the case if the $V_s$ are irreducible modules). We
have an evaluation morphism $\ev_{P_s}$ from $\bar\G(R) = \bar\G\otimes
R$ to $\bar\G$, defined as $id \otimes \chi_s$.   $\pi_s\circ\ev_s$
defines then a $\bar\G(R)$-module structure on $V_s$.  We denote by
$V_s^{(P_s)}$ the resulting $\bar\G(R)$-module. Define
$\cB_{\VV,(V_s)}^0$ as the sheaf over $\cM_{g,1^\infty,n}$, whose sections
over an open subset $U$ are maps $\psi : U \to \cV^*$, 
such that $\psi(m)$ is  $\bar\G(R)$-invariant for any $m$ in $U$, where 
$\cV = \VV\otimes (\otimes_s V_s^{(P_s)})$. 

Define $\cB_{\VV,(V_s)}$ as the pull-back sheaf $\pi_{2\to
1}^*(\cB_{\VV,(V_s)}^0)$ of  $\cB_{\VV,(V_s)}^0$ to
$\cM_{g,1^\infty,n\cdot 1^2}$.  Local sections of this sheaf are 
maps $\psi : [R,\bar\chi_1,\ldots,\bar\chi_n] \to 
[\VV \otimes (\otimes_s V_s)]^*$, such that for each 
$v\in\cV$, $m\mapsto \langle \psi(m),v\rangle$ is smooth, and each  
$\psi ([R,\bar\chi_1,\ldots,\bar\chi_n])$ is $\bar\G(R)$-invariant. 
Here $V_s$ is endowed with the $\bar\G(R)$-module 
structure given by $\pi_s\circ (id_{\bar\G} \otimes \{p\circ\bar\chi_s\}
)$ ($id_{\bar\G} \otimes \{p\circ\bar\chi_s\}$ is a map from  $\bar\G(R)
= \bar\G\otimes R$ to $\bar\G$, and $\pi_s$ maps $\bar\G$ to
$\End(V_s)$).

\subsubsection{The locus $\cW^\infty$ of Weierstrass points} \label{38}

Define $\cW$ as the subset of $\cM_{g,1}$ formed of the pairs 
$(X,P_0)$ such that $P_0$ is a Weierstrass point of $X$.

Let $\cM^{(a)}_{g,1^2}$, resp.\ $\cM^{(a)}_{g,1}$
 be the moduli space of systems $(X,P_0,v_0,
\{A_a\})$ of a genus $g$ curve with a marked point, a nonzero tangent
vector at this point, and a system of $a$-cycles (resp.\ of systems
$(X,P_0,\{A_a\})$).  Let $\pi'_{2\to 1}$
be the natural projection from $\cM^{(a)}_{g,1^2}$ to $\cM^{(a)}_{g,1}$. 

The Wronskian of the Abelian differentials defines a function $W(X,P_0,
v_0,\{A_a\})$ on $\cM^{(a)}_{g,1^2}$. This function is homogeneous
of degree ${{g(g+1)}\over 2}$ in $v_0$. 

Then $(\pi'_{2\to 1})^{-1}(\cW)$ is defined by the equation
$W(X,P_0,v_0,\{A_a\}) = 0$. It follows that  $(\pi'_{2\to 1})^{-1}(\cW)$
has codimension $1$ in the smooth part of  $\cM^{(a)}_{g,1^2}$, and that
 $\cW$ has codimension $1$ in the smooth part of  $\cM_{g,1}$. 

We will define $\cW^\infty$ as the preimage $\pi_{\infty}^{-1}(\cW)$
of $\cW$ by the natural projection map $\pi_\infty:\cM_{g,1^\infty, n \cdot 
1^2}$, such that
$$
\pi_\infty(X,t,P_i,v_i) = (X,P_0). 
$$

\subsubsection{The connection $\nabla^{\VV,(V_s)}$}

For $x$ in $U\G$, we will denote by $x^{(0)}$ the element
$\pi_{\VV}(x)\otimes(\otimes_s id_{V_s})$  of $\End(\VV \otimes
(\otimes_s V_s) )$ and by $x^{(s)}$ (sometimes $x^{(P_s)}$) 
the element $id_{\VV} \otimes (\otimes_{s'<s} id_{V_{s'}})
\otimes \pi_{V_s}(x)\otimes (\otimes_{s'>s} id_{V_{s'}})$. 

\begin{lemma}  \label{secu}
For $m\mapsto \psi(m)$ a local section of $\cB_{\VV,(V_s)}$,  and $\xi$ in $\CC((t)) \pa
/\pa t$, the value at $m$ of  $$ \pa_{[\xi]}\psi(m) -
\psi(m)\circ T_R[\xi]^{(0)} $$ belongs to $(\cV^*)^{\bar\G(R)}$ and
only depends on the value at $m$ of $[\xi]$.   
\end{lemma}

{\em Proof.} Let $\eps$ be a variable such that $\eps^2 = 0$. Let us set 
$R_\eps = (1 + \eps\xi)(R)$. Then $\bar\G(R_\eps) = \Ad(1 - \eps T[\xi])
[\bar\G(R)]$. It follmows that $\pa_{[\xi]}\psi(m) - \psi(m)\circ T_R[\xi]^{(0)}$ 
belongs to is $\bar\G(R)$-invariant. 

Let $R_{(P_s)}$ be the subspace of $R$ consisting of the
functions vanishing at each $P_s$.  Let $\bar\G(R_{(P_s)})$ be the image 
of $\bar\G\otimes R_{(P_s)}$ by the inclusion of $\bar\G\otimes R$ in $\G$. 
Let us show that for $\xi$ a vector
field in $\CC((t))\pa / \pa t$, with a regular prolongation to  $X -
\{P_0\}$ and double zeroes at each $P_s$, $T_R[\xi]$ belongs to
$\bar\G(R_{(P_s)}) U\G$.   Let $(r^i_{(P_s)} )$ be a basis of
$R_{(P_s)}$ and let $(\omega^i_{(P_s)})$ be a basis of the space of
forms $\Omega_{(P_s)}$ on $X$,  regular except at $P_0$ and the $P_s$,
and with at most simple poles at the $P_s$. $R_{(P_s)}$ and
$\Omega_{(P_s)}$ are the annihilators of each other by the residue
pairing.  We may therefore choose $(r_{i;(P_s)} )$ and
$(\omega_{i;(P_s)})$ in $\cK$ and  $\Omega$ so that $(r^i_{(P_s)} ;
r_{i;(P_s)} )$ and $(\omega_{i;(P_s)}; \omega^i_{(P_s)})$ are dual bases
of  $\cK$ and $\Omega$. As the $\xi \omega^i_{(P_s)}$ and the
$r^i_{(P_s)}$  belong to $R_{(P_s)}$, 
$$ 
T'[\xi] = {1\over{2\kappa}} \sum_{\al,i} x^\al[\xi
\omega^i_{(P_s)}] x^\al[r_{i;(P_s)}] +
x^\al[r^i_{(P_s)}]x^\al[\xi\omega_{i;(P_s)}]
$$
belongs to $\bar\G(R_{(P_s)}) U\G$.

On the other hand, let $\omega^{\prime i}_{(P_s)}$ and  $r'_{i;(P_s)}$
be lifts to $\Omega^{out}_{(P_s)}$ and $R^{out}$  of dual bases of
$\Omega^{out}_{(P_s)}  / \Omega^{out}$ and $R^{out} / R^{out}_{(P_s)}$. 
$T_R[\xi] - T'[\xi]$ is equal to 
$$ 
{1\over{2\kappa}} \sum_{\al,i} 
[x^\al[\xi \omega^{\prime i}_{(P_s)} ] ,  x^\al[r'_{i;(P_s)}]] . 
$$

This is proportional to $\sum_i  \langle \xi \omega^{\prime i}_{(P_s)} 
dr'_{i;(P_s)} \rangle$, which is zero because the  $\xi \omega^{\prime
i}_{(P_s)}$ belong to $\Omega^{out}$. It follows that  $T_R[\xi]$ is
equal to $T'[\xi]$, and therefore belongs to  $\bar\G(R_{(P_s)})
U\G$.  

For $x$ in $\bar\G(R_{(P_s)})$ and $v$ in $\VV$, 
$\psi(m)(xv)$  vanishes, therefore $\psi(m)\circ T_R[\xi]^{(0)}$
vanishes if $\xi$ belongs to the kernel of $\xi \mapsto \ev_m ([\xi])$.
This proves the Lemma. 
\hfill \qed \medskip

It follows from Lemma \ref{secu} that the formula 
$$
\nabla^{\VV,(V_s)}_{[\xi]}\psi(m) = \pa_{[\xi]}\psi(m) 
- \psi(m)\circ T_R[\xi]^{(0)} 
$$ 
defines a connection $\nabla^{\VV,(V_s)}$ over the restriction of
$\cB_{\VV,(V_s)}$ to $\cM_{g,1^\infty} - \cW^\infty$.

\subsubsection{Flatness of $\nabla^{\VV,(V_s)}$}

We will show: 

\begin{thm} \label{thm:flatness}
  $\nabla^{\VV,(V_s)}$ is a flat connection on the restriction of  
$\cB_{\VV,(V_s)}$ to $\cM_{g,1^\infty} - \cW^\infty$. 
\end{thm}

{\em Proof.}
Let us compute the curvature of $\nabla^{\VV,(V_s)}$. We have
$$ \left( 
[\nabla^{\VV,(V_s)}_\xi , \nabla^{\VV,(V_s)}_\eta] - \nabla^{\VV,(V_s)}_{[\xi,\eta]} 
\right) \psi
= \psi \circ \left( 
- T_R[[\xi,\eta]] - [T_R[\xi], T_R[\eta]] - \pa_\xi T_R[\eta] + \pa_\eta
T_R[\xi] \right)^{(0)} .
$$
We may assume that $r_i$ and $\omega_i$ are independent of $R$. Then 
$$ \pa_\xi r^i = \xi(r^i) - \sum_j \langle \xi(r^i), \omega_j\rangle
r^j,
$$
and
$$ \pa_\xi \omega^i = d(\xi\omega^i) - \sum_j \langle d(\xi\omega^i),
r_j \rangle \omega^j. 
$$
Then 
\begin{align*} 
  \pa_\xi T_R[\eta] = & {1\over{2\kappa}} \sum_\al\sum_{i} x^\al
  \left[\eta d(\xi\omega^i) - \sum_j \langle d(\xi\omega^i),r_j
    \rangle \eta \omega^j\right] x^\al[r_i] \\ & + x^\al\left[
    \xi(r^i) - \sum_j \langle \xi(r^i), \omega_j\rangle
    r^j\right]x^\al[\eta\omega_i].
\end{align*}
On the other hand, since $\xi(\eta\omega) - \eta(\xi\omega) =
[\xi,\eta]\omega$, we have 
\begin{align*} 
  [T_R[\xi], T_R[\eta]] + T_R[[\xi,\eta]] = - {1\over{2\kappa}}
  \sum_\al \sum_i & x^{\al}[\eta(\xi\omega^i)] x^\al[r_i] +
  x^{\al}[\eta\omega^i] x^\al[\xi(r_i)] \\ & + x^\al[\xi(r^i)]
  x^\al[\eta\omega_i] + x^\al[r^i] x^\al[\eta(\xi\omega_i)].
\end{align*} 
It follows that 
\begin{align*}
  & -2\kappa( \nabla^{\VV,(V_s)}_{[\xi,\eta]} - [\nabla^{\VV,(V_s)}_\xi,\nabla^{\VV,(V_s)}_\eta]
  ) \psi \\ & 
  =  \psi \circ \pi_{\VV} \lbrace - \sum_\al \sum_i \{ x^{\al}[\eta(\xi\omega^i)] x^\al[r_i] +
  x^{\al}[\eta\omega^i] x^\al[\xi(r_i)] + x^\al[\xi(r^i)]
  x^\al[\eta\omega_i] \\ & + x^\al[r^i] x^\al[\eta(\xi\omega_i)] \} 
  \\ & +
  \sum_\al\sum_{i} x^\al \left[\eta d(\xi\omega^i) - \sum_j \langle
    d(\xi\omega^i),r_j \rangle \eta \omega^j\right] x^\al[r_i] +
  x^\al\left[ \xi(r^i) - \sum_j \langle \xi(r^i), \omega_j\rangle
    r^j\right]x^\al[\eta\omega_i] \\ & - \sum_\al\sum_{i} x^\al
  \left[\xi d(\eta\omega^i) - \sum_j \langle d(\eta\omega^i),r_j
    \rangle \xi \omega^j\right] x^\al[r_i] + x^\al\left[
    \eta(r^i) - \sum_j \langle \eta(r^i), \omega_j\rangle
    r^j\right]x^\al[\xi\omega_i] \} \rbrace .
\end{align*}
For $f$ in $\CC((t))$ and $\omega$ in $\CC((t))dt$, denote by
$f_{(out)}$ and $f_{(in)}$, resp.\ $\omega_{(out)}$ and $\omega_{(in)}$ the
projections of $f$ and $\omega$ on the first and second spaces of
decompositions (\ref{decomp:funs}) and (\ref{decomp:forms}).
$-2\kappa(\nabla^{\VV,(V_s)}_{[\xi,\eta]} - [\nabla^{\VV,(V_s)}_\xi,\nabla^{\VV,(V_s)}_\eta])
\psi$
can then be expressed as
\begin{align*}
  \psi\circ \left( \right. \sum_\al \sum_i & - x^\al[\eta\omega^i] x^\al[\xi(r_i)] - x^\al[r^i]
  x^\al[\eta(\xi\omega_i)] \\ & - x^\al[\eta \{
  d(\xi\omega^i)\}_{(out)}] x^\al[r_i] - x^\al[\{ \xi(r_i)\}_{(out)}]
  x^\al[\eta \omega_i] \\ & - x^\al[\xi \{ d(\eta\omega^i)\}_{(in)}]
  x^\al[r_i] - x^\al[\{ \eta(r_i)\}_{(in)}] x^\al[\xi \omega_i]
  \left. \right)^{(0)} ,
\end{align*}
which is equal to 
\begin{align} \label{eetion}
  \psi \circ \left( \right. \sum_\al \sum_i & - x^\al[\eta \omega^i] x^\al[\{\xi(r_i)\}_{(out)}] -
  x^\al[r^i] x^\al[\eta\{d(\xi\omega_i)\}_{(out)}] \\ & \nonumber -
  x^\al[\xi\{d(\eta\omega^i)\}_{(in)}]x^\al[r_i] -
  x^\al[\{\eta(r^i)\}_{(in)}] x^\al[\xi\omega_i] \left. \right)^{(0)}.
\end{align}

Write (\ref{eetion}) in the form $\psi\circ ( \sum_\al \sum_i
x^\al[a_i]x^\al[b_i])^{(0)}$.  $\sum_\al \sum_i x^\al[b_i] x^\al[a_i]$ also
belongs to the completion of $U\G$ with respect to the family of left
ideals generated by the $x^\al[t^i],i\geq N$. Let us show
that $\sum_i a_i \otimes b_i + b_i \otimes a_i =0$: this identity is
written as
\begin{align} \label{tros}
  \sum_i & \eta\omega^i \otimes \{\xi(r_i)\}_{(out)} + r^i \otimes \eta
  \{d(\xi\omega_i)\}_{(out)} \\ & \nonumber + \xi \{d(\eta
  \omega^i)\}_{(in)} \otimes r_i + \{\eta(r_i)\}_{(in)} \otimes
  \xi\omega_i + (z \leftrightarrow w) = 0,
\end{align} 
which is proved as follows: for $\omega,\omega'$ in $\CC((t)) dt$, 
\begin{align*}
  & \langle \omega \otimes \omega', \on{lhs\ of\ }(\ref{tros})\rangle
  \\ & = \langle d(\{\eta \omega\}_{(in)}), \xi\omega'_{(in)}\rangle +
  \langle d(\eta \omega_{(in)}), \{\xi\omega'\}_{(in)}\rangle \\ & -
  \langle d(\{\xi\omega\}_{(out)}), \eta\omega'_{(out)}\rangle - \langle
  d(\eta\omega_{(out)}), \{\xi\omega'\}_{(out)} \rangle + (\omega
  \leftrightarrow \omega') = 0.
\end{align*}
(\ref{eetion}) is then equal to ${k\over 2} \dimm\bar\G\sum_i \langle
da_i, b_i\rangle$, which is
\begin{align*}
  {k\over 2} (\dimm\bar\G) \sum_i & \langle d(\eta\omega^i) ,
  \{\xi(r_i)\}_{(out)} \rangle + \langle \eta(r_i),
  \{d(\xi\omega_i)\}_{(out)} \rangle \\ & - \langle
  \{d(\eta\omega^i)\}_{(in)} , \xi(r_i)\rangle - \langle
  \{\eta(r^i)\}_{(in)}, d(\xi\omega_i) \rangle ,
\end{align*}
and is equal to zero. \hfill \qed\medskip


\begin{remark} {\it Descent to $\cM_{g,(n+1) \cdot 1^2}$.}
Let $\cM_{g,(n+1)\cdot 1^2}$ the the moduli space of systems
$(X,t^{(2)},P_i,v_i)$ of a curve of genus $g$, the jet to second order
of a coordinate at a point $P_0$ of $X$,  and of $n$ points $P_i$
together with tangent vectors at these points. $\cM_{g,(n+1)\cdot 1^2}$
can be identified with the moduli space of systems $(X,P_i,v_i)$ of a
genus $g$ curve, $n+1$ marked points and a tangent vector at these
points, with $v_0\neq 0$.  Let $\pi_2$, resp.\ $\pi_{\infty,2}$ be the
canonical projections from $\cM_{g,(n+1)\cdot 1^2}$ to $\cM_{g,n}$, and
from $\cM_{g,1^\infty,n\cdot 1^2}$ to $\cM_{g,(n+1)\cdot 1^2}$.

The map $\pi_{\infty,2}$ has simply connected fibers, and maps
$\cW^\infty$ to $\pi_2^{-1}(\cW)$. For any section $\sigma$ of
$\pi_{\infty,2}$, we may therefore consider the sheaf $\sigma^*
\cB_{\VV,(V_s)}$ on $\cM_{g,1^2}$; it is endowed, over
$\cM_{g,(n+1)\cdot 1^2} - \pi_2^{-1}(\cW)$, with the flat connection
$\sigma^* \nabla^{\VV,(V_s)}$. For $\sigma'$ another section of 
$\pi_{\infty,2}$, the connections  $(\sigma^*\cB_{\VV,(V_s)},
\sigma^*\nabla^{\VV,(V_s)})$ and $(\sigma^{\prime *}\cB_{\VV,(V_s)},
\sigma^{\prime *}\nabla^{\VV,(V_s)})$  are isomorphic (the connection
can be integrated from $\sigma(m)$ to $\sigma'(m)$ by ordered path
exponentials because the positive Fourier modes of the Sugawara  tensor
act on $\VV$ locally nilpotently).   
\end{remark}

\subsection{Behaviour around marked points}

Recall that $\cB_{\VV,(V_s)}$ is the pull-back sheaf $\pi_{2\to 1}^*(\cB^0_{\VV,(V_s)})$,
where $\pi_{2\to 1}$ is the natural projection of
$\cM_{g,1^\infty,n\cdot 1^2}$ on $\cM_{g,1^\infty,n}$.  Let us study
the variation of flat sections of $\nabla^{\VV,(V_s)}$ along the fibres of
$\pi_{2\to 1}$.

\begin{prop}
  Let $\Delta_s$ be the value in $V_s$ of the Casimir element  $\sum_\al
  x_\al^2$. Let for $\la = (\la_s)$ in $(\CC^\times)^n$, $m\mapsto m_\la$
  be the transformation of $\cM_{g,1^\infty,n\cdot 1^2}$ associating to 
  $m = (X,t, P_s,v_s)$,  $m_\la = (X,t, P_s,\la_s v_s)$.   Assume that
  $m\mapsto \psi(m)$ is a flat section of $\cB_{\VV,(V_s)}$.  We have (using the natural
  identification of the fibres of  $\cB_{\VV,(V_s)}$ at $m$ and $m_\la$)  
  $$
  \psi(m_\la) = \prod_s \la_s^{ {\Delta_s\over{2\kappa}} } \psi(m) . 
  $$
\end{prop}
 
{\em Proof.} Let $\xi_s$ be a vector field on $X$, regular on  $X -
\{P_0\}$, with expansion at $P_t$ given by $\xi_s = \delta_{st}  z_t
\pa/\pa z_t + o(z_t)\pa/\pa z_t$, where $z_t$ is a local coordinate at
$P_t$. Then $\pa_{[\xi_s]}\psi(m) =  - {{\pa \psi(m_\la)} \over
{\pa\la_s}}_{|\la = 1}$. 

On the other hand, since the $\psi(m)\circ x^\al[r^i_{(P_s)}]^{(0)}$ 
and $\psi(m)\circ x^\al[\xi_s \omega^i]^{(0)}$ vanish, 
$\psi(m)\circ T_R[\xi_s]^{(0)}$ is equal to 
$$
{1\over{2\kappa}}\psi(m)\circ \left( \sum_{\al,i} 
x^\al[\xi_s \omega^{\prime i}_{(P_i)}] x^\al[r'_{i;(P_s)}]
\right)^{(0)} , 
$$
which is equal to 
$$ 
{1\over{2\kappa}} \sum_t \sum_i [\xi_s \omega^{\prime i}_{(P_s)} ]
(t) r'_{i;(P_s)}(t) \psi(m) \circ [\sum_\al (x^{\al(t)})^2] .
$$
We may choose the bases $(\omega^{\prime i}_{(P_s)})$ and 
$(r'_{i;(P_s)})$ as follows: these bases are indexed by $s$ in 
$\{1,\ldots,n\}$, $r'_{s;(P_s)}(t) = \delta_{st}$, and 
$\omega^{\prime s}_{(P_s)}$ has the expansion 
$$
\omega^{\prime s}_{(P_s)}  (z_t)  = 
- \delta_{st}{{dz_t}\over{z_t}} + O(1) 
$$
near $P_t$. Therefore, $\psi(m)\circ T_R[\xi_s]$ is equal to 
${{\Delta_s}\over{2\kappa}}\psi(m)$. Therefore, we get  
$$
{{\pa \psi(m_\la)} \over {\pa\la_s}}_{|\la = 1} 
= 
{{\Delta_s}\over{2\kappa}}\psi(m)
$$ 
for each $s = 1,\ldots,n$, which implies the Proposition.  
\hfill \qed \medskip

\subsection{Singularity around the locus of Weierstrass points}

Let us analyze the restriction of $\sigma^*\nabla^{\VV,(V_s)}$ to a fiber of the projection
$\pi_2 : \cM_{g,(n+1)\cdot 1^2} \to \cM_{g,n\cdot 1^2}$ such that 
$\pi_2(X,P_0,v_0,P_i,v_i) = (X,P_i,v_i)$, more precisely, its singularity 
around the subvariety $(\pi_{2\to 1})^{-1}(\cW)$.  

Recall the we defined in function $W$ on $\cM_{g,1^2}^{(a)}$ (sect.\ \ref{38}). 
Since $W$ gets multiplied by a constant under modular transformations, 
${{dW}\over{W}}$ is the pull-back of a form over $\cM_{g,1^2}$, that we denote 
$({{dW}\over{W}})_1$. Let $\mu$ be the map from $\cM_{g,(n+1)\cdot 1^2}$ to 
$\cM_{g,1^2}$ defined by $\mu(X,P_0,v_0,P_i,v_i) = (X,P_0,v_0)$; we set 
$$
\left( {{dW}\over{W}} \right)_{n+1} = \mu^* [\left({{dW}\over{W}}
\right)_{1}] .  
$$ 

\begin{prop} \label{comp:log}
  Let $U$ be an open subset of the smooth part of $\cM_{g,(n+1)\cdot 1^2}$.  
 $\sigma^*\nabla^{\VV,(V_s)}$ has the expansion
  \begin{align} \label{sing} 
  & \sigma^*  \nabla^{\VV,(V_s)} = d - {{k\dimm\bar\G}\over{2\kappa}}
    \left({{dW }\over W }\right)_{n+1} id_{\cB_{\VV,(V_s)}} +
    \on{regular\ form\ on\ }U\on{\ with\ values} 
   \\ & \nonumber
   \on{in \ End}(\cB_{\VV,(V_s)}).
  \end{align} 
  On the other hand, the projectivization of $\sigma^*\nabla^{\VV,(V_s)}$ is regular on
  $U$.
\end{prop}

{\em Proof.}  A Weierstrass point $P_{0}$ is called normal if there
exists a basis of the space of holomorphic differentials on $X$ with
valuations at $0,1,\ldots,g-2,g$ at $P_{0}$ (see \cite{Cornalba}).
According to \cite{Griffiths}, normal Weierstrass points form an open
subset of $\cW$. Therefore it is sufficient to analyze $\sigma^*\nabla^{\VV,(V_s)}$
around normal Weierstrass points.

Let $X$ be a fixed curve and let $P_0$ be a normal Weierstrass point on
$X$. The restriction of $\sigma^*\nabla^{\VV,(V_s)}$ to $X$ is then a connection on $X$ with
singularities at the Weierstrass points. Let us analyze the singularity
at $P_0$.

Let $z$ be a local coordinate at $P_0$. Let $P$ be a point of $X$,
close to $P_0$. Denote $\bar\omega_{(P)}(z,z')$ the bidifferential
form (\ref{bar:omega}) computed at point $P$. We have
$\bar\omega_{(P)}(z,z') = \wt \omega(z,z') + \sum_a \omega_a(z) d_z
\pa_{a}\ln\Theta(gP - z - \Delta)$. The second term of this sum
has a simple pole at $P = P_0$.  The value of the connection at the
point $P$ is
$$ \sigma^*\nabla^{\VV,(V_s)}_{\pa / \pa z} = {\pa\over{\pa z}} - \langle
T_{\bar\omega_P}(z) (dz)^2, \pa/\pa z \rangle,  
$$ 
which we express as 
$$ \sigma^*\nabla^{\VV,(V_s)}_{\pa / \pa z} = \pa_z + {{k\dimm\bar\G}\over{2\kappa}}
\langle \pa/\pa u, \sum_a \omega_a(u) dr_a^{(P)}(u) \rangle + O(1),
$$where $u$ is the local coordinate $z - z_P$ at $P$. Let us compute
the pole at $P = P_0$ of $\langle \pa/\pa u, \sum_a \omega_a(u)
dr_a^{(P)}(u) \rangle$.  For this, let us change the dual basis
$(\omega_a), (r_a^{(P)})$ to $(\omega^{(P)}_i), (r^{\prime(P)}_i)$,
with $\omega^{(P)}_i = u^{i-1}du + \sum_{j\geq 0} b_{ij}(P)
u^{g+i}du$, and $r^{\prime(P)}_i(u) = u^{-i} + O(1)$. We have then
$$ \langle \sum_a \omega_a(u) dr_a^{(P)}(u) \rangle = \sum_i \langle
\omega^{(P)}_i(u) dr_i^{\prime(P)}(u) \rangle = - g b_{g0}(z).
$$

Let us show that
$$ b_{g0}(z) = {1\over{gz}} + O(1).
$$ Let $(\beta_i)$ be the basis of the space of holomorphic differentials
on $X$, with expansions at $P_0$ 
$$ \beta_i(z) = z^{i-1}dz + \la_i z^{g-1}dz + \sum_{j \geq 0} \la_{ij}
z^{g+1+j}dz, \quad i = 1, \cdots, g-1, 
$$ 
$$ \beta_g(z) = z^{g}dz + \sum_{j \geq 0} \la_{gj}z^{g+1+j}dz.
$$ Let $x = z_P$.  We have 
$$
\pmatrix \beta_1 \\ \vdots \\ \beta_g
\endpmatrix = \pmatrix A(x) & \la(x) \\ c(x) & gx + o(x)
\endpmatrix \pmatrix du \\ \vdots \\ u^{g-1} du \endpmatrix + O(u^g)du,
$$ with $A(x)\in M_{(g-1) \times (g-1)}(\CC)$, $A(x) = id + O(x)$,
$\la(x)\in M_{(g-1)\times 1}(\CC)$, $\la(x) = {}^{t}(\la_1,
\ldots,\la_{g-1}) + O(x)$, $c(x)\in M_{1 \times(g-1)}(\CC)$, $c(x) =
(O(x^g),\ldots,O(x^2))$.

It follows that 
$$ \omega_i^{(P)} = \beta_i - {{\la_i} \over{gx}} \beta_g + \sum_{j =
  1}^g \al_{ij}(x) \beta_j, 
$$
$i = 1,\ldots,g-1$, $\al_{i,j}(x) = O(1)$, and 
$$ \omega_g^{(P)} = ({1\over{gx}} + O(1)) \beta_g + \sum_{j=1}^{g-1}
\al_{gj}(x) \beta_j, 
$$ $\al_{g,j}(x) = O(x)$. It follows that $b_{g0}(x) = {1\over{gx}} +
O(1)$.

We therefore find the expansion 
$$
\sigma^*\nabla^{\VV,(V_s)}_{\pa / \pa z} = \pa_z - {{k\dimm\bar\G}\over{2\kappa}}
{{dz}\over z} id_{\cB_{\VV,(V_s)}} + O(1). 
$$
\hfill \qed \medskip

\newpage

\section{Green functions and connections on $\Omega_X$} 
\label{sect:2}

\subsection{Theta-functions} \label{sect:theta}

For $\tau = (\tau_{ab})_{a,b=1,\ldots,\leq g}$ a symmetric matrix with negative definite 
real part, and $\la$ in $\CC^g$, let us set
$$ \Theta(\la | \tau) = \sum_{\mmm\in\ZZ^g} \exp({1\over 2} \mmm \tau
\mmm^t + \mmm \la^t).
$$ Let us denote by $\delta_j$ the $j$th canonical basis vector of $\CC^g$. We
have
$$ \Theta(\la + 2i\pi \delta_j | \tau) = \Theta(\la |\tau), 
\quad \Theta(\la +
\tau\delta_j|\tau) = e^{- {1\over 2} \tau_{jj} - z_j} \Theta(\la|\tau).
$$

Let $\cM_g^{(a,b)}$ the moduli space of systems $(X,A_a,B_a)$ of a
complex curve $X$  and of a canonical basis $(A_a,B_a)$ of $H_1(X,\ZZ)$.
  Let $m = $ class of $(X,A_a,B_a)$ be a point of the moduli space $\cM_g^{(a,b)}$.  Let
$(\omega_a)_{a=1,\ldots,g}$ be the basis of holomorphic
differentials on $X$, such that $\int_{A_a} \omega_b = 2i\pi
\delta_{ab}$, and let $\tau(m)$ be the matrix of periods defined by 
$\tau(m)_{ab} = \int_{B_a} \omega_b$.  In what follows, we will write
$\Theta(\la | m)$ or $\Theta(\la)$ instead of  $\Theta(\la | \tau(m))$. 
We also write $\pa_a \Theta( \la| m) = {{\pa \Theta}\over{\pa \la_a}}(\la|m)$,
$\pa_a \ln\Theta(\la|m) = {{\pa \ln\Theta}\over{\pa \la_a}}(\la|m)$, etc.

The degree zero part $J^0(X)$ of the Jacobian of $X$ is the quotient
$\CC^g / L$, where $L$ is the lattice $(2i\pi\ZZ)^g + \tau(m) \ZZ^g$. The 
Abel-Jacobi map sends a degree zero divisor $\sum_i n_i P_i$  to the class of
the vector  $(\sum_i \int_{x_0}^{P_i} \omega_a)_{a = 1, \ldots,g}$ 
(this map is independent of the choice of a  point $x_0$ of $X$). We
will identify classes of degree zero divisors with their images by the  
Abel-Jacobi map.  

The vector of Riemann constants $\Delta$ is then defined as the element
$\Delta$ of the degree $g-1$ part of the Jacobian $J^{g-1}(X)$, such that 
the identity $\Theta(x_1 + \ldots + x_{g-1} - \Delta) = 0$ holds identically, 
for $x_1,\ldots,x_{g-1}$ points of $X$. 

The function $z\mapsto \Theta(gP_0 - z - \Delta)$ is not identically
zero, because $P_0$ is not a Weierstrass point of $X$. We set $$ r_a(z)
= - \pa_{a} \ln \Theta(gP_0 - z - \Delta); $$ then $r_a$ satisfies
$\on{val}_{P_0}(r_a) \geq -g$, and $$ 
r_a(\gamma_{A_b}z) = r_a(z), \quad
r_a(\gamma_{B_b}z) = r_a(z) -
\delta_{ab}, $$
where $\gamma_{A_a}$ and $\gamma_{B_a}$ are deck transformations 
of the universal cover of $X$, corresponding to the cycles $A_a$ and $B_a$. 

\begin{remark} In the case where $X$ is an elliptic curve $\CC / (\ZZ +
\tau'\ZZ)$,  $\tau$ is the $1\times 1$ matrix equal to $2i\pi \tau'$.
The vector of Riemann constants is $\Delta = (1 + \tau') / 2$. 
\end{remark}

\subsection{Green functions}

Let $\omega^{out}_i$ be the forms such that $(\omega_a,\omega^{out}_i)_{a
  = 1,\ldots,g, i\geq 1} $ and $(r_a,t^i)_{a = 1,\ldots,g, i\geq 1}$
are dual bases of $\Omega^{out}$ and $(\oplus_a \CC r_a) \oplus
t\CC[[t]]$. Set
$$ 
G(z,z')_{z' \ll z} = \sum_{i\geq 1} \omega_i^{out}(z) (z')^i, \quad 
\bar G(z,z')_{z' \ll z} = \sum_{i\geq 1}
\omega_i^{out}(z) (z')^i + \sum_{a = 1,\ldots,g} \omega_{a}(z) r_{a}(z'),
$$ then
$G(z,z')_{z' \ll z}$ and $\bar G(z,z')_{z' \ll z}$ are the expansions for 
$z'$ near $P_0$ of the rational forms
\begin{align} \label{expr:G}
  & G(z,z') = d_z \ln\Theta(z' - z + (g-1)P_0 - \Delta) - d_z
  \ln\Theta(gP_0 - z - \Delta) \nonumber \\ & = - \sum_a \omega_a(z) [\pa_{a}
  \ln\Theta(z'-z+ (g-1)P_0 - \Delta) - \pa_{a}
  \ln\Theta(gP_0 - z - \Delta)]
\end{align}
and 
\begin{align} \label{bar:G}
  & \bar G(z,z') = d_z \ln\Theta(z' - z + (g-1)P_0 - \Delta)- d_z
  \ln\Theta( gP_0 - z - \Delta) + \sum_a \omega_a(z) r_a(z') \\ & 
  \nonumber =
  \sum_a \omega_a(z) [ - \pa_{a} \ln\Theta(z'-z+ (g-1)P_0 -
  \Delta) + \pa_{a}\ln\Theta(gP_0 - z - \Delta) 
  \\ & \nonumber - \pa_{a}\ln\Theta(gP_0 - z' - \Delta)] . 
\end{align}
$G(z,z')_{z' \ll z}$ and $\bar G(z,z')_{z' \ll z}$ belong therefore to 
$\CC((z))((z'))dz$. We will denote the expansions of $G(z,z')$ and $\bar
G(z,z')$ for $z'$ near $P_0$ by  $G(z,z')_{z' \gg z}$ and $\bar
G(z,z')_{z' \gg z}$ (they are elements of $\CC((z'))((z))dz$).  
We will use the same convention for the functions $G_\la(z,w)$. 

For $\la = (\la_a)_{a = 1, \ldots, g}$ in $\CC^g$, let $\cL_\la$ be the 
 line bundle over $X$, defined by $\Gamma(U,\cL_\la) = \{ f : \wt U \to
X | f(\gamma_{A_a}z) = f(z),   f(\gamma_{B_a}z) = e^{\la_a}f(z) , a = 1,
\ldots, g\}$,  where $U$ is any open subset of $X$ and    $\wt U$ is the
preimage of $U$ by the projection $\sigma:\wt X\to X$.  

Let $\Omega_X$ be the canonical line bundle of $X$. We identify elements of
$H^0(X - \{P_0\},\Omega_X \otimes\cL_\la)$ with their Laurent expansions at $P_0$ in
$\CC((t))dt$. $H^0(X,\Omega_X \otimes \cL_\la)$ is then spanned by a basis
$(\omega^{out}_{i,\la})_{i\geq 1-g}$, dual to the family $(t^i)_{i\geq
1-g}$ of $\CC((t))$.

$G_\la(z,w)$ is the twisted Green function defined by
\begin{align} \label{twisted:Green}
  & G_\la(z,w)   = \nonumber \\ & {{\Theta(z-w + (g-1) P_0 - \la - \Delta)}\over{\Theta(z-w+
      (g-1) P_0 - \Delta)\Theta((g-1)P_0 - \la - \Delta)}}
  \sum_{i=1}^g {{\pa \Theta}\over{\pa \la_a}}((g-1)P_0 - \Delta)
  \omega_a(z) .
\end{align}
We have then 
$$
G_\la(z,w)_{w\ll z} =  \sum_{i\geq 1-g} \omega^{out}_{i,\la}(z)w^i . 
$$

\subsubsection{Behaviour of the Green functions}

Let us emphasize the dependence of $G(z,w),\bar G(z,w)$
and $G_\la(z,w)$  in $P_0$ by writing them $G_{P_0}(z,w),\bar
G_{P_0}(z,w)$ and $G_{\la, P_0}(z,w)$.  Let us summarize the 
main functional properties for the Green functions: 

\begin{prop} 
$G_{P_0}(z,w)$ is defined for any $P_0$; it is regular in $z$, except for 
simples poles with residues $1$ and $-1$ when $z$ meets
$w$ and $P_0$. It is multivalued in $w$, vanishes for $w = P_0$, 
and single-valued in $z$.

$\bar G_{P_0}(z,w)$ is defined for any $P_0$ which is not a Weierstrass 
point. It is single-valued both in $z$ and $w$.  It has
simples poles with residues $1$ and $-1$ when $z$ meets
$w$ and $P_0$. As a function of $w$, it regular except for a
simple pole when $w$ meets $z$ and a pole at $P_0$
of order $g$ (when $P_0$ is Weierstrass). 

$G_{\la,P_0}(z,w)$ is defined for $\la$ generic ($(g-1)P_0 - \la - \Delta$
should not lie in the theta-divisor). It is multivalued both in 
$z$ and $w$. It is a regular form in $z$, except for a simple pole at $z = w$; 
it vanishes to 
order $g-1$ when $z$ meets $P_0$. Also it has a pole of order
$g-1$ when $w$ meets $P_0$. 
\end{prop}

{\em Proof.} Let $\ddelta$ be a nonsingular odd theta-characteristic.
Let us set $\wt \Theta  = \Theta[\ddelta]$ (see
(\ref{def:theta:delta})). We have then  $$  G_{P_0}(z,w) = d_z \ln
\wt\Theta(z-w)   - d_z \ln \wt\Theta(z-P_0)  .   $$    This proves that
$G_{P_0}(z,w)$ is well-defined for any $P_0$ in $X$. 

Let $Q_1, \ldots , Q_{g-1}$ be the zeroes (other than $P_0$) of the
holomorphic differential with $g-1$ zeroes at $P_0$. Both
$\Theta(z-w+(g-1)P_0 - \Delta)$ and  $\sum_a \pa_a \Theta((g-1)P_0 -
\Delta) \omega_a(z)$ vanish when $z$ equals $Q_i$, so that
$G_{\la,P_0}(z,w)$ has no pole when $z$ meets some  $Q_i$. 

The other points are direct checks. \hfill \qed\medskip

\begin{remark} 
  The functions $r_a(z)$ and the form $\bar G(z,z')$ only depend on a
  choice of a basis of $a$-cycles. When the basis of $a$-cycles is
  changed, $\bar G(z,z')$ gets transformed by the addition of a form
  in $z$. On the other hand, recall that we defined 
  $\bar\omega(z;z')$ by (\ref{bar:omega}). We have 
  $\bar\omega(z,z') = d_{z'}\bar G(z,z')$,
  because $(\oplus_a \CC r_a \oplus t\CC[[t]]) \oplus \CC 1 = 
  \left( \oplus_{i=1}^{g} \CC t^{-i} \oplus t \CC[[t]]\right) 
  \oplus \CC 1$. It follows that $\bar\omega(z,z')$ independent 
  of the choice of $a$-cycles (it is ``modular invariant'').
\end{remark}

\subsection{Connections on $\Omega_X$}

Let $\bar D$ and $D^{(\la)}$ be the connections defined over $\Omega_X$ by the rule
$$ (\bar D \omega)(z) = -\limm_{z'\to z} [\bar G(z,z') \omega(z') +
\bar G(z',z)\omega(z)]
$$
and
$$
(D^{(\la)}\omega)(z) = - \limm_{z'\to z}\left( G_\la(z,z')\omega(z') 
+ G_\la(z',z)\omega(z) \right) . 
$$


\begin{prop} \label{pestalozzi}
 1)  If $\omega$ is a one-form defined on an open subset $U$ of $X$, on
 which a coordinate $z_U$ is fixed, and we set $\omega_U(z) =
 \omega(z) / dz_U$, we have
\begin{align} \label{explicit:barD}
  & \bar D\omega(z) = \omega_U'(z) (dz_U)^2 \\ & \nonumber -
  {{\sum_a \pa_a \Theta((g-1)P_0 - \Delta) (\omega_a)_U'(z) (dz_U)^2 -
      \sum_{a,b} \pa^2_{ab} \Theta((g-1)P_0 - \Delta )
      \omega_a\omega_b(z) }\over{ \sum_a \omega_a(z) \pa_a
      \Theta((g-1)P_0 - \Delta)}} \omega(z) .
 \end{align} 
 Let $R_0$ be a point of $X$ and set
 $$ \bar\al_{R_0}(z) = {{\Theta(R_0 - z + (g-1)P_0 -
     \Delta)}\over{\Theta(z - R_0 + (g-1)P_0 - \Delta)}} (\sum_a
 \omega_a(z) \pa_a\Theta((g-1)P_0 - \Delta) ) .
 $$ The forms ${\bar\al}_{R_0}(z)$ are well defined and nonzero when
 $P_0$ is generic, and are all proportional. They have no poles and
 no zeroes except for a zero of order $2(g-1)$ at $P_0$ and we have,
 when $\bar\al_{R_0}$ is nonzero
 \begin{equation} \label{id:barD} 
  \bar D\omega = \bar\al_{R_0} d({\omega \over{\bar\al}_{R_0}}).
 \end{equation}

2) We have 
$$ D^{(\la)} \omega = \bar D\omega - 2 [\sum_a \omega_a(z) \pa_a
\ln\Theta(-\la + (g-1)P_0 - \Delta)] \omega. 
$$
\end{prop}

{\em Proof.} Let us prove 1). 
Let us set
$$
\wt{\wt G}(z,w) = d_z \ln \Theta(w - z + (g-1)P_0 - \Delta) . 
$$
Since $(\wt{\wt G} - \bar G)(z,z) = 0$, we have 
$$
(\bar D\omega)(z) = \limm_{z'\to z} [\wt{\wt G}(z,z')\omega(z') 
+ \wt{\wt G}(z',z)\omega(z) ]. 
$$ (\ref{explicit:barD}) is then a consequence of the expansion
\begin{align*}
  & \wt{\wt G}(z,w) = {{dz_U}\over{z_U - w_U}} \\ & + {1\over 2} {{
      \sum_a (\omega_a)_U'(z) \pa_a \Theta((g-1)P_0 - \Delta)(dz_U)^2
      - \sum_{a,b} \omega_a\omega_b(z) \pa^2_{ab} \Theta((g-1)P-0 -
      \Delta)}\over{\sum_a \omega_a(z) \pa_a \Theta((g-1)P_0 -
      \Delta)}} \\ & + o(z_U-w_U).
\end{align*} 
The fact that the forms $\bar\al_{R_0}$ are proportional follows
from the identity (\ref{id:separation}).

Let us prove (\ref{id:barD}). This equality is equivalent to 
\begin{align*} 
  & d_z \ln (\bar\al_{R_0})_U(z) \\ & = {{ \sum_a (\omega_a)_U'(z)
      \pa_a \Theta((g-1)P_0 - \Delta)(dz_U)^2 - \sum_{a,b}
      \omega_a\omega_b(z) \pa^2_{ab} \Theta((g-1)P_0 -
      \Delta)}\over{\sum_a \omega_a(z) \pa_a \Theta((g-1)P_0 -
      \Delta)}},
\end{align*}
that is 
\begin{equation} \label{to:prove} 
 d_z \ln {{\Theta(z - R_0 + (g-1)P_0 - \Delta) }\over{\Theta(R_0 - z
    + (g-1)P_0 - \Delta)}} = {{\sum_{a,b} \pa^2_{ab}\Theta((g-1)P_0 -
    \Delta) \omega_a\omega_b(z)}\over{\sum_a \pa_a \Theta((g-1)P_0 -
    \Delta)\omega_a(z)}} .
\end{equation}
Let us show (\ref{to:prove}). As we have seen, its left side is
independent of $R_0$. Let us compute the expansion of this left side
when $P_0$ tends to $z$.  We find
\begin{align*} 
  & [d_z\ln\Theta(z - R_0 + (g-1)P_0 - \Delta)]_U = {1\over{z - R_0}}
  \\ & + {1\over 2} {{\sum_a (\omega_a)_U'(z) (dz_U)^2
      \pa_a\Theta((g-1)P_0 - \Delta) + \sum_{a,b} \omega_a\omega_b(z)
      \pa^2_{ab}\Theta((g-1)P_0 - \Delta) }\over{\sum_a \omega_a(z)
      \pa_a\Theta((g-1)P_0 - \Delta)}} \\ & + o(z_U - (R_0)_U)
\end{align*}
and
\begin{align*} 
  & [d_z\ln\Theta(R_0 - z + (g-1)P_0 - \Delta)]_U = {1\over{z - R_0}}
  \\ & + {1\over 2} {{\sum_a (\omega_a)_U'(z) (dz_U)^2
      \pa_a\Theta((g-1)P_0 - \Delta) - \sum_{a,b} \omega_a\omega_b(z)
      \pa^2_{ab}\Theta((g-1)P_0 - \Delta) }\over{\sum_a \omega_a(z)
      \pa_a\Theta((g-1)P_0 - \Delta)}} \\ & + o(z_U - (R_0)_U).
\end{align*}
The limit of the left side of (\ref{to:prove}) when $P_0$ tends to $z$
is therefore equal to its right side. (\ref{to:prove}) then follows
from the fact that its left side is independent of $P_0$.

Let us show 3). Let us set
$$
\wt G(z,w) = d_z \ln\Theta(z-w+(g-1)P_0 - \Delta), 
$$
then we have 
\begin{align*} 
 & (G_\la - \wt G)(z,w)_{| z = w} = \sum_a \omega_a(z) \pa_a
\ln\Theta( -\la + (g-1)P_0 - \Delta ) \\ & - {{\sum_{a,b} \omega_a
    \omega_b(z) \pa^2_{ab} \Theta((g-1)P_0 - \Delta)}\over{\sum_a
    \omega_a(z) \pa_a \Theta((g-1)P_0 - \Delta)}};
\end{align*}
therefore, if we set 
$$ \wt D \omega = -\limm_{z'\to z}(\wt G(z,z') \omega(z') + \wt
G(z',z) \omega(z)),
$$
we get 
\begin{align} \label{saboche}
  D^{(\la)}\omega = & \wt D\omega + 2 [ - \sum_a \omega_a(z)
  \pa_a\ln\Theta( -\la + (g-1)P_0 - \Delta ) \\ & \nonumber +
  {{\sum_{a,b} \omega_a \omega_b(z) \pa^2_{ab} \Theta((g-1)P_0 -
      \Delta)}\over{\sum_a \omega_a(z) \pa_a \Theta((g-1)P_0 -
      \Delta)}}] \omega.
\end{align}
 A proof similar to that of Prop.\ \ref{pestalozzi}, 1) shows that
 \begin{align} \label{explicit:wtD}
   & \wt D\omega(z) = \omega_U'(z) (dz_U)^2 \\ & \nonumber - {{\sum_a
       \pa_a \Theta((g-1)P_0 - \Delta) (\omega_a)_U'(z) (dz_U)^2 +
       \sum_{a,b} \pa^2_{ab} \Theta((g-1)P_0 - \Delta )
       \omega_a\omega_b(z) }\over{ \sum_a \omega_a(z)
       \pa_a \Theta((g-1)P_0 - \Delta)}} \omega(z) .
 \end{align} 
 3) follows then from (\ref{saboche}), (\ref{explicit:barD}) and
 (\ref{explicit:wtD}).  \hfill \qed \medskip

\begin{remark} 
  The proof of Prop.\ \ref{pestalozzi}, 1) also shows that if $R_0$ is
  a point of $X$ and
  $$ \wt\al_{R_0}(z) = {{\Theta(z - R_0 + (g-1)P_0 -
      \Delta)}\over{\Theta(R_0 - z + (g-1)P_0 - \Delta)}} (\sum_a
  \omega_a(z) \pa_a\Theta((g-1)P_0 - \Delta) ) .
  $$ The forms $\wt\al_{R_0}(z)$ are well-defined and nonzero when
  $R_0$ is generic, and are all proportional. They have no poles and
  no zeroes except for double zeroes at the points $Q_1, \ldots,
  Q_{g-1}$ and we have, when $\wt \al_{R_0}$ is nonzero
 \begin{equation} \label{id:wtD} 
   \wt D\omega = \wt\al_{R_0} d({\omega \over{\wt\al_{R_0}}}).
 \end{equation}
\end{remark}

\begin{remark}
  Another proof of (\ref{id:wtD}) can be found in \cite{Fay},
  Prop.\ 2.2. 
\end{remark}

\newpage

\section{Functional expression of the KZB connection}
\label{sect:3}

In this section, we will assume that $\bar \G = \SL_2$. 

\subsection{The modules} \label{sect:modules}

Let $e,f,h$ be the Cartan basis of $\bar\G$ and $\bar\HH = \CC h$, 
$\bar\N_+ = \CC e$  and $\bar\N_- = \CC f$.  Let $\ell$ be an
integer such that $\ell \geq g-1$.  Let us set
$$
\G_+^{(\ell)} = \left( \bar\HH \otimes \CC[[t]] \right) \oplus  
\left(  \bar\N_+ \otimes t^\ell\CC[[t]] \right) \oplus 
\left( \bar\N_- \otimes t^{-\ell}\CC[[t]] \right) \oplus \CC K;
$$  
$\G_+^{(\ell)}$ is a Lie subalgebra of $\G$. Let  $\chi$ be the
character of $\G_+^{(\ell)}$ defined by  $\chi(K) = k$, $\chi(h[1]) =
-\ell k$, and $\chi$ is zero when restricted to  $\bar\HH \otimes t
\CC[[t]] \oplus   \bar\N_+ \otimes t^\ell\CC[[t]] \oplus \bar\N_-
\otimes t^{-\ell}\CC[[t]]$ (recall that $h[1] = (h\otimes 1,0)$ under
the isomorphism  $\G = \bar\G\otimes\CC((t)) \oplus\CC K$).  Let
$\VV_{-\ell k}$ be the induced module $U\G\otimes_{U\G_+^{(\ell)}}
\CC_\chi$,  and let $v$ be the vector $1\otimes 1$ of this module.  
$\VV_{-\ell k}$ belongs to the class (C) of $\G$-modules defined in sect.\
\ref{class:C}.  When $k$ is a positive integer and $\ell$ is even,
$\VV_{-\ell k}$ has a quotient isomorphic to the vacuum integrable
module of level $k$.

Denote by $(\rho_{-\La},V_{-\La})$ the lowest weight Verma module over
$\bar\G$, freely generated by the vector $v_{-\La}$ such that
$\rho_{-\La}(f) (v_{-\La}) = 0$ and $\rho_{-\La}(h)(v_{-\La}) = - \La
v_{-\La}$.

In what follows, we will set $\VV = \VV_{-\ell k}$ and $V_s =
V_{-\La_s}$, and denote $\cB_{\VV,(V_s)}$ by $\cB_{\ell,(\La_i)}$.  In the
notation of sect.\ \ref{sect:before}, $\cV$ is then the
$\bar\G(R)$-module  $\VV_{-\ell k} \otimes (\otimes_i
V_{-\La_i}^{(P_i)})$; the spaces $\cB_{\ell,(\La_i)} = (\cV^*)^{\bar\G(R)}$
form a sheaf over $\cM_{g,n+1}$. Moreover,  the lift of
$\cB_{\ell,(\La_i)}$ to $\cM_{g,1^\infty, n \cdot 1^2}$ is endowed with the
 flat connection $\nabla^{\cB_{\ell,(\La_i)}}$.

\subsection{The connection $(\wt\cB_{\ell,(\La_i)},\wt\nabla^{\ell,(\La_i)})$}
\label{bunuel}

Let $\cM_g^{(a)}$ be the moduli space of genus $g$ curves with marked
$a$-cycles and $\cM^{(a)}_{g,1^\infty,n \cdot 1^2}$ be the fibered
product $\cM_{g}^{(a)} \times_{\cM_g} \cM_{g,1^\infty,n\cdot 1^2}$.

Let $\pi^{(a)}$ be the canonical projection from
$\cM^{(a)}_{g,1^{\infty}, n \cdot 1^2}$ to $\cM_{g,1^{\infty}, n \cdot
  1^2}$. This projection is a local isomorphism. Let 
$(\cB^{(a)}_{\ell,(\La_i)}, \nabla^{\ell,(\La_i), (a) })$ be the lift of
$(\cB_{\ell,(\La_i)}, \nabla^{\cB_{\ell,(\La_i)}})$ to
$\cM^{(a)}_{g,1^{\infty}, n \cdot 1^2}$.

For $\la = (\la_a)_{a = 1, \ldots,g}$ a $g$-uple of complex numbers, let us denote by 
$\G^{out}_\la$ the subalgebra of $\G$ formed by the Laurent expansions at a point 
of $\sigma^{-1}(P_0)$ of the functions $x$ on $\wt{X} - \sigma^{-1}(P_0)$, such that 
$$
x(\gamma_{A_a}(z))  = x(z) , \quad 
x(\gamma_{B_a}(z))  = \Ad(e^{\la_a h}) (x(z)) ; 
$$
recall that $\gamma_{A_a}$  and $\gamma_{B_a}$ are the deck transformations 
associated to the $a$- and $b$-cycles of $\wt X$.  

Let $\cO_{P_i}$ be the local ring of $X$ at $P_i$. Let $\G_{reg\ at\
P_i}$ be the extension of the Lie algebra $\bar\G \otimes [\CC((t))
\oplus (\oplus_i \cO_i)]$ by the Kac-Moody cocycle.  There is a unique Lie
algebra morphism $\ev$ from  $\G_{reg\ at\ P_i}$ to $\G\oplus (\oplus_i
\bar\G)$, induced by  evaluation at points $P_i$. 
Let us denote by $\pi_{\VV \otimes (\otimes_i V_{-\La_i})}$ the representation of
$\G_{reg\ at\ P_i}$ in $\VV \otimes (\otimes_i V_{-\La_i})$ equal to 
$\pi_{\VV} \otimes (\otimes_i \pi_{V_{-\La_i}}^{(P_i)})$. 

In what follows, we will view $\G^{out}_\la$ as a Lie subalgebra of 
$\G_{reg\ at\ P_i}$. 

Let $\psi(m_\infty,P_i,v_i)$ be a flat section of $(\cB_{(a),\ell,(\La_i)}, 
\nabla^{(a),\ell,(\La_i)})$ and  let $(\la_a^f)_{a = 1, \ldots,g}$ be formal
variables. Let us set  $\psi(m_\infty,P_i,v_i | \la_a^f) =
\psi(m_\infty,P_i,v_i)\circ  \pi_{\VV \otimes (\otimes_i V_{-\La_i}) }
(e^{\sum_a \la_a^f h[r_a]})$. Then 
$\psi(m_\infty,P_i,v_i | \la_a^f)$ is a $\G^{out}_{(\la_a^f)}$-invariant on 
$\VV \otimes (\otimes_i V_{-\La_i})$, and
is a solution of equations
$$ 
\pa_{\la_a^f} \psi(m_\infty,P_i,v_i | \la_a^f)  
= \psi(m_\infty,P_i,v_i | \la_a^f) \circ \pi_{\VV \otimes (\otimes_i
V_{-\La_i}) }(h[r_a]) , 
$$ 
$$ 
\pa_{[\xi]} \psi(m_\infty,P_i,v_i | \la_a^f)   = 
\psi(m_\infty,P_i,v_i | \la_a^f)   \circ T_{\bar\omega}[\xi]^{(0)} . 
$$ 

We will therefore define $\wt\cB_{\ell,(\La_i)}$ as the sheaf over 
$\cM^{(a)}_{g,1^{\infty},n\cdot 1^2}$ whose fibre over $(m_\infty,P_i,v_i)$
is the space of maps $\phi$ from an open subset of $\CC^g$ to 
$[\VV \otimes(\otimes_i V_{-\La_i})]^*$, such that for any $\la$
in $\CC^g$, $\phi(\la)$ is $\G^{out}_\la$-invariant, and 
satisfies the differential equation
$$
\pa_{\la_a} \phi(\la) 
= \phi(\la) \circ \pi_{\VV \otimes (\otimes_i
V_{-\La_i}) }(h[r_a]) . 
$$
We define also a connection $\wt\nabla^{\ell,(\La_i)}$ over $\wt\cB_{\ell,(\La_i)}$  
by 
$$
\wt\nabla^{\ell,(\La_i)}_{[\xi]} \psi(m_\infty,P_i,v_i | \la_a)   =  
\pa_{[\xi]} \psi(m_\infty,P_i,v_i | \la_a)  
- \psi(m_\infty,P_i,v_i | \la_a)   \circ T_{\bar\omega}[\xi]^{(0)} .  
$$
It is then easy to see that $\wt\nabla^{\ell,(\La_i)}$ is flat. 
A flat section of $(\wt\cB_{\ell,(\La_i)},\wt\nabla^{\ell,(\La_i)})$ 
therefore consists of a function $\psi(m_\infty,P_i,v_i | \la)$ 
defined on an open subset of $\cM_{g,1^\infty,n,\cdot 1^2}$, such that 
\begin{equation} \label{nagila1}
\pa_{\la_a} \psi(m_\infty,P_i,v_i | \la)  
= \psi(m_\infty,P_i,v_i | \la) \circ \pi_{\VV \otimes (\otimes_i
V_{-\La_i}) }(h[r_a])  
\end{equation}
and 
\begin{equation} \label{nagila2}
\pa_{[\xi]} \psi(m_\infty,P_i,v_i | \la_a)   = 
\psi(m_\infty,P_i,v_i | \la_a)   \circ T_{\bar\omega}[\xi]^{(0)} . 
\end{equation}

\subsection{The functional connection $(\cF_{\ell,(\La_i)},
\nabla^{\cF_{\ell,(\La_i)}})$} \label{corr}

Set 
$$
N = {1\over 2} [k\ell + \sum_{i}\La_i] .
$$ 
Define a sheaf $\cF_{\ell,(\La_i)}$ over
$\cM_{g,1^\infty, n \cdot 1^2}^{(a)}$ as follows. 
The fiber of $\cF_{\ell,(\La_i)}$ 
over $(m_\infty,P_i,v_i)$ is the space of all forms 
$f(\la|z_1,\ldots,z_N)$
defined on $U\times (\wt{X} - \sigma^{-1}\{P_0,P_i\})^N$, 
where $U$ is an open subset of $\CC^g$, with the conditions that 
$f(\la|z_1,\ldots,z_N)$ is a regular one-form in each variable $z_i$,
has at most simple poles at the $P_i$ ($i\neq 0$), satisfies the
transformation properties
$$
(\gamma_{A_a}^{(i)})^*
f(\la|z_1,\ldots,z_N) = 
f(\la|z_1,\ldots,z_N), \quad 
(\gamma_{B_a}^{(i)})^* 
f(\la|z_1,\ldots,z_N) = 
e^{-2\la_a} f(\la|z_1,\ldots,z_N), 
$$ 
and is symmetric under permutation of the variables $z_i$. Here 
$\gamma_{A_a}^{(i)}$ and $\gamma_{B_a}^{(i)}$ denote the action of
$\gamma_{A_a}$ and $\gamma_{B_a}$ on the $i$th variable $z_i$. 

For $f(m_\infty,P_i,v_i|\la,z_\al)$ a section of $\cF_{\ell,(\La_i)}$, set 
\begin{align} \label{operator:KZ:moduli}  
  \nabla^{\cF_{\ell,(\La_i)}}_{[\xi]} f(m_\infty,P_i,v_i|\la,z_\al) = 
  (\pa_{[\xi]}f)(m_\infty,P_i,v_i|\la,z_\al) - (K_{\xi}f)(m_\infty,P_i,v_i|\la,z_\al),
\end{align} 
where $[\xi]$ is the variation of moduli induced by a formal vector
field $\xi$ at $P_0$, and
\begin{equation} \label{KZB:moduli}  
  (K_{\xi}f)(m_\infty,P_i,v_i|\la,z_\al) = \res_{z = P_0}(\xi(z)
  \left(T^{\cF}_{\bar\omega}(z)f\right)(m_\infty,P_i,v_i|\la,z_\al) ),
\end{equation}
where 
\begin{align} \label{expr:conn:cF}
  & 2\kappa \left(T^{\cF}_{\bar\omega}(z)f\right)(m_\infty,P_i,v_i|\la,z_\al) = \\ & = \left[
    {1\over 2} \left(\sum_a \omega_a(z) \pa_{\la_a} + 2 \sum_{\al}
      G(z,z_{\al}) - \sum_i \La_i G(z,P_i)\right)^2 \right.\\ & 
   \nonumber \left.
    + \sum_a D^{(2\la)}_{z} \omega_a(z) \pa_{\la_a} + 2\sum_{\al}
    (D^{(2\la)}_{z} \otimes 1) G(z,z_\al) - \sum_{i} \La_i
    (D^{(2\la)}_z \otimes 1) G(z,P_i) + k \bar\omega_{2\la}(z)
  \right] \\ \nonumber & f(m_\infty,P_i,v_i|\la,z_\al) \\ \nonumber & + \sum_{\al = 1}^N \left[ -
    2G_{2\la}(z,z_{\al}) \left( \sum_a \omega_a(z_\al) \pa_{\la_a} + 2
      \sum_{\beta \neq \al}G(z_\al,z_\beta) - \sum_{i}\La_i
      G(z_\al,P_i) \right) \right.  \\ \nonumber & \left.  + \left(- 4
      G_{2\la}(z,z_\al)G(z_\al,z) +2 k d_{z_\al}G_{2\la}(z,z_\al)
    \right) \right] f(m_\infty,P_i,v_i|\la,z,(z_\beta)_{\beta\neq\al}).
\end{align}
Here $\kappa = k+2$,
\begin{equation} \label{bar:omega:la}
  \bar\omega_{2\la}(z) = \limm_{{z'\to z}} 2 d_{z'}[G_{2\la}(z,z')-
  G(z,z') ] - 3 \sum_a \omega_a(z) dr_a(z) . 
\end{equation}

\begin{prop} (see \cite{EF}) $\nabla^{\cF_{\ell,(\La_i)}}$ defines a flat connection 
over $\cF_{\ell,(\La_i)}$. 
\end{prop}

\subsection{The correlation functions morphism}

For $\psi(m_\infty,P_i,v_i|\la)$ a section of $\wt B_{\ell,(\La_i)}$, let us set 
$$
[\corr(\psi)](m_\infty,P_i,v_i)
(z_1,\ldots,z_N|\la) = 
\langle \psi(m_\infty,P_i,v_i|\la), \pi_{\VV}(\prod_{i=1}^N e(z_i))
(v\otimes (\otimes_i v_{-\La_i}) ) \rangle,  
$$
where $e(z)$ is the generating series $\sum_{i\in\ZZ} e[t^i] z^{-i-1}
dz$. 

\begin{prop} \label{conchita} (see  \cite{EF}) 
$\corr$ defines a morphism of sheaves with connection from
$(\wt\cB^{(a)},\wt\nabla^{(a)})$ to $(\cF_{p,n}, \nabla^{\cF_{p,n}})$.  
\end{prop}

\begin{remark} \label{rem:cb}
In the case where $k$ is a positive integer, one may translate the condition 
that $f$ is in the image by $\corr$ of the subbundle of conformal blocks into 
functional conditions on $f$, see \cite{EF}.  
\end{remark}

\subsection{Properties of $\bar\omega_{\la}(z)$}

The proof of Prop.\ \ref{conchita} uses the following results, that we stated 
without proof in \cite{EF}.  
Let us set 
\begin{equation} \label{omega:la}
  \omega_\la(z) = \limm_{z'\to z} 2 d_{z'}(G_\la(z,z') - G(z,z')).
\end{equation}

\begin{lemma} \label{chinois}
  $\omega_\la(z)$ is a quadratic differential on $X$, regular outside
  $P_0$.  Its expansion at $P_0$ is
  \begin{equation} \label{expansion} 
    \omega_\la(z) = - g(g-1) {{(dz)^2}\over{z^2}} + O({1/z}).
 \end{equation}
\end{lemma}

{\em Proof.} As we have seen, $G_{\la}(z,w)$ and $G(z,w)$ are rational
forms defined on $X\times X$. The poles of $G_{\la}(z,w)$ are the
diagonal and a pole of order $g-1$ at $w = P_0$; the poles of $G(z,w)$
are the diagonal and $z = P_0$ [$\Theta(z' - z + (g-1)P_0 - \Delta)$ has a
zero at $z' = Q_1, \cdots, Q_{g-1}$, however, logarithmic derivative
with respect to $z$ gives no singularity at this point]. Moreover,
$(G_\la - G)(z,w)$ is regular at the diagonal, therefore $(G_\la -
G)(z,w)$ has its only singularities at $w = P_0$ and $z = P_0$. It
follows that $\omega_\la(z)$ is regular on $X - \{P_0\}$.

Let us prove (\ref{expansion}). Since $(G_\la - G)(z,w)$ has its only
singularities at $w = P_0$ and $z = P_0$, the difference of the
expansions of $G_\la(z,w)$ and $G(z,w)$ (which belong to
$\CC((z))((w)) dz$) belong to $\CC[[z,w]][z^{-1},w^{-1}]dz$. Applying
$d_w$ and setting $w = z$ in the latter expansion will therefore give
us the expansion at $P_0$ of $\omega_\la(z)$.

Let us write 
$$ G_\la(z,w) = \sum_{i = 1-g}^{0} \omega_{i,\la}^{out}(z) w^i +
\sum_{i\geq 1} (\omega_{i,\la}^{out} - \omega_{i}^{out}) w^i. 
$$
It follows from the identities $\langle \omega_{i,\la}^{out}, w^j \rangle
= \delta_{ij}$ if $j\geq 1-g$ and $\langle \omega_{i}^{out}, w^j \rangle
= \delta_{ij}$ if $j\geq 1$ that 
$$
\omega_{i,\la}^{out}(z) = {{dz}\over{z^{i+1}}} + O(z^{g-1})\ \on{if}\ 
 i\geq 1-g,  
\quad 
(\omega_{i,\la}^{out} - \omega_{i}^{out})(z) = O(z^{-1}) \ \on{if}
\ i\geq 1. 
$$ 
Therefore, we have 
$$ (G_\la - G)(z,w) = \sum_{i = 1-g}^{0} {{dz}\over{z^{i+1}}} (w^i +
\sum_{j\geq g-1} a_{ij} w^j) + \sum_{i\geq 1, j\geq -1} b_{ij} z^i w^j. 
$$
It follows that 
$$ [d_w(G_\la - G)(z,w)]_{|w = z} = (\sum_{i = 1-g}^0
i){{(dz)^2}\over{z^2}} + O(1/z) = -
{{g(g-1)}\over{2}}{{(dz)^2}\over{z^2}} + O(1/z).
$$ (\ref{expansion}) follows.  \hfill \qed \medskip

\begin{cor} \label{omega:la:local}
  $\bar\omega_{\la}(z)$ is a quadratic differential on $X$, regular
  outside $P_0$. Its expansion at $P_0$ is
  $$ \bar\omega_{\la}(z) = {{g(g+5)}\over{2}} {{(dz)^2}\over{z^2}} +
  O(1/z).
  $$
\end{cor}

{\em Proof.} As we have seen, we can find dual bases
$(\omega'_i)_{i = 1, \ldots,g}$ and $(r'_i)_{i = 1, \ldots,g}$ of
$\oplus_{a=1}^g \CC \omega_a$ and $\oplus_{a = 1}^g \CC r_a$, with
$\omega'_i(z) = z^{i-1}dz + O(z^i)$, and $r'_i(z) = z^{-i} +
O(z^{-i+1})$. It follows that $\sum_{a=1}^g \omega_a dr_a = \sum_{i =
  1}^g \omega'_a dr'_a = (\sum_{i=1}^g -i){{(dz)^2}\over{z^2}} +
O(1/z)$, so 
 
\begin{equation} \label{expansion:omega:dr}
  \sum_{a=1}^g \omega_a dr_a = - {{g(g+1)}\over 2}
  {{(dz)^2}\over{z^2}} + O(1/z).
\end{equation}

The Corollary now follows from (\ref{expansion:omega:dr}) and
Lemma \ref{chinois}.
\hfill \qed \medskip

\newpage

\section{Flat sections of the sheaf $\nabla^{\cF_{\ell,(\La_i)}}$ in the case $n = N = 0$}
\label{sect:N=0:n=0}

In this section, we will restrict ourselves to the case  $n = N = 0$ and
 write $\cF$ instead of $\cF_{\ell,(\La_i)}$. A section of the sheaf
$\cF$ is a function $f(m_\infty|\la)$, where $m_\infty$ belongs to
$\cM_{g,1^\infty}$ and $\la$ is a formal variable on $\CC^g$.

We have 
$$ (T^{\cF}_{\bar\omega}(z) f)(m_\infty|\la) = {1\over{2\kappa}}\left[{1\over
    2}(\sum_a \omega_a (z) \pa_{\la_a})^2 + (D^{(2\la)} \omega_a)(z)
  \pa_{\la_a}+ k \bar\omega_{2\la}(z) \right] f,
$$
and the flatness condition is written as 
\begin{equation} \label{eq:flatness} 
  \pa_{[\xi]} f(m_\infty|\la) = \langle \xi(z) ,
  (T^{\cF}_{\bar\omega}(z)f)(m_\infty|\la) \rangle_{z}.
\end{equation}

\begin{thm} \label{thm:no:point} 
  Let us set
 \begin{align} \label{def:S} 
   & S(m,P_0|z) = \\ & \nonumber
   {{\sum_{a,b,c}\omega_a\omega_b\omega_c(z)
       \pa^{3}_{abc}\Theta((g-1)P_0 - \Delta)}\over{ \sum_a
       \omega_a(z) \pa_{a}\Theta((g-1)P_0 - \Delta)}} - {{\sum_{a,b}
       \omega_a\omega'_b(z) \pa^2_{ab}\Theta((g-1)P_0 -
       \Delta)}\over{\sum_a \omega_a(z) \pa_{a}\Theta((g-1)P_0 -
       \Delta)}} \\ & \nonumber - {{\sum_{a,b} \omega_a\omega_b(z)
       \pa^2_{ab}\Theta((g-1)P_0 - \Delta) }\over{(\sum_a \omega_a(z)
       \pa_{a}\Theta((g-1)P_0 - \Delta))^2}} \cdot \\ & \nonumber
   \cdot [ \sum_{a,b} \omega_a\omega_b(z) \pa^2_{ab}\Theta((g-1)P_0 -
   \Delta) - \sum_a \omega'_a(z) \pa_{a}\Theta((g-1)P_0 - \Delta) ] .
 \end{align}
 and $\bar S(m,P_0|z) = S(m,P_0|z) - 3\sum_a \omega_a(z)
 dr^{(P_0)}_a(z)$. $\bar S$ is a quadratic differential on $X$,
 regular outside $P_0$, and with expansion at this point given by
 $$ \bar S(m,P_0|z) = {{g(g+5)}\over 2} {{(dz)^2}\over{z^2}} + O(1/z).
 $$ 

 There exists a nonzero function $\al(m,P_0,v_0)$, locally defined on
 $\cM_{g,1^2}$, homogeneous of degree ${{g(g+5)}\over{2}}$ in $v_0$,
 such that
  \begin{equation} \label{eq:for:alpha} 
    \pa_{[\xi]} \al(m,P_0,v_0) = \langle \bar S(m,P_0|z) , \xi(z) 
    \rangle_z \al(m,P_0,v_0),
  \end{equation} 
  where $[\xi]$ is the variation of moduli induced by the vector field
  $\xi\in \CC((t)){\pa\over{\pa t}}$.

  The space of local flat sections of $\nabla^{\cF}$ is the space of
  functions $f(m_\infty | \la)$ of the form
 \begin{equation} \label{fim}
   f(m_\infty|\la) = \bar f(m_2|\la),
 \end{equation} where $m_2$ is the
 projection of $m_\infty$ by the projection $\cM_{g,1^\infty} \to
 \cM_{g,1^2}$, and $\bar f(m_2|\la)$ is of the form
 \begin{align} \label{f}
   \bar f(m,P_0,v_0 | \la) = & \Theta(-2\la + (g-1)P_0 -
   \Delta|m)^{-1} \varphi(m | \kappa \la + \Delta - (g-1) P_0) \cdot
   \\ & \nonumber \cdot \al(m,P_0,v_0)^{{k\over{2\kappa}}} ,
 \end{align}
 where $\varphi$ is a function defined locally on $\cM_g \times \CC^g$
 such that
 \begin{equation} \label{eq:varphi} 
   \delta_{\delta m} \varphi(\la | m) = -{\kappa\over{8i\pi}}
   \sum_{a,b} \delta \tau_{ab} \pa^2_{\la_a\la_b} \varphi(\la | m) ,
  \end{equation} 
  where $\delta m$ is a variation of moduli and $\delta \tau_{ab}$ is
  the corresponding variation of periods.
  (Let $\Sigma_g$ be the space of Siegel matrices $\tau = (\tau_{ab})_{a,b
  = 1,\ldots, g}$,  that is the space of matrices $\tau$ which are
  symmetric and with positive imaginary part.  $\varphi(m|\la)$ may be of
  the form $\wt\varphi(\tau(m)|\la)$,  where $\wt\varphi$ is a function on
  $\Sigma_g\times \CC^g$, solution of the equation
  ${{\pa}\over{\pa\tau_{ab}}} \wt\varphi(\tau|\la) = 
  {{\kappa}\over{8i\pi}} \pa^2_{\la_a\la_b}\varphi(\tau|\la)$.)
\end{thm}

{\em Proof.} We first prove:

\begin{lemma} \label{var:moduli} 
  Let $\xi$ be a formal vector field defined around $P_0$.  The
  variation of periods induced by $\xi$ is
\begin{equation} \label{delta:tau}
  \delta \tau_{ab} = 2i\pi \res_{P_0}(\xi \omega_a \omega_b ).
\end{equation}
\end{lemma}

{\em Proof.} Let $\cK = \CC((z))$ and $\cO = \CC[[z]]$. Let $\delta m$
be the variation of moduli induced by $\xi$.  The complex structure
associated with $m + \delta m$ is defined by the inclusions of $\cO$
in $\cK$ (which plays the role of the local ring of the new structure)
and $R_\eps = (1 + \eps \cL_\xi)(R)$ in $\cK$ (which is the new ring
of regular functions). The space of formal differential forms at $P_0$
is $\Omega_\cK = \CC((z))dz$, the space of regular forms at $P_0$ is
$\Omega = \CC[[z]]dz$; denote by $\Omega_R$ the space of forms regular
outside $P_0$. Then we have $\Omega_\eps = (1 + \eps
\cL_\xi)\Omega_R$, where $\cL_\xi$ is the Lie derivative
$\cL_{\xi}(\omega) = d_z(\xi\omega)$.

The deformation of $\omega_a$ is then  $$ \omega'_a = (1 +
\eps\cL_\xi)(\omega_a) + \delta\omega_a,  $$ with $\delta\omega_a$
regular outside $P_0$ and with  vanishing $A_a$-periods, such that
$\eps\cL_\xi(\omega_a)  + \delta \omega_a$ is regular at $P_0$. Since
$\delta\omega_a$ has vanishing $A_a$-periods, we have  $\delta\omega_a =
- {1 \over{2i\pi}}\sum_b \delta\tau_{ab} dr_b + dr$, for some $r$ in $R$
(the $r_b$ are such that $r_b(\gamma_{B_a}z) = r_b(z) - 2 i \pi
\delta_{ab}$; they are dual to the $\omega_a$). Since $\omega'_a$
belongs to $\Omega_\cO$, its pairing with all functions of $\cO$
vanishes; pairing it with the $\int_{P_0}^z \omega_a$, we find
(\ref{delta:tau}).  \hfill \qed \medskip

Let $\delta$ be an odd theta-divisor of $X$ (see e.g.\ \cite{Fay}).
\begin{lemma} \label{ajax}
  For $a = 1,\ldots,g$, let $Q_1^{(a)},\ldots,Q_{2g-2}^{(a)}$ be the
  zeroes of $\omega_a$; we have
  $$\bar D \omega_a(z) = \left(\sum_{i=1}^{2g-2} d_z
    \ln\Theta(z-Q_i^{(a)} + \delta - \Delta) - 2(g-1) d_z
    \ln\Theta(z-P_0 + \delta - \Delta) \right) \omega_a(z).
  $$
\end{lemma}

{\em Proof.} It follows from Prop.\ \ref{pestalozzi} that when
$\bar\al_{R_0}$ is not identically zero,
$$ {{\bar D \omega_a}\over{\omega_a}} = {{d(\omega_a /
    \bar\al_{R_0} )}\over{\omega_a / \bar\al_{R_0} }}.
$$
On the other hand, we have 
$$ {{\omega_a} \over{\bar\al_{R_0}}}(z) = c(R_0) {{
    \prod_{i=1}^{2g-2} \Theta(z-Q_i^{(a)} + \delta -
    \Delta)}\over{\Theta(z - P_0 + \delta - \Delta)^{2(g-1)}}},
$$ where $c(R_0)$ is a constant, because both sides are meromorphic
functions on $\wt X$ with zeroes and poles of the same orders at the
same points, and identical transformation properties along cycles of
$X$ (because of the equality $\sum_i Q_i^{(a)} = 2 \Delta$ in
$J^{2(g-1)}(X)$).  \hfill \qed \medskip

\begin{lemma} \label{var:int}
  We have
  $$ \pa_{[\xi]} \int_{(g-1)P_0}^\Delta \omega_a = {1\over 2}
  \langle \bar D \omega_a, \xi \rangle . 
  $$ 
\end{lemma}

{\em Proof.} Since $2\Delta = \sum_i Q_i^{(a)}$, we have
$$ \int_{(g-1)P_0}^{\Delta} \omega_a = {1\over 2} \sum_{i}
\int_{P_0}^{Q_i^{(a)}} \omega_a.
$$
Now $$
 \pa_{[\xi]}(\int_{P_0}^{Q_i^{(a)}} \omega_a) 
 = \int_{P_0}^{Q_i^{(a)}} \delta_{[\xi]} \omega_a,  
 $$ because $\omega_a$ vanishes at each $Q_i^{(a)}$.  It follows from
 the proof of Lemma \ref{var:moduli} that
$$
\delta_{[\xi]} \omega_a = d(\xi\omega_a) - [d(\xi\omega_a)]_{out},  
$$ where $\al\mapsto \al_{out}$ is the projection on the second
component of the direct sum decomposition $\CC((t)) dt = K \oplus
t^{-1}\CC[[t]]dt$, where $K$ is the space of Laurent expansions at
$P_0$ of all forms on $X$, regular outside $P_0$, with vanishing
integrals along $a$-cycles.

Therefore 
$$ \delta_{[\xi]} \omega_a(z) = d_z(\xi\omega_a)(z) + d_z\res_{w =
  P_0} \left( d_w \ln\Theta(w - z + (g-1)P_0 - \Delta) \xi(w)
  \omega_a(w) \right) . 
$$ By the properties of sheaf cohomology, we may assume that $\xi$ is
defined inside a loop $\gamma$ deforming the standard loop encircling
$P_0$ and containing both $P_0$ and the $Q_i^{(a)}$.  We have then
\begin{align*}
  \pa_{[\xi]} (\int_{(g-1)P_0}^{\Delta} \omega_a) & = {1\over 2}
  \sum_i \int_{P_0}^{Q_i^{(a)}} d_z(\xi(z) \omega_a(z) ) - \res_{w =
    P_0} G(z,w) d_w(\xi(w) \omega_a(w)) \\ & = {1\over 2} \sum_i
  \int_{P_0}^{Q_i^{(a)}} \int_{\gamma} d_z d_w \ln\Theta(z-w+\delta -
  \Delta) \xi(w) \omega_a(w)
\end{align*}
by the properties of the Green function. Integrating in $z$, we find
$$ {1\over 2} \int_{\gamma} [-2(g-1)d_w \ln\Theta(w - P_0 + \delta -
\Delta) + \sum_i d_w \ln\Theta(w - Q_i^{(a)} + \delta - \Delta) ]
\omega_a(w) \xi(w).
$$ By Lemma \ref{ajax}, this is ${1\over 2} \langle \bar D
\omega_a, \xi \rangle$.  \hfill \qed \medskip

\begin{lemma}\label{summer}
  Recall that $ \omega_\la(z) = \limm_{z'\to z} 2 d_{z'}(G_\la(z,z') -
  G(z,z'))$. We have
\begin{align}
  & \langle \omega_{2\la},\xi \rangle = \sum_a 2
  \pa_{a}\ln\Theta(-2\la + (g-1)P_0 -
  \Delta)(\pa_{[\xi]}\int_{(g-1)P_0}^\Delta \omega_a) \\ &
  \nonumber + \sum_{a,b}{{\delta\tau_{ab}}\over{-2i\pi}}
  {{\pa^2_{ab}\Theta(-2\la + (g-1)P_0 -
      \Delta)}\over{\Theta(-2\la + (g-1)P_0 - \Delta)}} \\ & \nonumber
  + \langle S(P_0,m|z) , \xi \rangle,
\end{align}
with $S(m,P_0|z)$ given by (\ref{def:S}). 
\end{lemma}

{\em Proof.} Recall that $\wt G(z,w) = d_z \ln \Theta(z - w + (g-1)P_0
- \Delta)$, and that $Q_1, \ldots, Q_{g-1}$ are the zeroes (other than
$P_0$) of $\sum_{a} \omega_a(z) \pa_a \Theta((g-1)P_0 - \Delta)$.  We
have $\Theta(z - w + (g-1)P_0 - \Delta) = \Theta(w - z + \sum_{i =
  1}^{g-1} Q_i - \Delta)$.  There exist functions $f$ and $g$ of one
variable such that
\begin{equation} \label{id:separation}
  {{\Theta(w - z + \sum_{i = 1}^{g-1} Q_i - \Delta)}\over{\Theta(w - z
      + (g-1)P_0 - \Delta)}} = f(z) g(w) 
\end{equation} 
(see \cite{Fay}), therefore $\wt G(z,w) - G(z,w)$ is constant in $w$.
It follows that $d_w \wt G(z,w) = d_w G(z,w)$, so that
$$ \omega_{\la}(z) = 2 d_{w} [G_{\la}(z,w) - \wt G(z,w)]_{| w = z}.
$$
Therefore, 
\begin{align*} 
  & \omega_{\la}(z) = [d_w \left\{ { { {{\Theta(z - w + (g-1)P_0 - \la
              - \Delta)}\over{ \Theta((g-1)P_0 - \la- \Delta)}} -
          1}\over{\Theta(z - w + (g-1)P_0 - \Delta)}} \sum_a
      \omega_a(z) \pa_a \Theta((g-1)P_0 - \Delta) \right\} \\ & + d_w
      ({{\sum_a \omega_a(z) [\pa_a \Theta((g-1)P_0 - \Delta) - \pa_a
          \Theta(z - w + (g-1)P_0 - \Delta) ]}\over{\Theta(z - w +
          (g-1) P_0 - \Delta)}})]_{w = z}.
\end{align*}
Expanding both terms around $w = z$, we find 
\begin{align*}
  & \omega_{2\la}(z) = \\ & - {{ - \sum_a \omega'_a(z) \pa_{a}
      \Theta((g-1)P_0-2\la-\Delta) + \sum_{a,b} \omega_a (z)\omega_b(z)
      \pa^2_{ab} \Theta((g-1)P_0 - 2\la -
      \Delta)}\over{\Theta((g-1)P_0 - 2\la - \Delta)}} \\ & + {{-
      \sum_a \omega'_a(z) \pa_{a} \Theta((g-1)P_0-\Delta) +
      \sum_{a,b} \omega_a (z)\omega_b(z) \pa^2_{ab}
      \Theta((g-1)P_0 - \Delta)}\over{ \sum_a \omega_a(z)
      \pa_{a}\Theta((g-1)P_0 - \Delta)}} \cdot \\ & \cdot
  {{\sum_a \omega_a(z) \pa_{a}\Theta((g-1)P_0 -2\la -
      \Delta)}\over{\Theta((g-1)P_0 - 2\la - \Delta)}} +
  S(m,P_0|z) .  
\end{align*}
It follows from (\ref{id:barD}) that 
\begin{align} \label{sachs}
  \omega_{2\la}(z) = & S(m,P_0|z) + \sum_a \bar D \omega_a(z) \pa_{a}
  \ln\Theta((g-1)P_0 - 2\la - \Delta ) \\ & \nonumber - \sum_{a,b}
  \omega_a\omega_b(z) {{\pa^2_{ab}\Theta((g-1)P_0 - 2\la - \Delta
      )}\over{\Theta((g-1)P_0 - 2\la - \Delta )}} ;
\end{align} the Lemma now follows from Lemma \ref{var:int}.
\hfill \qed \medskip

Let us prove Thm.\ \ref{thm:no:point}. (\ref{sachs}) and Cor.\ 
\ref{omega:la:local} imply the statements on $\bar S(m,P_0|z)$.

(\ref{delta:tau}) is equivalent to  
\begin{align*}
  \pa_{[\xi]} f(m_\infty | \la) & = {1\over{8i\pi \kappa}}
  \sum_{a,b} \delta\tau_{ab} \pa^2_{\la_a\la_b} f (m_\infty | \la) +
  {1\over {2\kappa}} \sum_a \langle D^{(2\la)} \omega_a (z), \xi
  \rangle \pa_{\la_a} f (m_\infty | \la) \\ & + {k \over{2\kappa}}
  \langle \bar\omega_{2\la}(z), \xi \rangle f(m_\infty | \la) .
\end{align*}
It follows from
$$ \pa_{[\xi]} \Theta(\eps|m) = {1\over{4i\pi}} \sum_{a,b} \delta
\tau_{ab} \pa^2_{ab} \Theta(\eps|m)
$$ 
that if we set 
$$ g(m_\infty | \la) = f(m_\infty | \la) \Theta(-2\la + (g-1)P_0 -
\Delta|m) ,
$$ the latter equality is equivalent to
\begin{align} \label{interm}
  & \pa_{[\xi]} g(m_\infty | \la) = {1 \over{8i\pi\kappa}}
  \sum_{a,b} \delta \tau_{ab} \pa^2_{\la_a\la_b} g(m_\infty | \la) +
  {1\over {2 \kappa}} \sum_a \langle \bar D \omega_a, \xi\rangle
  \pa_{\la_a} g(m_\infty | \la) \\ & \nonumber + \left( \sum_{a,b} -
    {1 \over{2i\pi}}{{\delta\tau_{ab}}\over\kappa}
    {{\pa^2_{ab}\Theta(-2\la + (g-1)P_0 - \Delta|m)}\over{\Theta(-2\la
        + (g-1)P_0 - \Delta|m)}} + \pa_{[\xi]} \ln\Theta(-2\la +
    (g-1)P_0 - \Delta|m) \right) \\ & \nonumber g(m_\infty | \la) \\ &
  \nonumber + \left( \sum_a {1\over{2\kappa}} \langle \wt D \omega_a ,
    \xi\rangle 2 \pa_{a} \ln\Theta(-2\la + (g-1)P_0 - \Delta|m) +
    {k\over{2\kappa}} \langle \bar\omega_{2\la}, \xi\rangle \right)
  g(m_\infty | \la).
\end{align}
It follows from Lemma \ref{summer}, that the coefficient of $g(m_\infty | \la)$ in
the right side of this equality is equal to 
$$ {{k}\over{2\kappa}}\langle \bar S(m,P_0|z), \xi(z) \rangle_z , 
$$ 
so that (\ref{interm}) is rewritten as
\begin{align} \label{interm'}
  & \pa_{[\xi]} g(m_\infty | \la) = {1 \over{8i\pi\kappa}}
  \sum_{a,b} \delta \tau_{ab} \pa^2_{\la_a\la_b} g(m_\infty | \la) +
  {1\over {2 \kappa}} \sum_a \langle \bar D \omega_a, \xi\rangle
  \pa_{\la_a} g(m_\infty | \la) \\ & \nonumber 
 + {{k}\over{2\kappa}}\langle \bar S(m,P_0|z), \xi(z) \rangle_z  g(m_\infty | \la) . 
\end{align}

Set then 
$$ h(m_\infty | \la') = g(m_\infty | {1\over\kappa}(\la' + (g-1)P_0 -
\Delta)).
$$ Set $\la = {1\over\kappa}(\la' + (g-1)P_0 - \Delta)$. We have 
\begin{align*}
  & \pa_{[\xi]} h(m_\infty | \la') = \pa_{[\xi]} g(m_\infty |
  \la ) \\ & + \sum_a \langle \pa_{[\xi]} [{1\over\kappa}( (g-1)P_0
  - \Delta)], \omega_a \rangle \pa_{\la_a} g(m_\infty | \la ) \\ & =
  {1\over{8i\pi\kappa}} \sum_{a,b} \delta\tau_{ab} \pa^2_{\la_a\la_b}
  g(m_\infty | \la) + {1\over{2\kappa}} \sum_a \langle \bar D\omega_a,
  \xi \rangle \pa_{\la_a} g(m_\infty | \la) \\ & + {1\over\kappa}
  \sum_a \langle \pa_{[\xi]} ( (g-1)P_0 - \Delta), \omega_a\rangle
  \pa_{\la_a} g(m_\infty | \la) +  {{k}\over{2\kappa}}
  \langle \xi(z),  \bar S(m,P_0|z)  \rangle_z 
  g(m_\infty | \la),
\end{align*}
Lemma \ref{var:int} implies that the second and third terms of this
expression are opposite.  (\ref{interm'}) is therefore equivalent to
\begin{equation} \label{eq:for:h}
 \pa_{[\xi]} h(m_\infty | \la') = {{\kappa}\over{8i\pi}}
\sum_{a,b} \delta \tau_{ab} \pa^2_{\la'_a \la'_b} h(m_\infty | \la') +
{{k}\over{2\kappa}}\langle \xi(z) , \bar S(m,P_0|z) \rangle_z h(m_\infty | \la').
\end{equation} 
Consider the connection over the trivial bundle over $\cM_{g,1^2}$
with fiber $\CC[[\la']]$, defined by 
$$ \nabla^{aux}_{[\xi]} = \pa_{\xi} -
\left({{\kappa}\over{8i\pi}}\sum_{a,b} \delta \tau_{ab} \pa^2_{\la'_a
    \la'_b} + \langle \bar \xi(z),  S(m,P_0,z) \rangle_z \right).
$$ It follows from Thm.\ \ref{thm:flatness} that this connection is
flat. Therefore, the connection over the trivial bundle over $\cM_{g,1^2}$
with fiber $\CC$, defined by 
$$ \nabla^{aux \prime}_{[\xi]} = \pa_{\xi} - \langle \bar
S(m,P_0,v_0), \xi\rangle ,
$$ is flat. This implies the existence of a nonzero local solution to
(\ref{eq:for:alpha}).

Moreover, since (\ref{eq:for:h}) is separated, its local solutions are
exactly the functions $h(m_\infty | \la')$ of the form $\al(m,P_0,v_0)
\varphi(\la')$, where $\al(m,P_0,v_0)$ is a solution of
(\ref{eq:for:alpha}) and $\varphi(\la')$ is a solutions to 
$$ \pa_{[\xi]} \varphi(\la') = {{\kappa}\over{8i\pi}} \sum_{a,b}
\delta\tau_{ab} \pa^2_{\la'_a\la'_b} \varphi(\la').
$$

This ends the proof of Thm.\ \ref{thm:no:point}. \hfill \qed \medskip

\begin{remark} \label{conjectures} {\it Conjecture on the form of $\al(m,P_0,v_0)$.}
  The ${g\over 2}$-form $\sigma(z)$ is determined up to multiplication
  by a $(g-1)$th root of $1$, by the identity
  $$ \Theta(g x - y - \Delta) = W(x) E(x,y)^{g}
  {{\sigma(y)}\over{\sigma(x)^g}} ,
  $$ where $W(x)$ is the Wronskian of the forms $(\omega_a)_{a = 1,
    \ldots, g}$, $E(x,y)$ is the prime form relative to the divisor
  $\delta$ (its is normalized by its behavior at the diagonal $E(x,y)
  \sim (x - y) (dx)^{-1/2}(dy)^{-1/2}$). We will emphazise the moduli
  dependence of $\sigma(z)$ by writing it $\sigma(m|z)$
  (the functional properties of $\sigma$ mean that it realizes the
  isomorphism of $\Omega_X^{g/2}$ with $\underline{\Theta}^{g-1}$, where 
  $\underline{\Theta}$ is the theta line bundle over $X$).

  We expect that the ``logarithmic primitives'' of the forms
  $S(m,P_0,z)$ and $\sum_a \omega_a dr_a^{(P_0)}$ will be respectively
  of the form $f(m) \sigma(m|z)^{-2(g-1)}$ and $W(m|z)^{-1}$, where $f(m)$
  has the form $\exp(\sum_{a,b} \al_{ab}\tau_{ab})$. More precisely, we should have
  $$ \varsigma(m,P_0,v_0) := ( f(m) \sigma(m|P_0)^{-2(g-1)}, v_0 ) \quad
  \on{and} \quad \varpi(m,P_0,v_0) := ( W(m|P_0)^{-1}, v_0 )
  $$ are functions on $\cM_{g,1^2}$, such that
  $$ \pa_{[\xi]} \ln \varsigma(m,P_0,v_0) = \langle S(m,P_0| \cdot ) , \xi
  \rangle \quad \on{and} \quad \pa_{[\xi]} \ln \varpi(m,P_0,v_0) =
  \langle \sum_a \omega_a dr_a^{(P_0)} , \xi \rangle .
  $$ 
  These statements imply that $\al(m,P_0,v_0)$ has the form
  $$ \al(m,P_0,v_0) = ( f(m) \sigma(m|P_0)^{-2(g-1)}
  W(m|P_0)^3, v_0 ).
  $$
  We checked the restrictions of these statements to the submanifold of
  Mumford curves.  
\end{remark}

\newpage

\section{Flat sections of $\cF_{\ell,(\La_i)}$ in the case $N = 0$,  $n$ arbitrary}
\label{sect:N=0:n:arb}

We assume now that $\ell$ and the $\La_i$ are such that $N = 0$ and $n$ is $\geq 0$.
We again write $\cF$ for the sheaf  $\cF_{\ell,(\La_i)}$.  
We have 
\begin{align*}
  & (T^{\cF}_{\bar\omega}(z) f)(m_\infty,P_i,v_i|\la) = \\ & 
  {1\over {2\kappa}}\left[ {1\over 2}
    (\sum_a \omega_a(z) \pa_{\la_a} - \sum_i \La_i G(z,P_i)_{z\ll P_i} )^2 +
    \right. \\ & \left.  \sum_{a} (D^{(2\la)} \omega_a)(z) \pa_{\la_a} - \sum_i \La_i
    D_z^{(2\la)} G(z,P_i)_{z\ll P_i}  + k \bar\omega_{2\la}(z)
  \right] f(m_\infty,P_i,v_i|\la).
\end{align*}

A flat section of the connection $\nabla^{\cF}$ is a function 
$f(m_\infty,P_i,v_i|\la)$ satisfying 
\begin{equation} \label{flatness:P_i}
 \pa_{[\xi]} f(m_\infty,P_i,v_i|\la) = (T^{\cF}_{\bar\omega}(z) f)
  (m_\infty,P_i,v_i|\la) . 
\end{equation}

Let us fix an odd nonsingular theta-characteristic $\ddelta$. $\ddelta + \Delta$ is
an effective divisor of $J^{g-1}(X)$. We have $\ddelta = 2i\pi \varepsilon 
+ \tau \varepsilon'$, with $\varepsilon,\varepsilon'$ in $({1\over 2}\ZZ)^g$, 
and  
\begin{equation} \label{def:theta:delta}
\Theta[\ddelta](\zz) = \Theta(\zz + \ddelta) \exp({1\over 2} \varepsilon'  
\tau \varepsilon^{\prime t} + (\zz + 2i\pi \varepsilon) \varepsilon^{\prime t}).  
\end{equation}
Let us set $\wt\Theta = \Theta[\ddelta]$.  $\wt\Theta$ is an odd function of 
$\zz$.  The divisor of the one-form $\sum_a \pa_a\wt\Theta(\ddelta) \omega_a(z)$ 
is $2\ddelta$ (this one-form is called a Prym form).  Let $\psi(z)$ be a
half-form such that 
\begin{equation} \label{def:psi}
\psi(z)^2 = \sum_a \pa_a\wt\Theta(\ddelta) \omega_a(z).  
\end{equation}
The prime-form $E(z,w)$ is then defined by 
$$
E(z,w) = {{\wt\Theta(z-w)}\over{\psi(z)\psi(w)}} . 
$$

The aim of Prop.\ \ref{koche} is to integrate the part of 
$T^{\cF}_{\bar\omega}(z)f$ which is linear in $\La_i$. 

\begin{prop} \label{koche}
1) For $Q$ in $X$, let us set 
\begin{align*}
& S'(m,P_0,Q|z) = - \sum_{a,b} {1\over 2} (\omega_a\omega_b)(z) 
{{\pa^2_{ab}\Theta}\over{\Theta}} (gP_0 - Q - \Delta)
\\ & + \sum_a {{\pa_a\Theta}\over{\Theta}}(gP_0 - Q - \Delta) [
{1\over 2} \bar D\omega_a(z) + \omega_a(z) G(z,Q)]
- {1\over 2} G(z,Q)^2 - {1\over 2} \bar D_z G(z,Q) .  
\end{align*}
The function $S'(m,P_0,Q|z)$ is independent of $Q$. We denote it 
$S'(m,P_0|z)$. It is a regular 
quadratic differential on $X - \{P_0\}$. It has the expansion at $P_0$
\begin{equation} \label{local:S'}
S'(m,P_0|z) = -g {{(dz_{P_0})^2}\over{z_{P_0}^2}} + O(1/z_{P_0}) , 
\end{equation}
where $z_{P_0}$ is a local coordinate at $P_0$. 

2) There exists a function $\beta(m,P_0,v_0)$, locally defined on $\cM_{g,1^2}$, 
homogeneous of degree $-g$ in $v_0$, satisfying the equation 
\begin{equation} \label{id:beta}
\pa_{[\xi]} \ln\beta(m,P_0,v_0)  = \langle \xi(z), S'(m,P_0, z)\rangle_z.  
\end{equation}
We will denote by $\beta(m,P_0)$ the $-g$-form such that 
$\beta(m,P_0,v_0) = (\beta(m,P_0) , v_0 )$. 

We have then 
\begin{equation} \label{eagle}
\pa_{[\xi]} \ln \{ {1\over{\Theta(gP_0 - P_i - \Delta) }} 
(E(P_0,P_i), v_0\otimes v_i) (\beta(m,P_0), v_0 )   \}
= {1\over 2}\langle \xi , \bar D_z G(z,P_i) \rangle . 
\end{equation}
\end{prop}

\begin{thm} \label{thm:pts}
The solutions of (\ref{flatness:P_i}) are all of the form 
\begin{align} \label{sol:P_i}
& f_\varphi(m_\infty, P_i,v_i | \la ) =
\prod_{i<j} \left( {{E(P_i,P_j)}\over{ E(P_0,P_i)E(P_0,P_j)}} , v_0
\right)^{{{\La_i\La_j}\over{2\kappa}}}
\prod_i (E(P_0,P_i) , v_0 \otimes v_i)^{ - {{\La_i(\La_i + 2)}\over{2\kappa}}} 
\cdot \\ & \nonumber \cdot 
\prod_i \left( {{ (\beta(m,P_0),v_0) }\over
{\Theta(gP_0 - P_i - \Delta)}}\right)^{ - {\La_i \over\kappa}}  
\al(m,P_0,v_0)^{k \over{2\kappa}} 
\cdot \\ & \nonumber \cdot 
\varphi(m | \kappa\la + \Delta - (g-1)P_0 - \sum_i \La_i (P_i - P_0) ) 
\Theta( -2\la + (g-1)P_0 - \Delta)^{-1} ,  
\end{align}
where $\varphi(m|\la)$ is a solution of 
$$ 
 \pa_{[\xi]} \varphi(m|\la) = 
  {\kappa \over{8 i \pi}} \sum_{a,b} \delta\tau_{ab} 
 \pa^2_{\la_a\la_b}  \varphi(m|\la) . 
$$
(At the end of Thm.\ \ref{thm:no:point}, we indicate how such functions $\varphi$
can be naturally obtained.)
\end{thm}

Before we turn to the proof of Prop.\ \ref{koche} and Thm.\
\ref{thm:pts},  we will compute logarithmic derivatives of functions
defined on $\cM_{g,1^2}$ (sect.\ \ref{sect:quadratic}).  These
computations will serve to integrate the  part of
$T^{\cF}_{\bar\omega}(z)f$ which is quadratic in the $\La_i$:  Prop.\
\ref{prop:var:ratio:E} will treat the terms in  $\La_i\La_j$ and Prop.\
\ref{var:E} will treat the terms in  $\La_i^2$.

\subsection{Integration of quadratic terms in the $\La_i$} 
\label{sect:quadratic}

\begin{prop} \label{prop:var:ratio:E}
Assume that $i\neq j$. We have 
$$
 \pa_{[\xi]} \ln ({{E(P_i,P_j)}\over{E(P_0,P_i) E(P_0,P_j)}}, 
 v_0 )  = \langle \xi(z), G(z,P_i)_{z\ll P_i} G(z,P_j)_{z\ll P_j}  
 \rangle_z . 
$$
\end{prop}

{\em Proof.}
We first show: 

\begin{lemma} \label{var:int:pts}
$$
\pa_{[\xi]}(\int_{P_0}^{P_i} \omega_a)  = \res_{z = P_0}
[(\xi\omega_a)(z) G(z,P_i)_{z \ll P_i}] . 
$$
\end{lemma}

{\em Proof of Lemma.} $\int_{P_0}^{P_i} \omega_a + \eps
\pa_{[\xi]}(\int_{P_0}^{P_i} \omega_a)$ is equal to 
\begin{equation} \label{verne}
\limm_{t \to P_0} \int_t^{P'_i} \omega'_a,
\end{equation} 
where $P'_i$ is the point corresponding to $P_i$ in the curve
$\Spec(R_\eps)  = \Spec(1 + \eps\xi)(R)$. Recall from the proof of Lemma
\ref{var:moduli} that we have $\omega'_a = (1 + \eps\cL_\xi)(\omega_a) + \eps
\delta_{[\xi]} \omega_a$, where 
$\delta_{[\xi]} \omega_a = - [d(\xi\omega_a)]_{out}$. Changing (\ref{verne}) 
into an integral over $\Spec(R)$, we find that (\ref{verne}) is equal  to 
\begin{equation} \label{dvd}
\limm_{t\to P_0} \int_{(1 - \eps\xi)(t)}^{P_i} ( \omega_a 
+ \eps \delta_{[\xi]}\omega_a).  
\end{equation}
(\ref{dvd}) is equal to 
$$
\int_{P_0}^{P_i} \omega_a + \eps[\limm_{t\to P_0} \int_{(1 - \eps \xi)(t)}^{P_i}
 - (d(\xi\omega_a))_{out} ] , 
$$
therefore the derivative $\pa_{[\xi]}(\int_{P_0}^{P_i} \omega_a)$ is equal to 
\begin{equation} \label{curie}
\pa_{[\xi]}(\int_{P_0}^{P_i} \omega_a) = 
\limm_{t\to P_0} [ - \int_{t}^{P_i} [d(\xi\omega_a)]_{out} + (\xi\omega_a)(t)] . 
\end{equation}
Let us set, for $f$ in $\cK$, 
$$
\phi_{P_i}(f) = \limm_{t\to P_0} [ - \int_{t}^{P_i} (df)_{out} + f(t)] . 
$$
We find that for $f$ in $R_a$, $\phi_{P_i}(f) = - f(P_i)$ and for $f$ in $z\CC[[z]]$, 
$\phi_{P_i}(f) = 0$. Therefore, we have 
\begin{equation} \label{expr:phi:i}
\phi_{P_i}(f) = - \langle G(z,P_i)_{z\ll P_i}, f(z) \rangle_z .   
\end{equation}
The Lemma now follows from 
(\ref{expr:phi:i}) 
and (\ref{curie}). 
\hfill \qed\medskip

\begin{lemma} \label{var:ratio:theta}
We have 
$$
\pa_{[\xi]}(\ln{{\wt\Theta(P_i - P_j)}\over{\wt\Theta(P_i - P_0)\wt\Theta(P_j - P_0)}}) = 
\langle \xi, \gamma(\cdot, P_i,P_j)  \rangle ,   
$$
where 
\begin{align*}
\gamma(z,P_i,P_j) = & \sum_{a,b = 1}^g  {1\over 2}
(\omega_a \omega_b)(z) \left( {{\pa^2_{ab}\wt\Theta}\over{\wt\Theta}}(P_i - P_j) 
- {{\pa^2_{ab}\wt\Theta}\over{\wt\Theta}}(P_i - P_0)
- {{\pa^2_{ab}\wt\Theta}\over{\wt\Theta}}(P_j - P_0)
\right)
\\ & + \sum_a \omega_a(z) G(z,P_i)_{z\ll P_i}  
\left( {{\pa_a \wt\Theta}\over{\wt\Theta}}(P_i - P_j)  
- {{\pa_a \wt\Theta}\over{\wt\Theta}}(P_i - P_0) 
\right)   
\\ & + \sum_a \omega_a(z) G(z,P_j)_{z\ll P_j}  
\left( {{\pa_a \wt\Theta}\over{\wt\Theta}}(P_j - P_i)  
- {{\pa_a \wt\Theta}\over{\wt\Theta}}(P_j - P_0) 
\right)  . 
\end{align*}
\end{lemma}

{\em Proof.} We have 
\begin{align} \label{7eme}
\pa_{[\xi]} \wt\Theta(P_i - P_0) & \nonumber = \sum_{a,b} 
\pa_{[\xi]} \tau_{ab} {1\over{4i\pi}}
\pa^2_{ab}\wt\Theta(P_i - P_0) + \sum_a \pa_{[\xi]}(\int_{P_0}^{P_i} \omega_a)
\pa_a \wt\Theta(P_i - P_0) 
\\ & = \langle \xi, {1\over 2}\sum_{a,b} (\omega_a \omega_b)(z)
\pa^2_{ab}\wt\Theta(P_i - P_0) + \sum_a \omega_a(z) G(z,P_i)_{z\ll P_i}  
\pa_a\wt\Theta(P_i - P_0) \rangle     
\end{align}
by virtue of Lemmas \ref{var:moduli} and \ref{var:int:pts}. 
In the same way, 
\begin{align*} 
\pa_{[\xi]} \wt\Theta(P_i - P_j) &  = \langle \xi, {1\over 2}\sum_{1\leq
a,b\leq g} (\omega_a \omega_b)(z) \pa^2_{ab}\wt\Theta(P_i - P_0) \\ &  
+ \sum_a \omega_a(z) [G(z,P_i)_{z\ll P_i} - G(z,P_j)_{z\ll P_j}] 
\pa_a\wt\Theta(P_i - P_0) \rangle.     
\end{align*}
The Lemma follows. 
\hfill \qed \medskip 

\begin{lemma} \label{var:omega:a}
We have
$$
\pa_{[\xi]} ( \omega_a(P_0), v_0)  
= \langle \xi(z), - \sum_{b,c = 1}^g  \pa^2_{bc} \ln\wt\Theta(z-P_0)
(\omega_a  \omega_b)(z) ( \omega_c(P_0), v_0) \rangle_z  .  
$$
\end{lemma}

{\em Proof.} $(\omega_a(P_0),v_0 )  
+ \eps \pa_{[\xi]} (\omega_a(P_0),v_0 )$ is equal to 
$(\omega'_a(P_0), v_0 )$, in the notation of the proof of 
Lemma \ref{var:moduli}. Therefore, if we set $\omega_{in} = 
\omega - \omega_{out}$, we have 
$$
\pa_{[\xi]} (\omega_a(P_0), v_0 ) = 
(\  (d(\xi\omega_a))_{in} (P_0) , v_0 \ ). 
$$
We have $(df)_{in} = 0$ if $f$ belongs to $R_a$ and $(df)_{in} = df$ if
$f$ belongs to $z\CC[[z]]$, therefore for $f$ in $\cK$, we have 
$(df)_{in}(u) = \langle f(t), d_u G(t,u)_{u \ll t} \rangle_t$ (the pairing is in
variable $t$). 
Therefore,
$$ 
\pa_{[\xi]} (\omega_a(P_0), v_0) = 
( \langle \xi(t), \omega_a(t) d_u G(t,u)_{u\ll t}\rangle_{t | u = P_0}
, v_0) . 
$$
The identity (\ref{id:separation}) implies that $d_u G(t,u) = d_u d_t
\ln\wt\Theta(t-u)$.  The Lemma now follows from the equalities  $d_u
G(t,u) = d_u d_t \ln\wt\Theta(t-u)  = - \sum_{b,c} 
\pa^2_{bc}\ln\wt\Theta(t-u) \omega_b(t)\omega_c(u)$.  \hfill
\qed\medskip 

\begin{lemma} \label{var:psi}
We have 
$$
\pa_{[\xi]} ( \psi(P_0)^2, v_0 ) = \langle \xi, \rho \rangle, 
$$
where
\begin{align*}
& \rho(z) =
\\ &  \sum_{a,b,c} \left( {1\over 2} \pa^3_{abc}\wt\Theta(0)
(\omega_b \omega_c)(z) ( \omega_a(P_0),  v_0 ) - 
\pa_a \wt\Theta(0) \pa^2_{bc} \ln\wt\Theta(z - P_0) 
(\omega_a\omega_b)(z) ( \omega_c(P_0),  v_0 ) \right) .  
\end{align*}
\end{lemma}

{\em Proof.} This follows from (\ref{def:psi}), Lemma \ref{var:omega:a}
and 
$$
\pa_{[\xi]} (\pa_a \wt\Theta(0)) = \langle \xi(z) , \sum_{b,c} {1\over 2}
\pa^{3}_{abc}\wt\Theta(0) (\omega_b \omega_c) (z) \rangle.  
$$
\hfill \qed\medskip 

\begin{lemma} \label{var:ratio:E}
We have 
$$
\pa_{[\xi]} \ln ({{E(P_i,P_j)}\over{E(P_i,P_0)E(P_j,P_0)}} , v_0)
= \langle \xi , \gamma'(\cdot, P_i,P_j) \rangle ,    
$$
where 
\begin{align*}
& \gamma'(z,P_i,P_j) = \gamma(z,P_i,P_j) + {1\over {\sum_a \omega_a(P_0) 
\pa_a\wt\Theta(0) }} \cdot \\ & \cdot \sum_{a,b,c}
\left( {1\over 2} \pa^3_{abc}\wt\Theta(0) (\omega_b \omega_c)(z) 
\omega_a(P_0) - \pa_a\wt\Theta(0) \pa^2_{bc}\ln\wt\Theta(z - P_0)
(\omega_a \omega_b)(z) \omega_c(P_0)   \right) . 
\end{align*}
\end{lemma}

{\em Proof.} This follows from Lemmas \ref{var:ratio:theta} and \ref{var:psi}.
\hfill \qed\medskip 

\begin{lemma} \label{id:var:rat:E}
$$
\gamma'(z,P_i,P_j)  = G(z,P_i)_{z\ll P_i} G(z,P_j)_{z\ll P_j}.   
$$
\end{lemma}

{\em Proof.} Let us set $\gamma''(z,P_i,P_j) =  \gamma(z,P_i,P_j)  
- G(z,P_i)_{z\ll P_i} G(z,P_j)_{z\ll P_j}$. $\gamma''(z,P_i,P_j)$ is equal to 
\begin{align} \label{expr:gamma'}
\sum_{a,b} (\omega_a\omega_b)(z)  & \left\{ \right. 
{1\over 2} (
{{\pa^2_{ab}\wt\Theta}\over{\wt\Theta}}(P_i - P_j)  
- {{\pa^2_{ab}\wt\Theta}\over{\wt\Theta}}(P_i - P_0) 
- {{\pa^2_{ab}\wt\Theta}\over{\wt\Theta}}(P_j - P_0) ) 
\nonumber\\ & 
+ [ {{\pa_a\wt\Theta}\over{\wt\Theta}}(P_i - P_j)  
- {{\pa_a\wt\Theta}\over{\wt\Theta}}(P_i - P_0) ]
[ {{\pa_b\wt\Theta}\over{\wt\Theta}}(z - P_i)_{z\ll P_i}  
- {{\pa_b\wt\Theta}\over{\wt\Theta}}(z - P_0) ]
\nonumber \\ &  
+ [ {{\pa_a\wt\Theta}\over{\wt\Theta}}(P_j - P_i)  
- {{\pa_a\wt\Theta}\over{\wt\Theta}}(P_j - P_0) ]
[ {{\pa_b\wt\Theta}\over{\wt\Theta}}(z - P_j)_{z\ll P_i}  
- {{\pa_b\wt\Theta}\over{\wt\Theta}}(z - P_0) ]
\nonumber \\ & 
- [ {{\pa_a\wt\Theta}\over{\wt\Theta}}(z - P_i)_{z\ll P_i}  
- {{\pa_a\wt\Theta}\over{\wt\Theta}}(z - P_0) ]
[ {{\pa_b\wt\Theta}\over{\wt\Theta}}(z - P_j)_{z\ll P_i}  
- {{\pa_b\wt\Theta}\over{\wt\Theta}}(z - P_0) ] \left. \right\} . 
\end{align}
(\ref{expr:gamma'}) has an obvious analytic prolongation to $\wt X$, 
which we still denote $\gamma''(z,P_i,P_j)$. 
Let us study its behaviour as a function of 
$P_i$. It has no poles at $P_j,P_0$ or $z$. Moreover, the rules 
$$
{{\pa_a\wt\Theta}\over{\wt\Theta}}(\gamma_{B_c}(z)  - w)  
= {{\pa_a\wt\Theta}\over{\wt\Theta}}(z  - w) - \delta_{ac},  
$$
$$
{{\pa^2_{ab}\wt\Theta}\over{\wt\Theta}}(\gamma_{B_c}(z) - w)  
= {{\pa^2_{ab}\wt\Theta}\over{\wt\Theta}}(z  - w) 
- \delta_{ac} {{\pa_b\wt\Theta}\over{\wt\Theta}}(z  - w)  
- \delta_{bc} {{\pa_a\wt\Theta}\over{\wt\Theta}}(z  - w) 
+ \delta_{ac}\delta_{bc}
$$
imply that $P_i \mapsto \gamma''(z,P_i,P_j)$ is univalued on $X$. It
follows that  $P_i \mapsto \gamma''(z,P_i,P_j)$ is a constant map.  
$\gamma''(z,P_i,P_j)$ is obviously symmetric in $P_i$ and $P_j$. 
Therefore, it is constant as a function of $(P_i,P_j)$. 

It follows that $\gamma''(z,P_i,P_j)$ is equal to $\limm_{P_i \to P_0}
\gamma''(z,P_i,P_j)$. Let us fix a coordinate at $P_0$ and let $\eps$ be the coordinate of 
$P_i$. We have the expansions 
$$
\wt\Theta(P_i - P_0) = \eps \sum_a \omega_a(P_0) \pa_a\wt\Theta(0) 
+ o(\eps), 
$$
and since $\pa^2_{ab}\wt\Theta(0) = 0$, 
$$
\pa^2_{ab}\wt\Theta(P_i - P_0) = \eps \sum_c \omega_c(P_0) \pa^2_{abc}\wt\Theta(0) 
+ o(\eps), 
$$
therefore 
$$
\limm_{P_i \to P_0} - {1\over 2}{{ \pa^2_{ab}\wt\Theta }\over{\wt\Theta}}(P_i - P_0) = 
- {1\over 2} {{ \sum_c \omega_c(P_0) \pa^2_{abc}\wt\Theta(0)  } \over
{ \sum_a \omega_a(P_0) \pa_a\wt\Theta(0)  }}. 
$$
On the other hand, we have 
$$
{{\pa_a\wt\Theta}\over{\wt\Theta}}(P_i - P_0) = {1\over \eps} {{\pa_a\wt\Theta(0) }\over
{\sum_a \omega_a(P_0) \pa_a\wt\Theta(0) }} + O(1),  
$$
and
$$
{{\pa_b\wt\Theta}\over{\wt\Theta}}(z - P_i)  
- {{\pa_b\wt\Theta}\over{\wt\Theta}}(z - P_0)
= - \eps \sum_a \omega_a(P_0) \pa_a({{\pa_b\wt\Theta}\over{\wt\Theta}})(z - P_0).   
$$ 
Therefore, 
\begin{align*}
& \limm_{P_i \to P_0}
[ {{\pa_a\wt\Theta}\over{\wt\Theta}}(P_i - P_j)  
- {{\pa_a\wt\Theta}\over{\wt\Theta}}(P_i - P_0) ]
[ {{\pa_b\wt\Theta}\over{\wt\Theta}}(z - P_i)  
- {{\pa_b\wt\Theta}\over{\wt\Theta}}(z - P_0) ] 
\\ & = 
{\pa_a \wt\Theta(0)} {{
\sum_c \omega_c(P_0) \pa_c ({{\pa_b\wt\Theta}\over{\wt\Theta}}) (z-P_0) }
\over{\sum_a \omega_a(P_0) \pa_a\wt\Theta(0)}} . 
\end{align*}
The limit of the two last lines of  (\ref{expr:gamma'}) is zero when $P_i$ tends 
to $P_0$. 
Therefore  $\gamma''(z,P_i,P_j)$ is equal to 
$$
\sum_{a,b} (\omega_a\omega_b)(z) 
\left[
- {1\over 2} {{ \sum_c \omega_c(P_0) \pa^2_{abc}\wt\Theta(0)  } \over
{ \sum_a \omega_a(P_0) \pa_a\wt\Theta(0)  }} 
+ {\pa_a \wt\Theta(0)} {{
\sum_c \omega_c(P_0) \pa_c ({{\pa_b\wt\Theta}\over{\wt\Theta}}) (z-P_0) }
\over{\sum_a \omega_a(P_0) \pa_a\wt\Theta(0)}} \right] . 
$$
This proves the Lemma.  \hfill \qed\medskip 

Prop.\ \ref{prop:var:ratio:E} now follows from Lemma \ref{var:ratio:E}
and Lemma \ref{id:var:rat:E}.  
\hfill \qed\medskip 

\begin{prop} \label{var:E}
Recall that elements of $\cM_{g,1^\infty,1^2}$ correspond to quadruples 
$(X,t,P_i,v_i)$, where $t$ is a formal of coordinate at a point $P_0$ of $X$,
$P_i$ is a point of $X - \{P_0\}$ and $v_i$ is a tangent vector at this
point.  We denote by $v_0$ the tangent vector at $P_0$ obtained from the
coordinate $t$. 
Then we have 
$$
\pa_{[\xi]} \ln (E(P_0,P_i)^{-2} , v_0 \otimes v_i) = 
\langle \xi(z) , G(z,P_i)_{z\ll P_i}^2 \rangle_z . 
$$
\end{prop}

{\em Proof.}
We first prove 
\begin{lemma} \label{var:psi:P_i}
Let us set 
$$
\rho_i(z) = \sum_{a,b,c} {1\over 2} \pa^3_{abc} \wt\Theta(0) (\omega_b\omega_c)(z)
(\omega_a(P_i) , v_i ) - \pa_a\wt\Theta(0) \pa^2_{bc} \ln\wt\Theta(z - P_i)
(\omega_a\omega_b)(z) (\omega_c(P_i), v_i ).          
$$
Then we have 
$$
\pa_{[\xi]} (\psi(P_i)^2, v_i) = \langle \xi , \rho_i \rangle .  
$$
\end{lemma}

{\em Proof.} We have 
\begin{align} \label{anov}
& \nonumber \pa_{[\xi]} (\psi(P_i)^2, v_i ) =  
\pa_{[\xi]} \{\sum_a \pa_a\wt\Theta(0) ( \omega_a(P_i), v_i ) \} 
\\ & \nonumber  = \sum_a (\pa_{[\xi]} \pa_a\wt\Theta(0))  ( \omega_a(P_i), v_i ) 
+ \sum_a \pa_a\wt\Theta(0) \pa_{[\xi]}  ( \omega_a(P_i), v_i )  
\\ & = 
\sum_a \langle \xi , {1\over 2} \sum_{b,c} (\omega_b\omega_c)(z)
\pa^3_{abc}(0) \rangle ( \omega_a(P_i), v_i )  
+ \sum_a \pa_a\wt\Theta(0) \pa_{[\xi]} ( \omega_a(P_i), v_i ) . 
\end{align}
We have $\pa_{[\xi]} \omega_a = - [d(\xi\omega_a)]_{out}$ so 
$$
\pa_{[\xi]} (\omega_a(P_i),v_i)  = - ( [d(\xi\omega_a)]_{out} , v_i ) . 
$$
As we have seen, for any $f$ in $\cK$, we have 
$$
(df)_{out}(P_i) = \sum_{b,c} \omega_c(P_i) \langle f (z), 
\pa^2_{bc}\ln\wt\Theta(z - P_i)_{z\ll P_i} \omega_b(z) \rangle_z .   
$$
It follows that 
\begin{equation} \label{bogd}
\pa_{[\xi]} (\omega_a(P_i),v_i)  = - \langle \xi ,  \sum_{b,c}
\pa^2_{bc} \ln\wt\Theta(z - P_i)_{z\ll P_i} (\omega_a \omega_b)(z)
(\omega_c(P_i), v_i  )  \rangle .   
\end{equation}
The Lemma now follows from (\ref{bogd}) and (\ref{anov}). 
\hfill \qed\medskip 

We have then 
\begin{align} \label{chev}
& \pa_{[\xi]} \ln (E(P_0,P_i)^{-2} , v_0 \otimes v_i) = 
\pa_{[\xi]} \ln (\psi(P_0)^2 , v_0)  + \pa_{[\xi]} \ln (\psi(P_i)^2 , v_i)
- 2 \pa_{[\xi]} \ln\wt\Theta(P_0 - P_i)
\nonumber \\ & = 
\langle \xi , {{\rho}\over{ (\psi(P_0)^2 , v_0) }} 
\rangle    
+ \langle \xi , {{\rho_i}\over{ (\psi(P_i)^2, v_i ) }}  
\rangle   
\nonumber \\ & + 
\langle \xi, -\sum_{a,b} (\omega_a\omega_b)(z) \pa^2_{ab}\wt\Theta(P_i - P_0)
- 2 \sum_a \omega_a(z) G(z, P_i)_{z\ll P_i} \pa_a\wt\Theta(P_i - P_0) 
\rangle       
\end{align}
by virtue of (\ref{7eme}) and Lemmas \ref{var:psi} 
and \ref{var:psi:P_i}. 

We have then 
\begin{lemma} \label{new:id}
We have 
\begin{align*}
& {{\rho(z) }\over{ (\psi(P_0)^2 , v_0) }}  
+ {{\rho_i(z) }\over{ (\psi(P_i)^2, v_i ) }} 
\\ & -\sum_{a,b} (\omega_a\omega_b)(z) \pa^2_{ab}\wt\Theta(P_i - P_0)
- 2 \sum_a \omega_a(z) G(z, P_i)_{z\ll P_i} \pa_a\wt\Theta(P_i - P_0) 
= G(z,P_i)_{z\ll P_i}^2 . 
\end{align*}
\end{lemma}

{\em Proof.} It follows from Lemma \ref{id:var:rat:E} that $\gamma'(z,P_i,P_j)
= G(z,P_i)_{z\ll P_i}G(z,P_j)_{z\ll P_j}$. Therefore, $G(z,P_i)_{z\ll P_i}^2 
= \limm_{P_j\to P_i}
\gamma'(z,P_i,P_j)$. The latter expression is  equal to 
$$
\limm_{P_j\to P_j} \gamma(z,P_i,P_j)  
+ {{\rho(z) }\over{ (\psi(P_0)^2 , v_0) }}  . 
$$
Now 
\begin{align*}
& \limm_{P_j\to P_j} \gamma(z,P_i,P_j)  \\ & = 
- \sum_{a,b} (\omega_a\omega_b)(z) \pa^2_{ab}\wt\Theta(P_i - P_j) 
- 2 \sum_a \omega_a(z) G(z,P_i)_{z\ll P_i} \pa_a\wt\Theta(P_i - P_0)   
\\ & + \limm_{P_j \to P_i} \left( 
\sum_{a,b} {1\over 2} (\omega_a\omega_b)(z)
{{\pa^2_{ab}\wt\Theta}\over{\wt\Theta}}(P_i - P_j)
\right. \\ & + \left. 
\sum_a \omega_a(z) [G(z,P_i)_{z\ll P_i} - G(z,P_j)_{z\ll P_j}] {{\pa_a\wt\Theta}\over
{\wt\Theta}}(P_i - P_j)   
\right) . 
\end{align*}
It follows from the computations of Lemma \ref{id:var:rat:E} (with $P_0$
 replaced by $P_i$) that the latter limit is equal to  ${{\rho_i(z)
}\over{ (\psi(P_i)^2, v_i ) }}$. The Lemma follows.  \hfill \qed\medskip

Prop.\ \ref{var:E} now follows from (\ref{chev}) and Lemma \ref{new:id}. 
\hfill \qed\medskip 

\begin{remark} It is natural that when $i\neq j$, the logarithmic primitive of 
$$G(z,P_i)_{z \ll P_i}G(z,P_j)_{z\ll P_j}$$ has degrees $1$ in $v_0$ and
zero in $v_i,v_j$,  because $G(z,P_i)G(z,P_j)$  has expansions 
$$
G(z,P_i)G(z,P_j) = {{(dz_{P_0})^2}\over{(z_{P_0})^2}} + O (1/z_{P_0}) 
\on{\ near \ } z = P_0, 
$$
and
$$
G(z,P_i)G(z,P_j) = O (1/z_{P_k}) \on{\ near \ } z = P_k, k = i,j, 
$$
$z_{P_k}$ being a local coordinate on $X$ at each $P_k$. 
For the same reason, the logarithmic primitive of 
$G(z,P_i)^2_{z \ll P_i}$ has degree $1$ both in $v_0$ and $v_i$. 
\hfill \qed\medskip 
\end{remark}

\subsection{Proof of Prop.\ \ref{koche}}

Let us prove 1). Let us study the behavior of $Q\mapsto S'(m,P_0,Q|z)$. 
We have clearly $S'(m,P_0,\gamma_{A_a}(Q)|z) =  S'(m,P_0,Q|z)$.
Moreover, it follows from the identities 
$$
{{\pa_a\Theta}\over{\Theta}}(gP_0 - \gamma_{A_c}(Q) - \Delta)  
= {{\pa_a\Theta}\over{\Theta}}(gP_0 - Q - \Delta) + \delta_{ac}, \
G(z, \gamma_{B_c}(Q) )  = G(z, Q) + \omega_c(z)
$$
and
$$   
{{\pa^2_{ab}\Theta}\over{\Theta}} (gP_0 - \gamma_{A_c}(Q) - \Delta)   
= {{\pa^2_{ab}\Theta}\over{\Theta}} (gP_0 - Q - \Delta)   
+ \delta_{ac} {{\pa_a\Theta}\over{\Theta}} (gP_0 - Q - \Delta)   
+ \delta_{bc} {{\pa_a\Theta}\over{\Theta}} (gP_0 - Q - \Delta)  
+ \delta_{ac} \delta_{bc}
$$
that we have 
$S'(m,P_0,\gamma_{B_a}(Q) | z )  = S'(m,P_0,Q | z )$. 

The possible poles of $Q \mapsto S'(m,P_0,Q|z)$ are at $Q = z$ and 
$Q = P_0$. Let us study the pole at $Q = z$.
Fix a local coordinate on $X$. Let us write 
$$
G(z,w) = {{dz}\over{z-w}} + \al(z,w) dz, \quad    
\bar G(z,w) = {{dz}\over{z-w}} + \bar\al(z,w) dz, 
$$
with 
$$
(\bar\al - \al)(z,w) = \sum_a \omega_a(z) r_a(w) . 
$$
Then if $\omega(z) = f(z)dz$ is a one-form, we have  
$$
(\bar D \omega)(z) = [f'(z) - 2\bar\al(z,z) f(z) ](dz)^2  , 
$$
therefore 
$$
\bar D_z G(z,w) = [ - {{1}\over{ (z-w)^2 }} 
- 2 {{ \bar\al(z,z) }\over{ z - w}} + O(1)](dz)^2 . 
$$
On the other hand, since 
$$
G(z,w)^2 = [{1\over{ (z-w)^2 }} + 2{{\al(z,z)}\over{z - w}} + O(1) ] (dz)^2,  
$$
we find that 
$$
{1\over 2}\bar D_z G(z,w) +  {1\over 2}G(z,w)^2 
+ \sum_a r_a(w) \omega_a(z)  G(z,w) 
$$
has no poles when $w$ tends to $z$. It follows that 
$Q\mapsto S'(m,P_0,Q|z)$ has no pole at $Q = z$.  

On the other hand, the function $Q \mapsto \Theta(gP_0 - Q - \Delta)$
has a zero of order $g$ at $P_0$; the function $Q\mapsto G(z,Q)$
vanishes  at $P_i$, therefore the pole of function $Q\mapsto 
S'(m,P_0,Q|z)$ at  $P_0$ is of order $\leq g$. Since $P_0$ is not a Weierstrass
point of $X$,  this implies that $Q\mapsto S'(m,P_0,Q|z)$ is a constant
function.  This proves the first part of 1). The expansion of
$S'(m,P_0|z)$ near $P_0$ follows from the expansions $G(z,Q) = -
{{dz_{P_0}}\over{ z_{P_0}}} + O(1)$ and $\bar D(f(z_{P_0}) dz_{P_0})=
[f'(z_{P_0}) +  \{ {{- 2(g-1)}\over{z_{P_0}}}+ O(1) \} f(z_{P_0})
](dz_{P_0})^2$ near $P_0$.

Let us prove 2). Consider the connection $\nabla^{\cF}$ in the case $p = 1$, 
$\La_1\neq 0$ and set $\La = \La_1$. 
Set 
\begin{align*}
& \wt g(m_\infty, P_1,v_1 |\la) = f(m_\infty, P_1,v_1 |\la)  
\Theta(m|-2\la + (g-1)P_0 - \Delta)  
(E(P_0,P_1),v_0 \otimes v_1)^{{\La(\La + 2)}\over{2\kappa}} 
\\ & \Theta(m| gP_0 - P_i - \Delta)^{ - {{\La}\over{\kappa}}} .  
\end{align*}
Then $f(m_\infty, P_1,v_1 |\la)$ is a flat section of $\nabla^{\cF}$ iff 
$\wt g(m_\infty, P_1,v_1 |\la)$ obeys the equation
\begin{align} \label{eq:1}
 & \pa_{[\xi]} \wt g(m_\infty, P_1,v_1 | \la ) =
 \left[ {1\over{8 i \pi\kappa}} \sum_{a,b} \delta\tau_{ab} 
 \pa^2_{\la_a\la_b}  
 + \sum_a  \langle \xi(z) , {1\over {2\kappa}} 
 (\bar D\omega_a)(z)  - {{\La}\over{\kappa}}
 \omega_a(z) G(z,P_1)  \rangle_z \pa_{\la_a}  \right. \\ & \nonumber 
 \left. + \langle \xi(z) , {k \over{2\kappa}}  \bar S(m_\infty,P_0|z) 
 - {\La\over\kappa} S'(m_\infty,P_0|z) 
\rangle_z \right]  
 \wt g(m_\infty, P_1,v_1 | \la ) . 
\end{align}
Let us set now 
$$
\wt h(m_\infty,P_1,v_1 |\la')  = g(m_\infty,P_1,v_1 | {1\over\kappa}
[\la' + (g-1)P_0 - \Delta + \La (P_1 - P_0) ]) 
\al(m,P_0,v_0)^{ - {k\over{2\kappa}}} . 
$$
Then (\ref{eq:1}) is equivalent to 
\begin{equation} \label{eq:2}
\pa_{[\xi]} \wt h(m_\infty,P_1,v_1 |\la')
= 
\left[ {1\over{8i\pi\kappa}} \sum_{a,b} \delta\tau_{ab} \pa^2_{\la'_a
\la'_b} -  {{\La}\over{\kappa}}
\langle \xi , S'(m_\infty,P_0|z)  \rangle_z \right] 
\wt h(m_\infty,P_1,v_1 |\la') . 
\end{equation}
Consider the connection over the trivial bundle $\cM_{g,1^2}$, with fibre 
$\CC[[\la']]$, defined by 
$$
\nabla^{aux \prime \prime} = \pa_{\xi} - 
\left[ {1\over{8i\pi\kappa}} \sum_{a,b} \delta\tau_{ab} \pa^2_{\la'_a
\la'_b} -  {{\La}\over{\kappa}}
\langle \xi , S'(m_\infty,P_0|z)  \rangle_z \right] . 
$$ 
This connection is flat. Therefore, the connection over the trivial bundle over 
$\cM_{g,1^2}$ with fiber $\CC$, defined by 
$$
\nabla^{aux \prime \prime\prime} = \pa_{\xi} + {{\La}\over{\kappa}}
\langle \xi , S'(m_\infty,P_0|z)  \rangle_z , 
$$ 
is also flat. The local existence of $\beta(m,P_0,v_0)$ follows. 
(\ref{local:S'}) implies that  $\beta(m,P_0,v_0)$ is homogeneous in $v_0$ of 
degree $-g$. 

(\ref{eagle}) follows then from adding up (\ref{id:beta}) and   
\begin{align*}
& \pa_{[\xi]} \ln \{ {1\over{\Theta(gP_0 - P_i - \Delta) }} 
(E(P_0,P_i), v_0\otimes v_i) \}
= \langle \xi, - \sum_{a,b} {1\over 2} (\omega_a\omega_b)(z) 
{{\pa^2_{ab}\Theta}\over{\Theta}} (gP_0 - P_i - \Delta)
\\ & + \sum_a {{\pa_a\Theta}\over{\Theta}}(gP_0 - P_i - \Delta) [
{1\over 2} \bar D\omega_a(z) + \omega_a(z) G(z,P_i)] \rangle . 
\end{align*}
\hfill \qed\medskip

\subsection{Proof of Thm. \ref{thm:pts}}

Let us set 
\begin{align*}
& f'(m_\infty, P_i,v_i | \la) = f(m_\infty, P_i,v_i | \la)   
\prod_{i<j} \left( {{E(P_i,P_j)}\over{E(P_0,P_i)E(P_0,P_j)}}, v_0
\right)^{ - {{\La_i\La_j}\over{2\kappa}}} 
\\ & \prod_i \left( E(P_0,P_i), v_0\otimes v_i \right)^{ {{\La_i^2}\over {2\kappa}} } 
\prod_i \{ {1\over{\Theta(gP_0 - P_i - \Delta) }} 
(E(P_0,P_i), v_0\otimes v_i) (\beta(m,P_0),v_0) \}^{ {\La_i \over \kappa}} . 
\end{align*}
Then it follows from Props.\ \ref{prop:var:ratio:E}, \ref{var:E} and 
(\ref{eagle}) that (\ref{flatness:P_i}) is equivalent to
\begin{align} \label{eq:3}
& 2\kappa \pa_{[\xi]} f'(m_\infty, P_i,v_i | \la) = 
\left[ {1\over 2} [\sum_a \omega_a(z) \pa_{\la_a} ]^2 
- \sum_{a,i} \La_i \omega_a(z) G(z,P_i) \pa_{\la_a}  
\right. \\ & \nonumber \left. + \sum_a (D^{(2\la)} \omega_a)(z) \pa_{\la_a}  
- \sum_i \La_i (D_z^{(2\la)} - \bar D_z) G(z,P_i) 
+ k \bar\omega_{2\la}(z) 
\right] f'(m_\infty, P_i,v_i | \la ) . 
\end{align}
Let us set now 
$$
g(m_\infty, P_i,v_i | \la ) = \Theta(m | - 2\la + (g-1) P_0 - \Delta)
f'(m_\infty, P_i,v_i | \la ) ,  
$$
then (\ref{eq:3}) is equivalent to 
\begin{align} \label{eq:4}
& \pa_{[\xi]} g(m_\infty, P_i,v_i | \la ) =  
{1\over{8 i \pi\kappa}} \sum_{a,b} \delta\tau_{ab} \pa^2_{\la_a\la_b} 
g(m_\infty, P_i,v_i | \la ) 
+ {1\over{2\kappa}} \sum_a \langle \bar D\omega_a, \xi\rangle \pa_{\la_a}
g(m_\infty, P_i,v_i | \la ) \\ \nonumber & + {k \over{2\kappa}} 
\langle \xi(z) , \bar S(m_\infty,P_0|z) 
\rangle_z g(m_\infty, P_i,v_i | \la ) 
\\ \nonumber & - \sum_{a,i} {{\La_i}\over{\kappa}}
\langle \xi(z), \omega_a(z) G(z,P_i) \rangle_z
{{\pa_{\la_a} ( g(m_\infty, P_i,v_i | \la ) \Theta(-2\la + (g-1)P_0 - \Delta )^{-1}) }
\over{\Theta(-2\la + (g-1)P_0 - \Delta )^{-1} }}
\\ \nonumber & - \sum_i 
{{\La_i}\over{\kappa}} \langle \xi(z) , (D^{(2\la)}_z - \bar D_z) G(z,P_i)  
\rangle_z g(m_\infty, P_i,v_i | \la ) . 
\end{align}
It follows from Prop.\ \ref{pestalozzi}, 2) that this equation is equivalent to 
\begin{align} \label{eq:for:g:pts}
 & \pa_{[\xi]} g(m_\infty, P_i,v_i | \la ) =
 \\ & \nonumber 
 \left[ {1\over{8 i \pi\kappa}} \sum_{a,b} \delta\tau_{ab} 
 \pa^2_{\la_a\la_b}   + \sum_a  \langle \xi(z) , {1\over {2\kappa}} 
 (\bar D\omega_a)(z)  - \sum_i {{\La_i}\over{\kappa}}
 \omega_a(z) G(z,P_i)  \rangle_z \pa_{\la_a}  
 \right. \\ & \nonumber \left. 
 + {k \over{2\kappa}} \langle \xi(z) , \bar S(m_\infty,P_0|z) 
 \rangle_z \right]  g(m_\infty, P_i,v_i | \la ) . 
\end{align}
Let us set 
$$
h(m_\infty,P_i,v_i |\la')  = g(m_\infty,P_i,v_i | {1\over\kappa}
[\la' + (g-1)P_0 + \Delta + \sum_i \La_i(P_i - P_0) ]) . 
$$
Then it follows from Lemmas \ref{var:int} and \ref{var:int:pts} that 
(\ref{eq:for:g:pts}) is equivalent to 
\begin{equation} \label{eq:6}
 \pa_{[\xi]} h(m_\infty, P_i,v_i | \la ) =
 \left[ {\kappa \over{8 i \pi}} \sum_{a,b} \delta\tau_{ab} 
 \pa^2_{\la_a\la_b}  + {k \over{2\kappa}} \langle \xi(z) , \bar S(m_\infty,P_0|z) 
 \rangle_z \right]  h(m_\infty, P_i,v_i | \la ) . 
\end{equation}
Let us set now 
$$
\wt h(m_\infty,P_i,v_i |\la') =  h(m_\infty,P_i,v_i |\la') 
\al(m,P_0,v_0)^{ - {k\over{2\kappa}}} , 
$$
then (\ref{eq:6}) is equivalent to 
$$ 
 \pa_{[\xi]} h'(m_\infty, P_i,v_i | \la ) =
  {\kappa \over{8 i \pi}} \sum_{a,b} \delta\tau_{ab} 
 \pa^2_{\la_a\la_b}  h'(m_\infty, P_i,v_i | \la ) . 
$$
Thm.\ \ref{thm:pts} follows. \hfill \qed\medskip

\newpage

\section{Flat sections of $\cF_{\ell,(\La_i)}$ in the general case}
\label{sect:general}

\subsection{Integral formulas for flat sections of $\cB_{\ell,(\La_i)}$}

Let $\VV$ be an arbitrary highest weight module of level $k$ over $\G$. 
Let $(\La_i)_{i = 1, \ldots, n}$ be weights of $\bar\G = \SL_2$ and let $p$
be an integer $\geq 0$. Let $(m_\infty,P_i,v_i,t_j,w_j) \mapsto 
\psi(m_\infty,P_i,v_i,t_j,w_j)$ be a flat section of the connection 
$(\cB_{\ell,(\La_i),-2, \ldots,-2},\nabla^{\cB_{\ell,(\La_i),-2, \ldots,-2}})$
($-2$ repeated $p$ times) 
over $\cM_{g,1^\infty, (n+p) \cdot 1^2}$. Recall that the fiber  
of $\cB_{\ell,(\La_i),-2, \ldots,-2}$ at $(m_\infty,P_i,v_i,t_j,w_j)$ is given by 
$$
[(\VV \otimes (\otimes_i V_{-\La_i}) \otimes V_{2\vpi}^{\otimes p} )^*]^{\bar\G(R)} , 
$$
where the action of $\bar\G(R)$ is dual to 
$\pi_{\VV} \otimes (\otimes_i \pi_{V_{-\La_i}}^{(P_i)})
\otimes (\otimes_j \pi_{V_{2\vpi}}^{(t_j)})$. Here $\vpi$ denote the fundamental
weight of $\SL_2$.  
Since the Casimir element of $\bar\G$ acts by zero on $V_{2\vpi}$, 
$\psi(m_\infty,P_i,v_i,t_j,w_j)$ is actually independent on the 
tangent vector $w_j$ at $t_j$. We therefore denote it simply 
$\psi(m_\infty,P_i,v_i,t_j)$.
We have  
\begin{equation} \label{psi:KZB}
\pa_{[\xi]} \psi(m_\infty,P_i,v_i,t_j) = \psi(m_\infty,P_i,v_i,t_j)
\circ T_{\bar\omega}[\xi]^{(0)} . 
\end{equation}

If $\rho$ is a rational function on $X$, regular at points $P_i$, we
will identify  it with its image in $\CC((t)) \oplus (\oplus_i \cO_i)$. 
In particular, for $t$ a point of $X - \{P_0,P_i\}$, $z\mapsto \bar
G(t,z)$  is a rational function in $z$, regular on $X - \{P_0,t\}$,
which we will denote $\bar G(t,\cdot)$, and identify with its image in
$\CC((t)) \oplus (\oplus_i \cO_i)$. Therefore $f[\bar G(t,\cdot)]$ belongs to 
$\G_{reg\ at\ P_i}$ (see sect.\ \ref{bunuel}). 

For $(m_\infty,P_i,v_i,t_j)$ in $\cM_{g,1^\infty,n\cdot 1^2,p}$ and $u$ in 
$\VV \otimes (\otimes_{i=1}^n V_{-\La_i})$, let us 
set 
$$
\langle \wt\psi(m_\infty,P_i,v_i,t_j) , u \rangle  = \langle
\psi(m_\infty,P_i,v_i,t_j) , \left( \prod_{j=1}^p \pi_{\VV \otimes (\otimes_i V_{-\La_i}) }
(f[ \bar G(t_j, \cdot )]) u \right) \otimes v_{2\vpi}^{\otimes p}\rangle . 
$$  

$\wt\psi(m_\infty,P_i,v_i,t_j)$ is a function on $\cM_{g,1^\infty,n\cdot
1^2,p}$, with values in  $\left( \VV \otimes (\otimes_{i=1}^n
V_{-\La_i}) \right)^*$.

\begin{prop} \label{freud}

1) There exists a function $\nu(m_\infty,P_i,v_i,t_j)$  
of $\cM_{g,1^\infty,n \cdot 1^2,p}$ with
values in $[\bar\G(R) \otimes \VV \otimes (\otimes_i V_{-\La_i} ) ]^*$, 
such that we have for $x$ in 
$\bar\G(R)$ and $u$ in $\VV \otimes (\otimes_i V_{-\La_i})$, 
\begin{equation} \label{quasi:invce}
\langle \wt\psi(m_\infty,P_i,v_i,t_j) 
, \pi_{\VV \otimes (\otimes_i V_{-\La_i})}(x)
(u) \rangle = d_t (\langle \nu(m_\infty,P_i,v_i,t_j) , x\otimes u\rangle ) ,  
\end{equation} 
where $d_t$ denotes
the partial differential with respect to variables $t_j$. 
Moreover, we have 
\begin{equation} \label{quasi:flatness}
\pa_{[\xi]}\wt\psi(m_\infty,P_i,v_i,t_j)   = 
\wt\psi(m_\infty,P_i,v_i,t_j)  \circ T_{\bar\omega}[\xi]^{(0)} 
\end{equation} 
for $\xi$ in $\CC((t)) {\pa\over{\pa t}}$. 

2) Assume that for $(m_\infty,P_i)$ fixed, the function defined by 
$(t_i)\mapsto \langle \psi(m_\infty,P_i,v_i,t_j), u \rangle$ extends to 
the complement of diagonals of $(X - \{P_0,P_i\})^p$ to a function with 
monodromy properties dependent on $(m_\infty,P_i)$, in other words, to 
the section of a bundle  $\cL(m_\infty,P_i)$. Then 
if $C$ belongs to $H_p((X - \{P_0,P_i\})^p - \on{diagonals}, \cL(m_\infty,P_i))$, 
$(m_\infty,P_i,v_i)\mapsto \int_C
\wt\psi(m_\infty,P_i,v_i,t_j)$ is a flat section of the connection 
$(\cB_{\ell,(\La_i)}, \nabla^{\cB_{\ell,(\La_i)}})$. 
\end{prop}

{\em Proof.} Let us prove (\ref{quasi:invce}). Define 
$\nu(m_\infty,P_i,v_i,t_j)$ as follows: for $x$ in $(\CC f \oplus \CC h)\otimes R$, 
$\langle \nu(m_\infty,P_i,v_i,t_j), x\otimes u\rangle = 0$; 
for $x = e[r], r\in R$,  
\begin{align*}
& \langle \nu(m_\infty,P_i,v_i,t_j), x\otimes u\rangle 
\\ & = - \kappa 
\sum_{j=1}^p r(t_j) \langle \psi(m_\infty,P_i,v_i,t_j), 
\{\prod_{j'\neq j} \pi_{\VV \otimes (\otimes_j V_{-\La_j}) }
(f[\bar G(t_{j'}, \cdot)])(u) \} \otimes v_{2\vpi}^{\otimes p}\rangle .
\end{align*}
 
Then we have, for $x = f[r]$, $r\in R$,   
\begin{align*}
& \langle \wt\psi(m_\infty,P_i,v_i,t_j) 
, \pi_{\VV \otimes (\otimes_i V_{-\La_i})}(x) (u) \rangle 
\\ & = \langle \psi(m_\infty,P_i,v_i,t_j) , 
\{ \prod_j \pi_{\VV \otimes (\otimes_i V_{-\La_i})} 
(f[\bar G(t_j,\cdot)])
\pi_{\VV \otimes (\otimes_i V_{-\La_i})}(f[r])(u)\} 
\otimes v_{2\vpi}^{\otimes p} \rangle 
\\ & = 
\langle \psi(m_\infty,P_i,v_i,t_j) , 
\{ \pi_{\VV \otimes (\otimes_i V_{-\La_i})}(f[r])
\prod_j \pi_{\VV \otimes (\otimes_i V_{-\La_i})} 
(f[\bar G(t_j,\cdot)])
(u) \}  \otimes v_{2\vpi}^{\otimes p} \rangle 
\end{align*}
because $f[r]$ commutes with all $f[\varrho]$, $\varrho$ in 
$\CC((t))\oplus (\oplus_i\cO_{P_i})$. 
Then 
\begin{align*}
& \langle \psi(m_\infty,P_i,v_i,t_j) , 
\{ \pi_{\VV \otimes (\otimes_i V_{-\La_i})}(f[r])
\prod_j \pi_{\VV \otimes (\otimes_i V_{-\La_i})} 
(f[\bar G(t_j,\cdot)])
(u) \}
\otimes v_{2\vpi}^{\otimes p} 
\\ & = 
- \sum_j 
\langle \psi(m_\infty,P_i,v_i,t_j) , 
\{ \prod_{j'} \pi_{\VV \otimes (\otimes_i V_{-\La_i})} 
(f[\bar G(t_{j'},\cdot)])
(u) \}  \otimes \pi_{V_{2\vpi}}(f[r])^{(j)} v_{2\vpi}^{\otimes p}
\end{align*}
by the $\bar\G(R)$-invariance of $\psi(m_\infty,P_i,v_i,t_j)$. 
Since we have $f[t^j]v_{2\vpi} = 0$ for all $j\geq 0$, each
term is the right side of this equality vanishes. Therefore
$$\langle
\wt\psi(m_\infty,P_i,v_i,t_j)  , \pi_{\VV \otimes (\otimes_i
V_{-\La_i})}(x) (u) \rangle$$ is zero. Therefore (\ref{quasi:invce})
is satisfied for $x = f[r]$.  

For $x = h[r]$, $r\in R$, we have 
\begin{align} \label{east}
& \langle \wt\psi(m_\infty,P_i,v_i,t_j) , 
\pi_{\VV \otimes (\otimes_i V_{-\La_i})}(x) (u) \rangle 
\\ & \nonumber 
= \langle \psi(m_\infty,P_i,v_i,t_j) , 
\{ \prod_j \pi_{\VV \otimes (\otimes_i V_{-\La_i})} 
(f[\bar G(t_j,\cdot)])
\pi_{\VV \otimes (\otimes_i V_{-\La_i})}(h[r])(u)\} 
\otimes v_{2\vpi}^{\otimes p} \rangle 
\\ & \nonumber = 
\langle \psi(m_\infty,P_i,v_i,t_j) , 
\{ \pi_{\VV \otimes (\otimes_i V_{-\La_i})}(h[r])
\prod_{j} \pi_{\VV \otimes (\otimes_i V_{-\La_i})} 
(f[\bar G(t_{j},\cdot)])
(u)\} 
\otimes v_{2\vpi}^{\otimes p} \rangle 
\\ & \nonumber 
+ 2 \sum_j 
\langle \psi(m_\infty,P_i,v_i,t_j) , 
\{ \pi_{\VV \otimes (\otimes_i V_{-\La_i})}(f[r\bar G(t_j,\cdot)])
\prod_{j'\neq j} \pi_{\VV \otimes (\otimes_i V_{-\La_i})} 
(f[\bar G(t_{j'},\cdot)])
(u)\} 
\otimes v_{2\vpi}^{\otimes p} \rangle .  
\end{align} 
Since $\psi(m_\infty,P_i,v_i,t_j)$ is $\bar\G(R)$-invariant, 
\begin{align} \label{west}
& \langle \psi(m_\infty,P_i,v_i,t_j) , 
\{ \pi_{\VV \otimes (\otimes_i V_{-\La_i})}(h[r])
\prod_{j} \pi_{\VV \otimes (\otimes_i V_{-\La_i})} 
(f[\bar G(t_{j},\cdot)])
(u)\} 
\otimes v_{2\vpi}^{\otimes p} \rangle 
\\ & \nonumber = 
- \sum_{j'} r(t_{j'}) 
\langle \psi(m_\infty,P_i,v_i,t_j) , 
\{ \pi_{\VV \otimes (\otimes_i V_{-\La_i})}(h[r])
\prod_{j} \pi_{\VV \otimes (\otimes_i V_{-\La_i})} 
(f[\bar G(t_{j},\cdot)])
(u)\} 
\otimes h^{(j')}v_{2\vpi}^{\otimes p} \rangle 
\\ & \nonumber = 
- 2\sum_{j'} r(t_{j'}) 
\langle \psi(m_\infty,P_i,v_i,t_j) , 
\{ \pi_{\VV \otimes (\otimes_i V_{-\La_i})}(h[r])
\prod_{j} \pi_{\VV \otimes (\otimes_i V_{-\La_i})} 
(f[\bar G(t_{j},\cdot)]) (u)\} 
\otimes v_{2\vpi}^{\otimes p} \rangle 
\end{align}
because $hv_{2\vpi} = 2v_{2\vpi}$. 

On the other hand, we have for any $u'$ in $\VV \otimes(\otimes_i V_{-\La_i})$
and $t$ in $X - \{P_0,P_i\}$, 
\begin{align} \label{wien}
& \langle \psi(p_\infty,P_i,v_i,t_j) , \{\pi_{\VV \otimes (\otimes_j V_{-\La_j}) }
( f[r \bar G(t,\cdot)]) u'\} \otimes v_{2\vpi}^{\otimes p} \rangle 
\\ \nonumber & 
=  r(t) 
\langle \psi(p_\infty,P_i,v_i,t_j) , \{\pi_{\VV \otimes (\otimes_j V_{-\La_j}) }
( f[\bar G(t,\cdot)]) u'\} \otimes v_{2\vpi}^{\otimes p} \rangle ; 
\end{align}
indeed, the difference of both sides of (\ref{wien}) is 
$$
\langle \psi(p_\infty,P_i,v_i,t_j) , \{\pi_{\VV \otimes (\otimes_j V_{-\La_j}) }
( f[(r - r(t)) \bar G(t,\cdot)]) u'\} \otimes v_{2\vpi}^{\otimes p} \rangle ; 
$$
since $(r - r(t))  \bar G(t,\cdot)$ belongs to $R$, this difference is equal to 
$$
- \sum_{j'} [(r - r(t)) \bar G(t,\cdot)](t_{j'})
\langle \psi(p_\infty,P_i,v_i,t_j) , 
u' \otimes f^{(j')}v_{2\vpi}^{\otimes p} \rangle , 
$$
which is zero because $f v_{2\vpi} = 0$.
(\ref{west}) and (\ref{wien}) imply that the right side of (\ref{east})
vanishes, therefore
$$
\wt\psi(m_\infty,P_i,v_i,t_j) 
, \pi_{\VV \otimes (\otimes_i V_{-\La_i})}(x) (u) \rangle =0
$$
for $x = h[r]$. So (\ref{quasi:invce}) is satisfied for $x = h[r]$.  

Let us assume that $x = e[r]$. We have 
\begin{align} \label{dimanche}
& \langle \wt\psi(m_\infty,P_i,v_i,t_j) , 
\pi_{\VV \otimes (\otimes_i V_{-\La_i})}(x) (u) \rangle 
\\ & \nonumber 
= \langle \psi(m_\infty,P_i,v_i,t_j) , 
\{ \prod_j \pi_{\VV \otimes (\otimes_i V_{-\La_i})} 
(f[\bar G(t_j,\cdot)])
\pi_{\VV \otimes (\otimes_i V_{-\La_i})}(e[r])(u)\} 
\otimes v_{2\vpi}^{\otimes p} \rangle 
\\ & \nonumber = 
\langle \psi(m_\infty,P_i,v_i,t_j) , 
\{ \pi_{\VV \otimes (\otimes_i V_{-\La_i})}(e[r])
\prod_{j} \pi_{\VV \otimes (\otimes_i V_{-\La_i})} 
(f[\bar G(t_{j},\cdot)])
(u)\} 
\otimes v_{2\vpi}^{\otimes p} \rangle 
\\ & \nonumber 
- \sum_j  \langle \psi(m_\infty,P_i,v_i,t_j) , 
\{ \prod_{j'< j} \pi_{\VV \otimes (\otimes_i V_{-\La_i})} 
(f[\bar G(t_{j'},\cdot)])
[\pi_{\VV \otimes (\otimes_i V_{-\La_i})}( h[r\bar G(t_j,\cdot)])
+ k dr(t_j)] \\ & \nonumber 
\prod_{j'> j} \pi_{\VV \otimes (\otimes_i V_{-\La_i})} 
(f[\bar G(t_{j'},\cdot)])
(u)\} 
\otimes v_{2\vpi}^{\otimes p} \rangle .  
\end{align}
The first term of the right side of (\ref{dimanche}) is equal to 
$$
- \sum_j r(t_j) \langle \psi(m_\infty,P_i,v_i,t_j) , 
\{ \prod_{j} \pi_{\VV \otimes (\otimes_i V_{-\La_i})} 
(f[\bar G(t_{j},\cdot)]) (u)\} 
\otimes e^{(j)} v_{2\vpi}^{\otimes p}
\rangle  . 
$$
On the other hand, we have 
\begin{align*}
& \sum_{j=1}^p  
\langle \psi(m_\infty,P_i,v_i,t_j) , 
\{  \prod_{j'< j} \pi_{\VV \otimes (\otimes_i V_{-\La_i})} 
(f[\bar G(t_{j'},\cdot)])
  \pi_{\VV \otimes (\otimes_i V_{-\La_i})}( h[r\bar G(t_j,\cdot)])
\\ &  
\prod_{j'> j}  
\pi_{\VV \otimes (\otimes_i V_{-\La_i})} 
(f[\bar G(t_{j'},\cdot)]) (u)\} 
\otimes v_{2\vpi}^{\otimes p} \rangle 
\\ &  \nonumber 
=  \sum_{j=1}^p  \langle \psi(m_\infty,P_i,v_i,t_j) , 
\{ \pi_{\VV \otimes (\otimes_i V_{-\La_i})}( h[r\bar G(t_j,\cdot)])
\prod_{j'\neq j} \pi_{\VV \otimes (\otimes_i V_{-\La_i})} 
(f[\bar G(t_{j'},\cdot)]) (u)\} 
\otimes v_{2\vpi}^{\otimes p} \rangle 
\\ & 
+ 2 \sum_{j' < j}  \langle \psi(m_\infty,P_i,v_i,t_j) , 
\{ \pi_{\VV \otimes (\otimes_i V_{-\La_i})}
(f[r\bar G(t_j,\cdot)\bar G(t_{j'},\cdot)])
\\ & \prod_{j''\neq j,j'} \pi_{\VV \otimes (\otimes_i V_{-\La_i})} 
(f[\bar G(t_{j''},\cdot)])(u)\}  
\otimes v_{2\vpi}^{\otimes p} \rangle  , 
\end{align*}
which is equal to 
\begin{align*}
& \sum_{j=1}^p  r(t_j) \langle \psi(m_\infty,P_i,v_i,t_j) , 
\{ \pi_{\VV \otimes (\otimes_i V_{-\La_i})}( h[\bar G(t_j,\cdot)])
\prod_{j'\neq j} \pi_{\VV \otimes (\otimes_i V_{-\La_i})} 
(f[\bar G(t_{j'},\cdot)]) (u)\} 
\otimes v_{2\vpi}^{\otimes p} \rangle 
\\ & 
- 2 \sum_{j=1}^p 
\left( \sum_{j'=1}^p [r(t_{j'}) - r(t_j)]\bar G(t_j,t_{j'}) \right) 
\\ & \langle  \psi(m_\infty,P_i,v_i,t_j) , 
\{ \prod_{j''\neq j} \pi_{\VV \otimes (\otimes_i V_{-\La_i})} 
(f[\bar G(t_{j''},\cdot)]) (u)\} 
\otimes v_{2\vpi}^{\otimes p} \rangle 
\\ & 
+ 2 \sum_{j'\neq j}
r(t_j) \bar G(t_{j'},t_j) 
\langle \psi(m_\infty,P_i,v_i,t_j) , 
\{ \prod_{j''\neq j'} \pi_{\VV \otimes (\otimes_i V_{-\La_i})} 
(f[\bar G(t_{j''},\cdot)])(u)\}  
\otimes v_{2\vpi}^{\otimes p} \rangle ,  
\end{align*}
by the identities $r\bar G(t_j,\cdot) \in 
r(t_j)\bar G(t_j,\cdot) + R$, $r\bar G(t_j,\cdot)\bar G(t_{j'},\cdot)
\in r(t_{j'})\bar G(t_j,t_{j'})\bar G(t_{j'},\cdot) 
+ r(t_j)\bar G(t_{j'},t_j)\bar G(t_j,\cdot) + R$, 
and invariance of  $\psi(m_\infty,P_i,v_i,t_j)$. 
Here for $j' = j$, we have $[r(t_{j'}) - r(t_j)]\bar G(t_j,t_{j'}) = - dr(t_j)$. 
Therefore 
\begin{align*}
& \sum_{j=1}^p  
\langle \psi(m_\infty,P_i,v_i,t_j) , 
\{  \prod_{j'< j} \pi_{\VV \otimes (\otimes_i V_{-\La_i})} 
(f[\bar G(t_{j'},\cdot)])
  \pi_{\VV \otimes (\otimes_i V_{-\La_i})}( h[r\bar G(t_j,\cdot)])
\\ &  
\prod_{j'> j}  
\pi_{\VV \otimes (\otimes_i V_{-\La_i})} 
(f[\bar G(t_{j'},\cdot)]) (u)\} 
\otimes v_{2\vpi}^{\otimes p} \rangle 
\\ &  \nonumber 
= 
\sum_{j=1}^p  r(t_j) \langle \psi(m_\infty,P_i,v_i,t_j) , 
\{ \pi_{\VV \otimes (\otimes_i V_{-\La_i})}( h[\bar G(t_j,\cdot)])
\prod_{j'\neq j} \pi_{\VV \otimes (\otimes_i V_{-\La_i})} 
(f[\bar G(t_{j'},\cdot)]) (u)\} 
\otimes v_{2\vpi}^{\otimes p} \rangle 
\\ & 
+ 2\sum_{j=1}^p dr(t_j)
\langle \psi(m_\infty,P_i,v_i,t_j) , 
\{ \prod_{j'\neq j} \pi_{\VV \otimes (\otimes_i V_{-\La_i})} 
(f[\bar G(t_{j'},\cdot)])(u)\}  
\otimes v_{2\vpi}^{\otimes p} \rangle 
\\ & 
+ 2 \sum_{j'\neq j} r(t_j) \bar G(t_j,t_{j'})
\langle \psi(m_\infty,P_i,v_i,t_j) , 
\{ \prod_{j''\neq j} \pi_{\VV \otimes (\otimes_i V_{-\La_i})} 
(f[\bar G(t_{j''},\cdot)])(u)\}  
\otimes v_{2\vpi}^{\otimes p} \rangle 
\end{align*}

It follows that we have 
\begin{align*} 
& \langle \wt\psi(m_\infty,P_i,v_i,t_j) , 
\pi_{\VV \otimes (\otimes_i V_{-\La_i})}(x) (u) \rangle 
\\ & \nonumber 
= - (k+2) \sum_j dr(t_j) \langle \psi(P_i,v_i,t_j), \prod_{j'\neq j}
\pi_{\VV \otimes (\otimes_i V_{-\La_i})} (f[\bar G(t_{j'},\cdot)])(u)\}  
\otimes v_{2\vpi}^{\otimes p} \rangle 
\\ & 
- \sum_j r(t_j) \langle \psi(m_\infty,P_i,v_i,t_j) , 
\{ \prod_{j} \pi_{\VV \otimes (\otimes_i V_{-\La_i})} 
(f[\bar G(t_{j},\cdot)]) (u)\} 
\otimes e^{(j)} v_{2\vpi}^{\otimes p}
\rangle  
\\ &  
- \sum_{j=1}^p  r(t_j) \langle \psi(m_\infty,P_i,v_i,t_j) , 
\{ \pi_{\VV \otimes (\otimes_i V_{-\La_i})}( h[\bar G(t_j,\cdot)])
\\ & \prod_{j'\neq j} \pi_{\VV \otimes (\otimes_i V_{-\La_i})} 
(f[\bar G(t_{j'},\cdot)]) (u)\} 
\otimes v_{2\vpi}^{\otimes p} \rangle 
\\ & 
- 2 \sum_{j'\neq j} r(t_j) \bar G(t_j,t_{j'})
\langle \psi(m_\infty,P_i,v_i,t_j) , 
\{ \prod_{j''\neq j} \pi_{\VV \otimes (\otimes_i V_{-\La_i})} 
(f[\bar G(t_{j''},\cdot)])(u)\}  
\otimes v_{2\vpi}^{\otimes p} \rangle . 
\end{align*}
To prove (\ref{quasi:invce}), it suffices therefore to show that 
for any $j$, we have 
\begin{align} \label{todo}
& - \kappa d_{t_j} 
\left( \langle \psi(m_\infty,P_i,v_i,t_j), 
\{\prod_{j'\neq j} \pi_{\VV \otimes (\otimes_j V_{-\La_j}) }
(f[\bar G(t_{j'}, \cdot)])(u) \} \otimes v_{2\vpi}^{\otimes p}\rangle 
\right) 
\\ & \nonumber = 
- \langle \psi(m_\infty,P_i,v_i,t_j) , 
\{ \prod_{j} \pi_{\VV \otimes (\otimes_i V_{-\La_i})} 
(f[\bar G(t_{j},\cdot)]) (u)\} 
\otimes e^{(j)} v_{2\vpi}^{\otimes p}
\rangle  
\\ &  \nonumber 
- \langle \psi(m_\infty,P_i,v_i,t_j) , 
\{ \pi_{\VV \otimes (\otimes_i V_{-\La_i})}( h[\bar G(t_j,\cdot)])
\prod_{j'\neq j} \pi_{\VV \otimes (\otimes_i V_{-\La_i})} 
(f[\bar G(t_{j'},\cdot)]) (u)\} 
\otimes v_{2\vpi}^{\otimes p} \rangle 
\\ & \nonumber 
- 2 \sum_{j'\neq j} \bar G(t_j,t_{j'})
\langle \psi(m_\infty,P_i,v_i,t_j) , 
\{ \prod_{j''\neq j} \pi_{\VV \otimes (\otimes_i V_{-\La_i})} 
(f[\bar G(t_{j''},\cdot)])(u)\}  
\otimes v_{2\vpi}^{\otimes p} \rangle . 
\end{align}

Recall that $\psi$ satisfies (\ref{psi:KZB}). Let in this equation,
$\xi$ be a vector field on $X$, regular outside $P_0$ and vanishing to
second order at the $P_i$ and $t_{j'},j'\neq j$. Let $(r^{\prime\prime i})$ 
and $(\omega^{\prime\prime i})$  be bases of $R_{(P_s,t_{j'},j'\neq j)}$ and
$\Omega_{(P_s,t_{j'},j'\neq j)}$,  and let $(r''_i)$ and $(\omega''_i)$ be
such that $(r^{\prime\prime i},r''_i)$ and $(\omega''_i, 
\omega^{\prime\prime i})$ are
dual bases of $\CC((t))$ and $\CC((t))dt$ (see the proof of Lemma \ref{secu}). 
Let us set 
$T'[\xi] = {1\over{2\kappa}} \sum_{\al} (x^{\al}[\xi\omega^{\prime\prime
i}]x^{\al}[r''_i] +  x^{\al}[r^{\prime\prime i}]x^{\al}[\xi \omega''_i])$. Since
$\xi$ is regular outside $P_0$, $T'[\xi] = T_{\bar\omega}[\xi]$ (see
the proof of Lemma \ref{secu}). Then the invariance of $\psi(m_\infty,P_i,v_i,t_j)$ 
implies that 
\begin{align*}
& \psi(m_\infty,P_i,v_i,t_j) \circ T'[\xi]^{(0)}
= - {1\over{2\kappa}} \sum_{i,s,\al} 
\psi \circ x^{\al}[r''_i]^{(0)} x^{\al}[\xi\omega^{\prime\prime i}]^{(P_s)}
+ \psi \circ x^{\al}[\xi \omega''_i]^{(0)}
x^{\al}[r^{\prime\prime i}]^{(P_s)}
\\ & - {1\over{2\kappa}} \sum_{i,s',\al} 
\psi \circ x^{\al}[r''_i]^{(0)} x^{\al}[\xi\omega^{\prime\prime i}]^{(t_{s'})} 
+ \psi \circ x^{\al}[\xi \omega''_i]^{(0)} x^{\al}[r^{\prime\prime i}]^{(t_{s'})} , 
\end{align*} 
where we set $\psi = \psi(m_\infty,P_i,v_i,t_j)$. All contributions associated to 
the $P_s$ or to the $t_{j'},j'\neq j$, are zero. Therefore 
\begin{align*}
&\psi(m_\infty,P_i,v_i,t_j) \circ T'[\xi]^{(0)}  \\ & 
= - {1\over{2\kappa}} \sum_{i,\al} 
(\xi\omega^{\prime\prime i})(t_j)
\psi \circ x^{\al}[r''_i]^{(0)} (x^{\al})^{(t_j)} 
+ r^{\prime\prime i}(t_j)
\psi \circ x^{\al}[\xi \omega''_i]^{(0)} (x^{\al})^(t_j) . 
\end{align*}
Let us set $\Gamma(z,w)_{w\ll z} 
= \sum_i \omega^{\prime\prime i}(z) r''_i(w)$. $\Gamma(z,w)_{w\ll z}$
is the expansion of a rational form $\Gamma(z,w)$.
The function $\xi(t_j)\Gamma(t_j,\cdot) - \xi(\cdot)\Gamma(\cdot,P_j)$ 
vanishes at all $P_i$ and $t_{j'},j'\neq j$. It follows that 
\begin{align*}
& - {1\over{2\kappa}} \sum_{i,\al} 
(\xi\omega^{\prime\prime i})(t_j)
\psi \circ x^{\al}[r''_i]^{(0)} (x^{\al})^{(t_j)} 
+ r^{\prime\prime i}(t_j)
\psi \circ x^{\al}[\xi \omega''_i]^{(0)} (x^{\al})^(t_j) 
\\ & = 
- {1\over{2\kappa}} \sum_{\al} 
\psi \circ 
\pi_{\VV\otimes (\otimes_i V_{-\La_i}) \otimes (\otimes_{j'\neq j} V_{2\vpi})}
(x^\al[\xi(t_j)\Gamma(t_j,\cdot) - \xi(\cdot)\Gamma(\cdot,P_j)]) 
(x^{\al})^{(t_j)} . 
\end{align*}
Now we have 
$$
2\xi(t_j)\bar G(t_j,\cdot) - 
\{\xi(t_j)\Gamma(t_j,\cdot) - \xi(\cdot)\Gamma(\cdot,P_j) \}
\in R ; 
$$
let us denote by $\dot\delta$ the value at $t_j$ of this difference. 
It follows that  
\begin{align*}
& - {1\over{2\kappa}} \sum_{\al} 
\psi \circ 
\pi_{\VV\otimes (\otimes_i V_{-\La_i}) \otimes (\otimes_{j'\neq j} V_{2\vpi})}
\{ (x^\al[\xi(t_j)\Gamma(t_j,\cdot) - \xi(\cdot)\Gamma(\cdot,P_j)]) 
(x^{\al})^{(t_j)} \}
\\ & 
= - {1\over{2\kappa}} \sum_{\al} 
\psi \circ 
\pi_{\VV\otimes (\otimes_i V_{-\La_i}) \otimes (\otimes_{j'\neq j} V_{2\vpi})}
\{ (x^\al[2\xi(t_j)\bar G(t_j,\cdot)]) (x^{\al})^{(t_j)} \}
- {1\over{2\kappa}} \dot\delta  
\psi \circ \sum_{\al} \{(x^{(\al)})^{(t_j)} \}^2.  
\end{align*}
Since the Casimir element of $\bar\G$ vanishes on $V_{2\vpi}$, this 
implies that that for any $u'$ in $\VV \otimes 
(\otimes_i V_{-\La_i})$ , we have 
\begin{align*}
& - \kappa d_{t_j} 
\left( \langle \psi(m_\infty,P_i,v_i,t_j), 
u' \otimes v_{2\vpi}^{\otimes p}\rangle 
\right) 
\\ & = 
- \langle \psi(m_\infty,P_i,v_i,t_j) , 
\{ \pi_{\VV \otimes (\otimes_i V_{-\La_i})} 
(f[\bar G(t_{j},\cdot)]) (u')\} 
\otimes e^{(j)} v_{2\vpi}^{\otimes p}
\rangle  
\\ &  
- \langle \psi(m_\infty,P_i,v_i,t_j) , 
\{ \pi_{\VV \otimes (\otimes_i V_{-\La_i})}
( h[\bar G(t_j,\cdot)]) (u')\} 
\otimes v_{2\vpi}^{\otimes p} \rangle 
\\ & 
- 2 \sum_{j'\neq j} \bar G(t_j,t_{j'})
\langle \psi(m_\infty,P_i,v_i,t_j) , 
u' \otimes v_{2\vpi}^{\otimes p} \rangle . 
\end{align*}
(\ref{todo}) follows. As we have seen, this implies 
(\ref{quasi:invce}) in the case $x  = e[r]$ and therefore 
ends the proof of (\ref{quasi:invce}).

Let us now prove (\ref{quasi:flatness}). Since $\bar G(t_j,\cdot)$ is a
holomorphic differential in $t_j$, $(t_j)_j\mapsto \langle 
\wt\psi(m_\infty,P_i,v_i,t_j) , u \rangle$ is a section of $H^0(
(X _ \{P_0,P_i\})^p - \on{diagonals} , \Omega^{top} \otimes \cL(m_\infty,P_0,P_i))$, 
so that $\int_C \wt\psi(m_\infty,P_i,v_i,t_j)$ makes sense. 
We have for $u$ in $\VV \otimes
(\otimes_i V_{-\La_i})$ and $\eps^2 = 0$, 
\begin{align*}
& \eps \left( \langle \pa_{[\xi]}\wt\psi - \wt\psi \circ T_{\bar\omega}[\xi]^{(0)} , 
u \rangle  \right) 
\\ & = 
\langle \psi[(1+\eps \xi)R\subset \CC((t)), \chi_{P_i}\circ (1 - \eps\xi), 
\chi_{t_j}\circ (1 - \eps\xi)]  ,
\\ & 
\prod_i [\pi_{\VV} \otimes (\otimes_i \pi_{V_{-\La_i}} \circ (1 - \eps\xi) )](
f[(1 + \eps\xi)\bar G(t_j, \cdot)]) (u \otimes v_{2\vpi}^{\otimes p}) \rangle 
\\ & 
- \langle \psi [R\subset \CC((t)), \chi_{P_i}, \chi_{t_j}]  ,
\prod_i [\pi_{\VV} \otimes (\otimes_i \pi_{V_{-\La_i}})](
f[\bar G(t_j, \cdot)]) ((1 + \eps T_{\bar\omega}[\xi]^{(0)}) 
(u) \otimes v_{2\vpi}^{\otimes p}) \rangle 
\end{align*}
because $(1 + \eps\xi)\bar G(t_j, \cdot)$ is a function on $\Spec((1 + \eps\xi)(R))$
with a simple pole at the image of $t_j$. (\ref{vf}) implies that this is equal to zero. 

Let us prove 2). The formula for $\nu(m_\infty,P_i,v_i,t_j)$ shows that
for  $x$ and $u$ fixed, $\langle \nu(m_\infty,P_i,v_i,t_j) , x\otimes
u\rangle$ is a section of $\cL(m_\infty,P_i)$, so that $\int_C
\wt\psi(m_\infty,P_i,v_i,t_j)$ is $\bar\G(R)$-invariant. It follows then from 
(\ref{quasi:invce}) that it is a flat section of $\nabla^{\cB_{\ell,(\La_i)}}$. 
\hfill \qed\medskip

\begin{prop} \label{adler}
Assume that $\psi(m_\infty,P_i,v_i,t_j | \la)$ is a flat section of 
$(\wt\cB_{\ell,(\La_i), -2,\ldots,-2},$ $\wt\nabla^{\ell,(\La_i), -2,\ldots,-2})$
($-2$ repeated $p$ times); 
this means that it is a $\G_\la^{out}$-invariant form 
on $\VV \otimes(\otimes_i V_{-\La_i}) \otimes V_{2\vpi}^{\otimes p}$, defined on an open 
subset of $\cM_{g,1^\infty, n \cdot 1^2, p\cdot 1} \times \CC^g$, and is
a solution  of (\ref{nagila1}), (\ref{nagila2}). 
Assume also that for $(m_\infty,P_i)$
fixed, there exists a bundle $\cL(m_\infty,P_i) $ over 
$(X - \{P_0,P_i\})^p -\on{diagonals}$ such that 
$(t_j) \mapsto \langle \psi(m_\infty,P_i,v_i,t_j | \la) , \wt u \rangle$ is a 
section of $\cL(m_\infty,P_i)  \otimes \cL_{-2\la}^{\boxtimes p}$, for any $\wt u$
in $\VV \otimes(\otimes_i V_{-\La_i}) \otimes V_{2\vpi}^{\otimes p}$. 
Let us define $\wt\psi(m_\infty,P_i,v_i,t_j | \la)$ by  
$$
\langle \wt\psi(m_\infty,P_i,v_i,t_j | \la) , u\rangle = 
\langle \psi(m_\infty,P_i,v_i,t_j | \la) , \left( \prod_{j=1}^p
\pi_{\VV \otimes (\otimes_i V_{-\La_i})}(f[G_{2\la}(t_i,\cdot)]) u \right) 
\otimes v_{2\vpi}^{\otimes p} \rangle  
$$
for $u$ in $\VV \otimes(\otimes_i V_{-\La_i})$

Then if $C$ belongs to $H_p((X - \{P_0,P_i\})^p - \on{diagonals}, \cL(m_\infty,P_i))$, 
set 
$$
I_C(m_\infty,P_i,v_i|\la) = \int_C \wt\psi(m_\infty,P_i,v_i,t_j|\la).  
$$
Then $I_C(m_\infty,P_i,v_i|\la)$ is a $\G^{out}_\la$-invariant form on 
$\VV \otimes (\otimes_i V_{-\La_i})$, defined
on an open subset of  $\cM_{g,1^\infty, n \cdot 1^2} \times \CC^g$, 
satisfying equations (\ref{nagila1}) and (\ref{nagila2}), that is a flat section of 
$(\wt\cB_{\ell,(\La_i)}, \wt\nabla^{\ell,(\La_i)})$.  
\end{prop}

The proof is similar to that of Prop.\ \ref{freud}. 

\begin{remark} It seems that the generalization of 
Props.\ \ref{freud} and \ref{adler} to a general Lie algebra
$\bar\G$ requires using the Wakimoto modules, see \cite{FF}. 
\end{remark}

\begin{remark} In the elliptic case $g=1$, the cycles $C$ are constructed in 
\cite{Felder:Silvotti}. In that case, $\varphi$ satisfies theta-like
conditions around the $a$- and $b$-cycle. The cycles $C$ consist of
points running on  8-shaped contours and correspond to singular vectors
in products of Verma and adjoint modules for the quantum group
$U_q\SL_2$, $q = e^{{{2i\pi}\over{\kappa}}}$.  
\end{remark}

\subsection{Flat sections of $\cF_{\ell,(\La_i)}$}

\subsubsection{The operators $\Phi_{(\La_i)}(w_1,\ldots,w_p|z_1,\ldots,z_p)$}

For $z$ in $X$, define $\wt h(z)$ as the operator acting on functions
defined on $\CC^g$ by 
$$
(\wt h(z)\Psi)(\la) = \sum_{a=1}^g \omega_a(z) \pa_{\la_a}\Psi(\la). 
$$
If $(w_i)_{i = 1, \ldots,p}$ and $(z_i)_{i = 1, \ldots,p}$ are $2p$
points of $X$, and $(\La_i)_{i= 1, \ldots,q}$ are complex numbers, 
define inductively $\Phi_{(\La_i)}(w_1,\ldots,w_p|z_1,\ldots,z_p)$
as the following operator acting on functions
defined on $\CC^g$: for $p=0$, $\Phi_{(\La_i)} = id$, and 
\begin{align*}
& \Phi_{(\La_i)} (w_1,\ldots,w_p|z_1,\ldots,z_p)
= \sum_{i=1}^p \{ - G_{2\la}(w_1,z_i)_{z_i\ll w_1}
[\wt h(z_i) - \sum_{j=1}^q \La_j G(z_i,P_j)_{P_j\ll z_i}]
\\ &   + k d_{z_i}G_{2\la}(w_1,z_i)_{z_i\ll w_1}
\} \circ 
\Phi_{(\La_i)} (w_2,\ldots,w_p|z_1,\ldots,\check z_i, \ldots, z_p)
\end{align*}

One can show by induction on $p$ that the operator $\Phi_{(\La_i)}
(w_1,\ldots,w_p|z_1,\ldots,z_p)$  is symmetric in both groups of
variables $(z_j)_{j = 1,\ldots,p}$ and  $(w_j)_{j = 1,\ldots,p}$.  The
operators $\wt h(z)$ and $\Phi_{(\La_i)}(w_1,\ldots,w_p|z_1,\ldots,z_p)$
enjoy the following properties.

\begin{lemma} \label{gardel} 
Let $\psi_\la$ be a family of $\G_\la^{out}$-invariant forms on $\VV\otimes 
(\otimes_i V_{-\La_i})$, such that 
$$
\pa_{\la_a} \psi_\la = \psi_\la \circ 
\pi_{\VV\otimes (\otimes_i V_{-\La_i})}(h[r_a]) . 
$$
Then 
$$
\langle \psi_\la, \pi_{\VV}(h(z)) (v\otimes (\otimes_i v_{-\La_i})) \rangle 
= 
\left( \wt h(z) - \sum_{i=1}^q \La_i G(z,P_i)_{P_i \ll z}\right) 
(\langle \psi_\la, v\otimes (\otimes_i v_{\La_i}) \rangle ) . 
$$
We have also 
\begin{align*}
& \langle \psi_\la, \pi_{\VV}(\prod_{j=1}^p f(w_j)\prod_{j=1}^p e(z_j)) 
(v\otimes (\otimes_i v_{-\La_i})) \rangle 
\\ & =  \Phi_{(\La_i)} (w_1,\ldots,w_p|z_1,\ldots,z_p) 
(\langle \psi_\la, v\otimes (\otimes_i v_{-\La_i}) \rangle ) . 
\end{align*}
\end{lemma}

{\em Proof.} See \cite{EF}. \hfill \qed\medskip 

If $g(w_1,\ldots,w_p|\la)$ is an analytic function of an 
open subset of $X^p\times\CC^g$, define $\Phi_{(\La_i)}(w_1,\ldots,w_p|
z_1,\ldots,z_p)_{an}g(w_1,\ldots, w_p|\la)$ the same way as 
$\Phi_{(\La_i)}$ is defined, replacing the $G_{2\la}(z,w)_{z\ll w}$
by their meromorphic prolongations.

\begin{prop} \label{fontainebleau}
Let us assume that in Prop.\ \ref{adler}, $p = N = {1\over 2} [k\ell +
\sum_{i}\La_i]$.  Let $(I_C\psi)(m_\infty,P_i,v_i|\la)$ be the flat section of 
$\cB^{\ell,(\La_i)}$ constructed in that
Proposition.  Its image by the map $\corr$ of sect.\ \ref{corr} is given by 
\begin{align} \label{expr:corr}
& [\corr(I_C\psi)](m_\infty,P_i,v_i)(z_1,\ldots,z_{N}|\la) = 
\\ & \nonumber \int_C \Phi_{\La_1,\ldots,\La_n,-2,\ldots,-2}
(t_1\ldots,t_N|z_1,\ldots,z_N)_{an}
\{\corr(\psi)(m_\infty,P_i,v_i,t_1,\ldots,t_N|\la) \}, 
\end{align}
where $\corr(\psi)(m_\infty,P_i,v_i,t_1,\ldots,t_N|\la)$ is expressed by 
the right side of (\ref{sol:P_i}), with $(P_i)_{i = 1, \ldots,n}$  and 
$(\La_i)_{i = 1, \ldots,n}$ replaced by $(P_1,\ldots,P_n, t_1,\ldots,
t_N)$  and 
$(\La_1,\ldots,\La_n,$ $-2,\ldots,-2)$ ($-2$ repeated $N$ times). 
\end{prop}

{\em Proof.} 
For $\rho_1,\ldots,\rho_p$ in $\CC((t))$, we will set 
\begin{align*}
& \Phi_{(\La_i)} [\rho_1,\ldots,\rho_p](z_1,\ldots,z_p) 
\\ & = \res_{w_1 = P_0} \cdots \res_{w_p = P_0}
[\Phi_{(\La_i)} (w_1,\ldots,w_p|z_1,\ldots,z_p) 
\rho_1(w_1) \cdots \rho_p(w_p)] . 
\end{align*}
The second identity of Lemma \ref{gardel} then implies  
that the image of $(I_C\psi)(m_\infty,P_i,v_i|\la)$
by the map $\corr$ of sect.\ \ref{corr} is given by 
\begin{align*} 
& [\corr(I_C\psi)](m_\infty,P_i,v_i)(z_1,\ldots,z_{N}|\la) = 
\\ &  \int_C \Phi_{\La_1,\ldots,\La_n,-2,\ldots,-2}
[G_{2\la}(t_1,\cdot), \ldots, G_{2\la}(t_N,\cdot)] (z_1,\ldots, z_N)
\\ &  \quad  \{\corr(\psi)(m_\infty,P_i,v_i,t_1,\ldots,t_N|\la) \} . 
\end{align*}
It follows then from 
$$
\int_{t'\on{\ around\ }P_0,z} G_{2\la}(t',z) G_{2\la}(t,t') = G_{2\la}(t,z)
$$
that this expression is equal to (\ref{expr:corr}). 
\hfill \qed\medskip

\begin{thm} \label{jahrzeit}
For any $\varphi(m|\la)$ satisfying the heat equation
$$ 
 \pa_{[\xi]} \varphi(m|\la) = 
  {\kappa \over{8 i \pi}} \sum_{a,b} \delta\tau_{ab} 
 \pa^2_{\la_a\la_b}  \varphi(m|\la)  
$$
(we indicate at the end of Thm.\ \ref{thm:no:point} how such functions
$\varphi$ may be  naturally obtained), define
$f_{\varphi}(m_\infty,P_i,v_i,t_1,\ldots,t_N|\la)$ as the right side of
(\ref{sol:P_i}), with $(P_i)_{i = 1, \ldots,n}$  and  $(\La_i)_{i = 1,
\ldots,n}$ replaced by $(P_1,\ldots,P_n, -2,\ldots,-2)$  and 
$(\La_1,\ldots,\La_n,-2,\ldots,-2)$ ($-2$ repeated $N$ times). 

For any cycle $C$ in $H_p((X - \{P_0,P_i\})^p - \on{diagonals},
\cL(m_\infty,P_i))$,  set
\begin{align*}
& f_{C,\varphi}(m_\infty,P_i,v_i)(z_1,\ldots,z_{N}|\la) = 
\\ & \nonumber \int_C \Phi_{\La_1,\ldots,\La_n,-2,\ldots,-2}
(t_1,\ldots,t_N | z_1,\ldots, z_N)_{an}\{ f_{\varphi}(m_\infty,P_i,v_i,t_1,\ldots,t_N|\la) \} . 
\end{align*}
Then $(m_\infty,P_i,v_i)\mapsto f_{C,\varphi}(m_\infty,P_i,v_i)$ defines
a flat section of $(\cF_{\ell,(\La_i)},\nabla^{\cF_{\ell,(\La_i)}})$.
\end{thm}

{\em Proof.} Let $\G_{1-g}$ be the subalgebra of $\G$ equal to 
$$
\G_{1-g} = \left( \bar\HH \otimes \CC[[t]] \right) 
\oplus \left( \bar\N_- \otimes t^{1-g}\CC[[t]] \right) \oplus \CC K.  
$$
Let $\chi$ be the character of $\G_{1-g}$ defined by  $\chi(K) = k,
\chi(f[t^i]) = 0, \chi(h[t^i]) = -\ell \delta_{i0}$.  Let $\WW$ be the
induced $\G$-module $\WW =  U\G \otimes_{U\G_{1-g}} \CC_\chi$.

Assume that $\la$ is generic. In \cite{EF}, we identified the functional
connection $(\cF_{\ell,(\La_i)}, \nabla^{\cF_{\ell,(\La_i)}})$ with the
version of the connection  $(\wt\cB_{\ell,(\La_i)} ,
\wt\nabla^{\ell,(\La_i)})$ in which  $\VV$ is replaced by $\WW$. 

It is clear that Prop.\ \ref{freud} holds when $\VV$ is replaced
by $\WW$. This implies the Theorem. 
\hfill \qed\medskip 

\begin{remark} One can define operators $\Phi_{\La_i,p}(z_1,\ldots,z_q
| w_1,\ldots,w_q)$ generalizing the $\Phi_{\La_i}(z_1,\ldots,z_p
| w_1,\ldots,w_p)$, such that 
\begin{align*}
&\langle \psi_\la, \pi_{\VV} 
\left( \prod_{j=1}^q f(w_j)\prod_{j=1}^p e(z_j) \right) 
[v \otimes (\otimes_i v_{-\La_i})] \rangle 
\\ & =
\Phi_{\La_i,p}(z_1,\ldots,z_q | w_1,\ldots,w_q)\{
\langle \psi_\la, 
\pi_{\VV}\left( \prod_{j=q+1}^{p} e(z_j)\right) 
[v \otimes (\otimes_i v_{-\La_i})] \rangle  \}
\end{align*}

Using these operators, it is then easy to generalize Prop.\ 
\ref{fontainebleau} to the case $p\neq N$. This also provides integral
formulas, expressing flat sections of $\cF_{\ell,(\La_i)}$ in terms of
flat sections of  $\cF_{\ell,(\La_i),-2,\ldots,-2}$ ($-2$ repeated $p$
times). However,  it is not clear that the family of flat sections
obtained in this  way from the flat sections with $N = 0$ is more
general than  that of Thm.\ \ref{jahrzeit}.    
\end{remark}

\begin{remark} \label{rem:integrable}
{\it Integrability conditions.} 
As we explained in Remark \ref{rem:cb}, it should be possible, by imposing
functional conditions on $f_{C,\varphi}(m_\infty,P_i,v_i)(z_1,\ldots,z_{N}|\la)$, 
to find the conditions on $\varphi$ and $C$ for this function to be a 
correlation function of conformal blocks. One probably finds this
way that $\varphi$ is a theta-function, satisfying some vanishing 
conditions.
\end{remark}


\subsubsection{Explicit formula in the case $g = 2,k=1,N=1$}
\label{sect:example}

In this simple case, we find the formula
\begin{align*}
& f(m,P_0,v_0)(z|\la)
= 
\al(m,P_0,v_0)^{1/6} (\beta(m,P_0),v_0)^{-1} \cdot 
\\ & 
\cdot \int_{t\in C}
[ - G_{2\la}(t,z) (\sum_a \omega_a(z)\pa_{\la_a}) + 
d_z G_{2\la}(t,z)]
\{ \varphi(m|3\la + \Delta + P_0 - 2 t)
\Theta(m|-2\la + P_0 - \Delta)^{-1}  
\\ & \quad \Theta(m|2P_0 - t - \Delta)^{2/3}  \}, 
\end{align*}
where $\varphi$ is a solution of (\ref{eq:varphi}).

 \end{document}